\documentclass[a4paper,12pt]{amsart}

\usepackage{amsmath, amsfonts, amsthm, amsfonts, amssymb}
\usepackage[T1]{fontenc}
\usepackage[latin1]{inputenc}

\def\]{\textup{\mbox{]\hspace{-.15em}]}}}
\def\[{\textup{\mbox{[\hspace{-.15em}[}}}

\def\got{\mathfrak}

\newenvironment{pf}
{\medskip\noindent {\it Proof --- \ }}
{\hfill\nobreak $\Box$ \par\bigbreak}

\newcommand{\M}{{\rm M}}
\newcommand{\F}{{\mathbb F}}

\newcommand{\Tr}{{\rm Tr}}
\newcommand{\isomo}{\overset{\sim}{\rightarrow}}
\newcommand{\cal}{\mathcal}

\newcommand{\Hom}{\mathrm{Hom}}

\newcommand{\Gal}{\mathrm{Gal}}

\newcommand{\tr}{\mathrm{tr}}

\newcommand{\Spec}{\mathrm{Spec}}

\newcommand{\GL}{\mathrm{GL}}

\newcommand{\rig}{\mathrm{rig}}

\newcommand{\ps}{\par \smallskip}

\newcommand{\End}{\mathrm{End}}
\newcommand{\N}{\mathbb{N}}
\newcommand{\Z}{\mathbb{Z}}

\newcommand{\Q}{\mathbb{Q}}

\newcommand{\C}{\mathbb{C}}

\newcommand{\rhob}{{\bar\rho}}

\newcommand{\Ker}{\mathrm{ker}}

\newcommand{\kb}{\overline{k}}

\newcommand{\spec}{\mathrm{Spec}}

\newcommand{\OO}{\mathcal{O}}

\newcommand{\D}{\mathrm{D}}

\newcommand{\TS}{{\rm TS}}
\newcommand{\Det}{{\rm D}}
\newcommand{\CH}{{\rm CH}}
\newcommand{\Rad}{{\rm Rad}}
\newcommand{\rad}{{\rm Rad}}

\newcommand{\Rb}{{\overline{\rm R}}}
\newcommand{\Detb}{{\overline{\Det}}}

\newcommand{\CX}{{\cal X}}

\newcommand{\TE}{\widetilde{E}} \newcommand{\plim}{{\underset{\leftarrow}{\lim}}}\newcommand{\An}{{\rm An}} 
\newcommand{\Spf}{{\rm Spf}}
\newcommand{\Dbar}{\overline{\Det}}
\newcommand{\Ext}{{\rm Ext}}

\newcommand{\cCH}{{\cal{C\!H}_d(G)}}
\newcommand{\St}{{\rm St}}
\newcommand{\Aff}{\rm Aff}
\newcommand{\RED}{{\rm Red}}

\newtheorem{thm}[subsection]{Theorem}
\newtheorem{lemma}[subsection]{Lemma}
\newtheorem{remark}[subsection]{Remark}
\newtheorem{cor}[subsection]{Corollary}
\newtheorem{prop}[subsection]{Proposition}
\newtheorem{example}[subsection]{Example}

\newtheorem{exo}[subsection]{Exercise}
\newtheorem{definition}[subsection]{Definition}
\newtheorem{defprop}[subsection]{Definition-Proposition}

\input xy
\xyoption{all}

\title[The $p$-adic analytic space of pseudocharacters of a profinite
group]{The $p$-adic analytic space of pseudocharacters of a profinite group
and pseudorepresentations over arbitrary rings}

\author{Ga\"etan Chenevier}

\begin{document}

\maketitle \medskip {\sc Abstract : } Let\footnote{The author is supported by the C.N.R.S., as well as by the
A.N.R. project ANR-10-BLAN 0114.} $G$ be a profinite
group which is topologically finitely generated\footnote{Actually, we only
assume that for any normal open subgroup $H \subset G$, there are only
finitely many continuous group homomorphisms $H \longrightarrow \Z/p\Z$.},
$p$ a prime number and $d\geq 1$ an integer.  We show that the functor from
rigid analytic spaces over $\Q_p$ to sets, which associates to a rigid space
$Y$ the set of continuous $d$-dimensional pseudocharacters $G
\longrightarrow \OO(Y)$, is representable by a quasi-Stein rigid analytic
space $X$, and we study its general properties.  \par Our main tool is a
theory of {\it determinants} extending the one of pseudocharacters but which
works over an arbitrary base ring ; an independent aim of this paper is to
expose the main facts of this theory.  The moduli space $X$ is constructed
as the generic fiber of the moduli formal scheme of continuous formal
determinants on $G$ of dimension $d$.\par As an application to number
theory, this provides a framework to study rigid analytic families of Galois
representations (e.g.  eigenvarieties) and generic fibers of
pseudodeformation spaces (especially in the "residually reducible" case,
including when $p\leq d$).\par

\ps\ps

\section*{Introduction}
Let $G$ be a group, $A$ a commutative ring with unit and let $$T : G \longrightarrow
A$$ be a map such that $T(gh)=T(hg)$ for all $g,h \in G$. For $n\geq 1$ an
integer and $\sigma \in \got S_n$, set\footnote{This expression has the following important interpretation, due to Kostant. Assume that $T : \GL_{m}(A) \rightarrow A$ is the trace map, 
if $g_{1},\dots, g_n\in \GL_{m}(A)$ and if $\sigma \in \got{S}_{n}$, then $T^{\sigma}(g_{1},\dots,g_{n})$ coincides with the trace of the 
element $(g_{1},\dots,g_{n})\sigma$ acting on $V^{\otimes_{A}n}$, where $V:=A^{m}$.}
$T^{\sigma}(g_1,g_2,\dots ,g_n)=T(g_{i_1}g_{i_2}\dots g_{i_r})$ if $\sigma$ is the cycle $(i_1 i_2 \dots i_r)$, and in general $$T^\sigma=\prod T^{c_i}$$ if
$\sigma=c_1\dots c_s$ is the cycle decomposition of $\sigma$. The {\it $n$-dimensional pseudocharacter identity} is the relation

\begin{equation}\label{npseudo} \forall g_1,\, g_2,\,\dots ,g_n, g_{{n+1}} \in G, \,\,\, \,
\sum_{\sigma \in \got S_{n+1}}
\varepsilon(\sigma)T^\sigma(g_1,\,g_2,\,\dots,\,g_{n+1})=0,
\end{equation}
where $\varepsilon(\sigma)$ denotes the signature of the permutation
$\sigma$. We say that $T$ is a {\it $d$-dimensional pseudocharacter} of $G$ with values in $A$ if $T$ satisfies the $d$-dimensional 
pseudocharacter identity, if $T(1)=d$ and if $d!$ is invertible in $A$. \ps

The main interest of pseudocharacters lies in the close relations they share
with traces of representations : by an old result of Frobenius~\cite[p. 50]{frobenius}, the trace of
a representation $G \longrightarrow \GL_{d}(A)$ is a $d$-dimensional
pseudocharacter\footnote{In view of the previous footnote, this simply
expresses the fact that $\Lambda^{n+1}A^{n}=0$.}, and it is known that the
converse holds when $A$ is an algebraically closed field with $d!  \in
A^{*}$ (Procesi \cite{Proc2}, Taylor \cite{tay} for $\Q$-algebras,
\cite{rou} in general\footnote{As already observed in~\cite[\S 1]{bch}, let us warn the reader that although Rouquier does not require $d!$ to be invertible in $A$ in~\cite{rou}, there is a gap in the proof of his Lemma 4.1, hence of his Theorem 4.2, without this assumption. Indeed, it is not clear that each element of his ring $R$ is algebraic over $k$, as asserted on p. 580 line 2, because the polynomial $P_x$ given by his Lemma 2.13 might be identically zero if $d!$ is not invertible in $A$.}) as well as in various other situations (see below). 
In particular, we obtain this way an interesting parametrization of the
isomorphism classes of semisimple representations of $G$ over such
algebraically closed fields.  As the covariant functor from the category of
$\Z[1/d!]-$commutative algebras with unit to the category ${\rm Ens}$ of sets, which
associates to $A$ the set of $d$-dimensional pseudocharacters $G \rightarrow
A$, is obviously representable\footnote{Consider the ring $B_0$ which is the quotient of the polynomial 
$\Z[1/d!]$-algebra over the indeterminates $X_g$ for all $g \in G$, by the
ideal generated by the elements $X_{gh}-X_{hg}$ for all $g,h
\in G$. For each $\sigma \in \got S_n$ and $g_1,\cdots,g_n \in G^n$, we
have a well defined element $X^{\sigma}(g_1,\dots
,g_n) \in B_0$ defined as the class of $X_{g_{i_1}g_{i_2}\dots g_{i_r}}$ if $\sigma$
is the cycle $(i_1 i_2 \dots i_r)$, and of $\prod_i X^{c_i}(g_1,\cdots,g_n)$
in general if $\sigma=\prod_i c_i$ is the cycle decomposition of $\sigma$.
Define $B_1$ as the quotient of $B_0$ by the ideal generated by $X_1-d$ and
the elements $\sum_{\sigma \in {\got S}_{d+1}} \varepsilon(\sigma)
X^\sigma(g_1,\,g_2,\,\dots,\,g_{d+1})$ for all $g_1,\cdots,g_{d+1} \in
G^{d+1}$. The map $G \rightarrow B_1$, $g \mapsto
\overline{X_g}$, is the universal $d$-dimensional pseudocharacter of $G$.}, it turned out to be an interesting
substitute for the quotient functor ${\rm Hom}(G,{\rm GL}_{d}(-))/{\rm
PGL}_d(-)$ of isomorphism classes of $d$-dimensional representations of $G$. 
Indeed, since they have been introduced in number theory by Wiles
\cite{wiles} (when $d=2$), and by Taylor~\cite{tay} under the form above
(sometimes under the name of {\it pseudorepresentations}), they have proved
to be a successful tool, first to actually construct some (Galois)
representations, and then to study Galois representations and
Hecke-algebras.  \ps Over $\Q$-algebras, most of the basic properties of
pseudocharacters follow actually from earlier work of Procesi on invariants
of $n$-tuples of $d\times d$-matrices \cite{Proc} and on the very close
subject of Cayley-Hamilton algebras \cite{Proc2}.  In relation to
deformation theory, pseudocharacters over local rings have also been studied
by Nyssen \cite{nyssen} and Rouquier \cite{rou} in the residually
irreducible case, and by Bella\"iche-Chenevier \cite[Ch.  1]{bch} in the
residually multiplicity free case.  \ps

The first part of this paper addresses the problem of setting a
definition for an $A$-valued pseudocharacter of dimension $d$ which works
for an {\it arbitrary ring}, i.e. without the assumption that $d! \in
A^*$, and to extend to this setting most of the aformentioned results. When $d!$ is invertible, the pseudocharacter 
identity of degree $d$ is very close to the Cayley-Hamilton identity of degree $d$ defined by
the pseudocharacter $T$ and this is actually the key to most of the interesting properties of
pseudocharacters\footnote{Precisely, relation (\ref{npseudo}) is exactly $T({\rm
CH}_T(g_1,g_2,\dots,g_d)g_{d+1})=0$ (see \cite{Proc}) where
${\rm CH}_{T}$ is the multi-linearization of the "characteristic polynomial of
degree $d$ associated to $T$", which is the homogeneous polynomial  $x^d-T(x)x^{d-1}+\frac{T(x)^2-T(x^2)}{2}x^{d-2} + \dots +
\det(x)$.}; it is certainly not surprising that the definition above of pseudocharacters does not
work well in general. The key notion turned out to be the one of
{\it multiplicative homogeneous polynomial laws} on algebras, which have been studied by
Roby (\cite{roby},\cite{roby2}), Ziplies \cite{ziplies}, Ferrand~\cite{ferrand} and more recently by Vaccarino 
(\cite{vac0},\cite{vac1},\cite{vac2}), and which immediately leads to a definition
for a "generalized" pseudocharacter. To avoid confusions, we rather call them
({\it law-}) {\it determinants}. Up to the language of polynomial laws of Roby~\cite{roby} that we recall
in a preliminary \S~\ref{notations}, our definition is surprisingly simple :\ps\ps

{\bf Definition:} {\it An $A$-valued determinant on $G$ of dimension $d$ is an
$A$-polynomial law $D: A[G] \longrightarrow A$ which is homogeneous of degree $d$ and
multiplicative}. \ps\ps

Of course, usual determinants of true $A$-algebra representations $A[G] \longrightarrow
M_d(A)$ are determinants in this setting, and we shall prove various converse results. By definition, it is equivalent to give a determinant as above and a finite collection of maps $G^d \longrightarrow A$ satisfying various identities, which are in general much more complicated than the pseudocharacter identity. We explicit this point of view in the special case $d=2$ (\S~\ref{exdim2}). Thanks to Roby's works, there is also an equivalent general definition for a determinant in terms of the divided power ring $\Gamma_\Z^d(\Z[G])$ of degree 
$d$ (\cite{roby2}), which is naturally isomorphic to the, maybe more
standard, ring of invariants $(\Z[G]^{\otimes d})^{\got{S}_d}$ : an $A$-valued
determinant on $G$ of dimension $d$ is simply a ring homomorphism
$$\Gamma_\Z^d(\Z[G]) \longrightarrow A.$$
In particular, {\it the natural functor associating to $A$ the set of $d$-dimensional $A$-valued determinants is representable } by the ring $\Gamma_\Z^d(\Z[G])^{\rm
ab}\isomo ((\Z[G]^{\otimes d})^{\got{S}_d})^{\rm ab}$. 
\ps

Thanks to works of many people (Amitsur, Procesi, Donkin, Zubkov, 
Vaccarino and certainly others), much of the deepest properties of 
determinants are actually known, although it is hard to extract from 
the literature a unified picture\footnote{We thank a referee for pointing out the recent paper~\cite{dcprr} for a collection of results on Cayley-Hamilton algebras.}. In the first half of this paper, 
which may be viewed mostly as an introduction to the subject, we make an attempt to expose
the theory from the narrow point of view of determinants, trying to remain as self-contained (and coherent) as possible. \ps

In the first section, we develop the most basic properties of determinants 
(\S~\ref{polid},\S~\ref{faithCH}): polynomial identities, {\it Kernel of a determinant}, {\it faithful and Cayley-Hamilton quotients}, and properties with respect to base
change. An important role is played there and in the whole theory by a polynomial identity which is formally analogue to {\it Amitsur's formula}
(\ref{amitsurhomogene}), which expresses the determinant of a sum of
elements $x_1+\dots+x_r$ in terms of the coefficients of the characteristic
polynomials of some explicit monomials in $x_i$ ; we give an elementary proof for this 
formula for any determinant by mimicking an elegant proof in the matrix case due to 
Reutenauer-Sch\"utzenberger\footnote{As F. Vaccarino pointed out to us, very similar results had also been obtained by Ziplies in \cite{zipoubli}.} \cite{rs}. Another important step is to show that the 
faithful quotient $A[G]/\Ker(D)$ satisfies the Cayley-Hamilton identity, and for that 
we have to rely for the moment on an important result of
Vaccarino~\cite{vac0} describing $\Gamma_{\Z}^d(\Z\{X\})^{\rm ab}$ when $\Z\{X\}$ is the free ring over the 
finite set $X$ (actually, we only use that this ring is torsion free, but Vaccarino's 
result is much stronger, see \S~\ref{polid}). Using results of Procesi, we show also 
that {\it over $\Q$-algebras, determinants and pseudocharacters coincide}, but we do 
not know if this holds under the weaker assumption that $d! \in A^*$ (we do however prove it when $d=2$ and in several other cases, see Remark~\ref{rmkQd}). These last two points are actually the only places in the paper 
where we are not self-contained (see Remarks \ref{questionun} and \ref{rmkQd}). 
\ps

In the second section, we prove the analogue for determinants of the standard aformentioned results 
of the theory of pseudocharacters. The approach we follow is inspired from the one in 
\cite[Ch. 1]{bch}, but we have to face with several extra difficulties inherent to the use of 
polynomial laws and also from the presence of some inseparable extensions which occur in
characteristic $p \leq d$.
\ps\ps
{\bf Theorem A : } {\it  Let $k$ be an algebraically closed field and
$D : k[G] \longrightarrow k$ be a determinant of dimension $d$. There exists
a unique semisimple representation $\rho : G \longrightarrow \GL_d(k)$ such that
for any $g \in G$, $\det(1+t\rho(g))=D(1+tg)$}.\ps\ps

(See Theorem \ref{thmcorps}) Of course, "unique" here means 
"unique up to $k$-isomorphism". In fact, {\it if $k$ is any 
perfect field, or any field of characteristic $p>0$ such that 
either $p > d$ or $[k:k^p]< \infty$}, we show the stronger fact 
that {\it $k[G]/\Ker(D)$ is a semisimple finite dimensional 
$k$-algebra} (Theorem \ref{thmcorps2}). \ps\ps

{\bf Theorem B : }  {\it Let $A$ be a henselian local
ring with algebraically closed residue field $k$, $D : A[G] \longrightarrow A$ a determinant
of dimension $d$, and let $\rho$ be the semisimple representation attached to $D \otimes_A k$
by Theorem $A$. If $\rho$ is irreducible, then there exists a unique representation
$\tilde{\rho} : G \longrightarrow \GL_d(A)$ such that for any $g \in G$,
$\det(1+t\tilde{\rho}(g))=D(1+tg)$. } \ps\ps

(See Theorem \ref{structurethm}) Actually, we show the stronger fact 
that the biggest Cayley-Hamilton quotient of $A[G]$ is the faithful one, 
and is isomorphic to $(M_d(A),\det)$. We consider also the more general 
case where, under the assumption of Theorem B, $\rho$ is only assumed to 
be multiplicity free, and we show then that any Cayley-Hamilton quotient 
of $A[G]$ is a generalized matrix algebra in the sense of \cite[Ch. 1]{bch}, 
extending a result there. \ps 

Let us stress here that Theorems A and B should not be considered as 
original, as they could probably be deduced from earlier works of 
Procesi (\cite{Proc0},\cite{Proc3}) via the relations between determinants 
and generic matrices established by Vaccarino, Donkin and Zubkov. \ps

The last part of the second section deals with the problem deforming 
a given determinant $D_0$ to $A[\epsilon]$ with $\epsilon^2=0$ 
(\S~\ref{deformations}). The set of such deformations of $D_0$ appears 
naturally as a relative tangent space and has a natural structure of 
$A$-module. When $G$ is a topological group and $A$ a topological ring, 
we say that an $A$-valued determinant on $G$ is continuous if the 
coefficients of the characteristic polynomial $D(t-g)$ are continuous 
functions of $g \in G$. The main result here is the following 
(Prop.~\ref{finitenesstgsp}) : \ps

{\bf Proposition C : } {\it Let $k$ be a discrete algebraically closed 
field, $G$ a profinite group, $\rho : G \longrightarrow \GL_d(k)$ a 
semisimple continuous representation, and $D_0 = \det \circ \rho$. Let $p\geq 0$ be the characteristic of the field $k$.
The space of continuous deformations of $D_0$ to 
$k[\epsilon]$ is finite dimensional in the following two cases : 
\begin{itemize}
\item[(a)] $p=0$ or $p>d$, and the  continuous cohomology group $H^1_c(G,{\rm ad}(\rho))$ is finite 
dimensional over $k$, \ps
\item[(b)] $0<p \leq d$ and for each open subgroup $H \subset G$, there are only finitely many continuous homomorphisms $H \rightarrow \Z/p\Z$.
\end{itemize}}
\ps\ps

All of this being done, we are perfectly well equipped 
to study rigid analytic families of pseudocharacters. Let us assume 
from now on that $G$ is a profinite topological group, fix $d \geq 1$ 
an integer, and let $p$ be a prime number. Assume moreover that $G$ 
satisfies the following finiteness condition  : {\it For any normal open 
subgroup $H \subset G$, there are only finitely many continuous group 
homomorphisms $H \longrightarrow \Z/p\Z$}. This holds for instance when 
$G$ is topologically finitely generated, when $G$ is the absolute Galois group of a local 
field of characteristic $\neq p$ (e.g. $\Q_p)$, or when $G$ is the absolute 
Galois group of a number field with finite restricted ramification. \ps
 
Let $\An$ be the category of rigid analytic 
spaces over $\Q_p$ in the sense of Tate (see~\cite{BGR}). If $X$ is such a space, we shall denote by $\OO_X$ its structural sheaf and by $\OO(X)$ the $\Q_p$-algebra of global sections of
$\OO_X$. We equip $\OO(X)$ with the topology of uniform convergence on the open affinoids of $X$. The main aim of this paper is to study the contravariant functor 
$E^{\rm an} : \An \longrightarrow  {\rm Ens}$, which associates to a rigid space $X$ 
the set $E^{\rm an}(X)$ of continuous $d$-dimensional pseudocharacters 
$G \longrightarrow \OO(X)$. \ps\ps

{\bf Theorem D }: {\it $E^{\rm an}$ is representable by a quasi-Stein 
rigid analytic space}. \ps\ps

(See Theorem~\ref{mainthm}) This rigid analytic space might be called the
{\it $p$-adic character variety of $G$ in dimension $d$}. To show this theorem we actually start with studying 
other natural functors. First, we fix a continuous semisimple representation 
$$\rhob : G \longrightarrow \GL_d(\overline{\F}_p)$$
and whose determinant $D$ takes values in some finite field $k \subset \overline{\F}_p$. We consider 
the continuous deformation functor $F$ of $D$ to discrete artinian local 
$W(k)$-algebras with residue field $k$. Here $W(k)$ denotes the Witt ring of $k$. We prove first the following 
(Prop. \ref{repF}) : \ps\ps

{\bf Proposition E } : {\it $F$ is prorepresentable by a complete local noetherian $W(k)$-algebra $A(\rhob)$ with residue field $k$.} \ps\ps

Of course, for the noetherian property we rely on Proposition C. The ring
$A(\rhob)$ is constructed as a certain profinite completion of
$$\Gamma_\Z^d(\Z[G])^{\rm ab} \otimes_\Z W(k).$$

We consider then the functor $E$ from the category of formal schemes over 
$\Spf(\Z_p)$ to sets, which associates to $\cal X$ the set of continuous 
$d$-dimensional $\OO(\cal X)$-valued determinants on $G$. We can attach to 
each such formal determinant a subset of "residual determinants". The set 
$|G(d)|$ of all residual determinants is in natural bijection with the set 
of (determinants of the) continuous semisimple representations $\rhob$ as 
above, taken up to isomorphism and Frobenius actions on coefficients. It 
turns out that that {\it the subfunctor $E_{\rhob} \subset E$ parameterizing 
determinants which are residually constant and "equal to" $\det \circ \rhob$ 
is representable} and isomorphic to the affine formal scheme $\Spf(A(\rhob))$ 
over $\Spf(\Z_p)$ (where $A(\rhob)$ is equipped with the $m$-adic topology of 
given by its maximal ideal $m$). Our main second result is then the following 
(\S~\ref{analyticfinal}), which implies Thm. D :\ps\ps

{\bf Theorem F : } {\it The functor $E$ is representable by the disjoint union of the $\Spf(A(\rhob))$, for $\rhob \in |G(d)|$. 
The functor $E^{\rm an}$ is canonically isomorphic to the generic fiber of $E$ in the sense of Berthelot.}\ps\ps

If we fix an isomorphism $W(k)[[t_1,\dots,t_h]]/I \isomo A(\rhob)$, then we get a closed immersion
$$\Spf(A(\rhob))^{\rm rig} \hookrightarrow {\mathbb B}_{[0,1[}^h$$
as the closed subspace of the open $h$-dimensional unit ball defined by
$I=0$, and $X$ is then a disjoint union of such spaces.\ps

In section~\ref{complements}, we give some general complements about the rigid analytic space $X$ representing $E^{\rm an}$. 
For instance, consider the functor $$E^{\rm irr} : \An \rightarrow  {\rm Ens}$$ which associates to any rigid space $X$ the set of 
isomorphism classes of pairs $(R,\rho)$ where $R$ is an Azumaya $\OO_{X}$-algebra of rank $d^{2}$ and $\rho : G \rightarrow R^{*}$ is a 
continuous group homomorphism such that for all closed points $x \in X$, the evaluation $\rho_{x} : 
G \rightarrow R_{x}^{*}$ is absolutely irreducible (see \S~\ref{univcohCH}). \ps\ps

{\bf Proposition G : } {\it $E^{\rm irr}$ is representable by a Zariski-open subspace of $E^{\rm an}$ equipped with its universal Cayley-Hamilton representation.}
\ps\medskip
In the last section~\ref{sectexemple}, we give an application of some of the previous results to Galois representations. 
Let $G$ be the Galois group of a maximal algebraic extension of $\Q$ unramified outside $\{2,\infty\}$, and let $X$ be 
the $2$-adic analytic space parameterizing the $2$-dimensional rigid analytic pseudocharacters of $G$ (so $p=d=2$). This space $X$ is 
an admissible disjoint union of three open subspaces $X^{\rm odd}$, $X^{+}$ and $X^{-}$ over which the trace 
of a complex conjugation of $G$ is respectively $0$, $2$ and $-2$. \ps\ps

{\bf Theorem H : } {\it $X^{\rm odd}$ (resp. $X^{\pm}$) is the open unit ball of dimension $3$ (resp. dimension $2$) over $\Q_{2}$. }

\ps\ps\ps\ps

The author would like to thank Emmanuel Breuillard, Claudio Procesi and Francesco Vaccarino for some useful
discussions, a referee for his careful reading, as well as Jean-Pierre Serre and Jo\"el Bella\"iche for their remarks.\ps\ps\ps

\section{Determinants of algebras}

\subsection{Homogeneous multiplicative $A$-polynomial laws}\label{notations}
We need some preliminaries about polynomial laws between two modules. We refer to \cite{roby} and \cite{roby2} for the proofs of
all the results stated below. \ps \ps 

Let $A$ be a commutative unital ring, and let $M$ and $N$ be two $A$-modules. Let $\cal C_A$ be the
category of commutative $A$-algebras. Each $A$-module $M$ gives rise to a
functor $\underline{M}: \cal C_A \longrightarrow {\rm Ens}$ via the formula $B \mapsto M
\otimes_A B$. An {\it $A$-polynomial law} $P : M
\longrightarrow N$ is a natural transformation $\underline{M}
\longrightarrow \underline{N}$. In other words, it is a collection of maps $$P_B: M\otimes_A B
\longrightarrow N \otimes_A B,$$ where $B$ is any commutative $A$-algebra,
which commute with any scalar extension $B \rightarrow B'$ over $A$. By a
slight abuse of notations, if $B$ is a commutative $A$-algebra and $m \in
M\otimes_A B$ we shall often write $P(m)$ for $P_B(m)$. When
$B=A[T_1,\dots,T_s]$, we shall write $M[T_1,\dots,T_s]$ for $M\otimes_A
A[T_1,\dots,T_s]$. \ps\ps

We refer to \cite{roby} for the basic operations that we can do with polynomial
laws. If $B$ is a commutative $A$-algebra and $P : M \longrightarrow N$ is an
$A$-polynomial law, we will denote by $P \otimes_A B : M \otimes_A B \longrightarrow
N\otimes_A B$ the natural induced $B$-polynomial law\footnote{By definition, if $C$ is a commutative $B$-algebra,
$(P \otimes_A
B)_C=P_C$ via the isomorphism $(-\otimes_A B)\otimes_B C=-\otimes_A
C$.}.\ps\ps

We say that $P$ is homogeneous of degree $n$ (an integer $\geq
0$) if $P(bx)=b^nP(x)$ for all object $B$ in $\cal C_A$, $b \in B$ and $x \in M\otimes_A B$. 

\begin{example}\label{examplehomogene}{\rm Let $P: M \longrightarrow N$ be 
an homogeneous $A$-polynomial law of degree $n$. \begin{itemize}\ps 
\item[(i)] When $n=1$ (resp. $n=0$), $P_A$ is an $A$-linear\footnote{Let $X,Y,T$ be indeterminates. If $u,v \in M$, then $P(uX+vY) \in N[X,Y]$. As $P$ has degree $1$, sending $(X,Y)$ to $(XT,YT)$ shows that $P(uX+vY)$ is of the form $a(u,v)X+b(u,v)Y$ for some well-defined functions $a,b : M^2 \rightarrow N$. By evaluating $(X,Y)$ at $(1,0)$, $(0,1)$ and $(1,1)$, we obtain respectively $a(u,v)=P(u)$, $b(u,v)=P(v)$ and $P(u+v)=P(u)+P(v)$.} map and $P_B = P_A \otimes_A B$
(resp. $P_B=P_A(0) \otimes 1$ is a constant), and $P \mapsto P_A$ induces a bijection between $A$-polynomial laws of degree $1$ (resp. $0$) and $\Hom_A(M,N)$ (resp. $N$).
\ps \item[(ii)]When $n=2$, $P_B$ is again uniquely
determined by $P_A$, which is any map $q : M \longrightarrow N$ such that
$q(am)=a^2q(m)$ for all $a \in A$, $m \in M$, and such that $(m,m')\mapsto
q(m+m')-q(m)-q(m')$ is $A$-bilinear. \ps \item[(iii)] When $n\geq 3$, $P_A$ does not
determine $P_B$ in general. For instance, let $A$ be the finite field $\F_q$ with $q$ 
elements, $M=\F_q^2$ and let $X,Y$ be an $A$-basis of $\Hom_A(M,A)$. The $A$-polynomial law 
$P : M \rightarrow A$ defined by $P=XY^q-X^qY$ is homogeneous of degree $q+1$, we have
$P_A=0$ but $UV^q-U^qV \in P(M[U,V]) \, \neq 0$. \end{itemize}}\end{example} 

In any cases, a homogeneous $P$ of degree $n$ is uniquely
determined by $P_{A[T_1,\dots,T_n]}: M[T_1,\dots,T_n] \longrightarrow
N[T_1,\dots,T_n]$. Precisely, if $X \subset M$ generates $M$ as
$A$-module, then such a $P$ is uniquely determined by the (finite) set
of functions $$P^{[\alpha]}: X^n \longrightarrow N,$$ with $\alpha \in
I_n=\{(\alpha_1,\alpha_2,\dots,\alpha_n) \in \N^n, \,\,\alpha_1+\dots+\alpha_n=n\}$, defined by
the relation 
$$P ( \sum_{i=1}^n T_i x_i ) = \sum_{\alpha \in I_n}
P^{[\alpha]}(x_1,\dots,x_n)T^\alpha,$$
\noindent where $T^\alpha=\prod_{i=1}^n T_i^{\alpha_i}$. \ps \ps

We denote by ${\cal
P}^n_A(M,N)$ the $A$-module of homogeneous $A$-polynomial laws of degree $n$
from $M$ to $N$. The functor ${\cal P}^n_A(M,-) : {\rm Mod}(A) \rightarrow
{\rm Mod}(A)$ is representable by the
usual $A$-module $\Gamma^n_A(M)$ of divided powers of order $n$ on $M$ relative to $A$
(\cite[Thm. 4.1]{roby}). Let us recall that $\Gamma_A^n(M)$ is naturally isomorphic to the $n^{\rm
th}$-graded piece of the commutative $A$-algebra $\Gamma_A(M)$ which is
generated by the symbols $m^{[i]}$ for $m \in M$ and $i\geq 0$, with the usual
homogeneous relations : \begin{itemize}\ps
\item $m^{[0]}=1$ for all $m \in M$,\ps
\item $(am)^{[i]}=a^im^{[i]}$ for all $a \in A$ and $m \in M$,\ps
\item $m^{[i]}m^{[j]}=\frac{(i+j)!}{i!j!}m^{[i+j]}$ for all $i,j \geq 0$ and
$m \in M$,\ps
\item $(m+m')^{[i]}=\sum_{p+q=i}m^{[p]}{m'}^{[q]}$ for all $i\geq 0$ and $m,
m' \in M$. \ps
\end{itemize}\ps
The natural map $P^{\rm univ} : m \mapsto m^{[n]}$, $M \longrightarrow \Gamma_A^n(M)$, induces {\it 
the universal homogeneous $A$-polynomial law of degree $n$}. 
For $\alpha \in I_n$ as above, $(P^{\rm univ})^{[\alpha]}(m_1,\dots,m_d)=\prod_{j=1}^d m_j^{[\alpha_j]}$.
\ps\ps

Let $R$ and $S$ be two $A$-algebras\footnote{By an $A$-algebra
we shall always mean an associative and unital $A$-algebra (but not
necessarily commutative).}, and $P : R \longrightarrow S$
be a homogeneous $A$-polynomial law of degree $n$. We say that $P$ is {\it multiplicative} if
$P(1)=1$
and if $P(xy)=P(x)P(y)$ for all $B$ and $x, y \in R\otimes_A B$. For example, the homogeneous multiplicative $A$-polynomial laws of degree $1$ are the $A$-algebra homomorphisms. By \cite{roby2}, the
structure of $A$-algebra on $R$ induces an $A$-algebra structure on\footnote{This structure is not to be confused with the $A$-algebra
structure on $\Gamma_A(R)$, which is always graded and commutative. The
$A$-algebra $\Gamma_A^n(R)$ is commutative if $R$ is, its neutral element is $1_{R}^{[n]}$, and $\Gamma_A^1(R)=R$.} $\Gamma_A^n(R)$, 
and it turns out that the functor ${\cal M}_A^n(R,-)$, from $A$-algebras to sets, that associates to any $A$-algebra $S$ the set of ${\cal M}_A^n(R,S)$ of
$n$-homogeneous multiplicative $A$-polynomial laws from $R$ to $S$, is
representable by the $A$-algebra $\Gamma_A^n(R)$ (\cite[Th\'eor\`eme]{roby2}). In particular, the universal 
homogeneous $A$-polynomial law $$P^{\rm univ} : R \longrightarrow \Gamma_{A}^{n}(R), \, \, \, \, \, r \mapsto r^{[n]},$$
is multiplicative.
\begin{remark}\label{remTS}{\rm Let $M$ be an $A$-module and let $\TS^n_A(M)$ be the $A$-submodule of $M^{\otimes_A^n}$ invariant by the symmetric group $\got{S}_n$. The natural map $M \longrightarrow \TS_A^n(M)$, $m \mapsto m^{\otimes n}$, induces a homogeneous $A$-polynomial law of degree $n$, hence there is a natural $A$-linear map \begin{equation}\label{GTS}\Gamma_A^n(M) \longrightarrow \TS_A^n(M),\end{equation} 
which is actually an isomorphism if $M$ is free as $A$-module (\cite[Prop. IV.5]{roby}). 
When $M=R$ is an $A$-algebra, $\TS_A^n(R)$ has an obvious $A$-algebra structure and $r \mapsto r^{\otimes n}$ is clearly multiplicative, so (\ref{GTS}) is actually an $A$-algebra homomorphism. In particular, if $R$ is free as $A$-module, then 
$$\Gamma_A^n(R) \isomo \TS_A^n(R)$$ is an $A$-algebra isomorphism.}
\end{remark}

\ps

\begin{remark}\label{remextbase}{\rm If $B$ is a commutative $A$-algebra and $M$ an $A$-module, the homogeneous 
$A$-polynomial law of degree $n$ $$M \longrightarrow \Gamma_B^n(M\otimes_A B), \, \, m \mapsto  (m\otimes 1)^{[n]},$$ 
induces an isomorphism (\cite[Thm. III.3]{roby}) $\Gamma_A^n(M)\otimes_A B \isomo \Gamma_B^n(M\otimes_A B)$. 
When $M=R$ is an $A$-algebra, this latter isomorphism is a $B$-algebra homomorphism as the polynomial law above is multiplicative.}
\end{remark}

\ps\ps

\subsection{Definition of a determinant}\label{defdet}
Let $R$ be any $A$-algebra and $d \geq 1$ an
integer.\ps\medskip
{\bf Definition : }\label{defpseudodet} {\it A $d$-dimensional $A$-valued
determinant on $R$ is an element of ${\cal M}_A^d(R,A)$, {\it i.e. } a multiplicative $A$-polynomial law $\Det: R \longrightarrow
A$ which is homogeneous of degree $d$. When $R=A[G]$ for some group $G$ (or
unital monoid), we
say also that $\Det$ is a determinant on $G$. }\ps\medskip

Of course, if $R=M_d(A)$ (resp. any Azumaya algebra of rank $d^2$ over its
center $A$), the usual determinant $\det : M_d(A) \longrightarrow A$ (resp. the reduced norm) induces in the obvious way\footnote{For any commutative $A$-algebra $B$, define $\det_{B}$ as the determinant $M_{d}(B) \rightarrow B$.} 
a determinant of dimension $d$. In particular, for any $A$-algebra homomorphism $\rho : R
\longrightarrow M_d(A)$, $$\Det := \det \circ \rho$$ is a
$d$-dimensional $A$-valued determinant on $R$. In section~\ref{structurethm}, we will
prove some converse to this construction. For example, we will show that when $A$ is an algebraically closed field, any determinant of $R$
is of the form above, and we will also study the case when $A$ is a
local henselian ring. When $d=1$, a determinant $\D : R \longrightarrow A$ of dimension $1$ is by definition the same as an $A$-algebra homomorphism (see Example~\ref{examplehomogene} (i)). \ps\ps

Let $\det_A(R,d): \cal C_A \longrightarrow {\rm Ens}$ be the covariant
functor associating to any commutative $A$-algebra $B$, the set  of
$B$-valued determinants $R\otimes_A B \longrightarrow B$ of dimension $d$, which is the
same as the set of multiplicative homogeneous $A$-polynomial laws $R \longrightarrow B$ of dimension $d$
(recall that $\cal
M_A^d(R,B) \isomo \cal M_B^d(R \otimes_A B,B)$ by Remark~\ref{remextbase}).
It is equivalent to give such a law or an $A$-algebra homomorphism $\Gamma_A^d(R) \longrightarrow B$, which
necessarily factors through its abelianization\footnote{By definition, the
abelianization of a ring $R$ is the quotient of $R$ by the two-sided ideal
generated by the $xy-yx$ with $x,y \in R$.} $\Gamma_A^d(R)^{\rm ab}$, hence
we get the :

\begin{prop}\label{representabilite} $\det_A(R,d)$ is representable by the $A$-algebra $\Gamma_A^d(R)^{\rm ab}$. 
\end{prop}

In particular, when $R=\Z[G]$, then $\det_\Z(\Z[G],d)$ is representable by $$\Gamma_\Z^d(\Z[G])^{\rm
ab}\isomo \TS_\Z^d(\Z[G])^{\rm ab},$$ that we shall simply denote by
$\Z(G,d)$. This ring is nonzero thanks to the trivial representation of dimension $d$ of $G$.
We will set also $$X(G,d)={\rm Spec}(\Z(G,d)).$$
Of course, if $S$ is any scheme, we may define a determinant of dimension
$d$ on $G$ over $S$ as an $\OO(S)$-valued determinant of dimension $d$ on
$G$, and $X(G,d)$ obviously still represents this extended determinant functor. 

\begin{example}\label{exziplies}{\rm \begin{itemize}\ps
\item[(i)] When $R=A[X]$ is a polynomial ring in one variable, then $\TS_A^d(R)=A[X_1,\dots,X_d]^{\got S_d}=A[\Sigma_1,\dots,\Sigma_d]$ by the classical theorem on symmetric polynomials (with the obvious notations for $X_i$ and $\Sigma_j$). In particular, $\Gamma_A^d(R)=\Gamma_A^d(R)^{\rm ab} \simeq A[\Sigma_1,\dots,\Sigma_d]$ by Remark~\ref{remTS}. As we will see in \S~\ref{polid}, the universal determinant is the determinant of the regular representation of $A[X]$ on
$$\Gamma_A^d(R)[X]/(X^d-\Sigma_1 X^{d-1} + \Sigma_2 X^{d-2} - \dots + (-1)^d \Sigma_d).$$\ps
\item[(ii)] When $R$ is an Azumaya algebra of rank $d^{2}$ over its center $A$, a result of Ziplies \cite{ziplies} (see also Ex.~\ref{zipliesbis}) shows that the
reduced norm is the unique $A$-valued determinant of dimension $d$ of $R$, and even
that the reduced norm induces an $A$-algebra isomorphism $\Gamma_A^d(R)^{\rm
ab} \isomo A$. \ps
\item[(iii)] When $G$ is a finite group $\Z(G,d)$ is a finite $\Z$-algebra (as
$\Gamma_\Z^d(\Z[G])$ is free of finite type as $\Z$-module).\ps
\item[(iv)] Using Remark~\ref{remextbase}, we get that if $B$ is any commutative
$A$-algebra, the natural $A$-algebra
homomorphism $B \otimes_A \Gamma_A^d(R)^{\rm ab} \longrightarrow
\Gamma_B^d(B \otimes_A R)^{\rm ab}$ is an isomorphism.\ps
\end{itemize}
}\end{example}

\ps\ps

In the case $R=A[G]$, it is equivalent to give a determinant $A[G] \longrightarrow A$ of dimension $d$ 
and a $d$-homogeneous multiplicative polynomial law $\Z[G] \longrightarrow A$.
Such a law is uniquely determined by the set of functions 
$$\Det^{[\alpha]} : G^d
\longrightarrow A, \, \, \alpha \in I_d,$$ 
which satisfy a finite number of identities coming from the
requirement that the map $$\prod_{j=1}^dg_j^{[\alpha_j]} \mapsto
\Det^{[\alpha]}(g_1,\dots,g_d) \in A$$ extends to a ring homomorphism
$\Gamma^d_\Z(\Z[G]) \longrightarrow A$. \ps\ps

\begin{example}\label{exdim2} ({\it Determinants of dimension $2$ on a group $G$
(or a unital monoid)})
{\rm  As an example, let us specify a bit those relations when $d=2$. In this case,
we may write $${\Det}(gU+hV)=D(g)U^2+f(g,h)UV+D(h)V^2$$ for some
functions $D=\Det_{|G} : G \longrightarrow A$ and $f : G \times G \longrightarrow A$.
As we are in degree $2$, any pair of such functions determines a unique
homogeneous $\Z$-polynomial law of degree $2$ from $\Z[G]$ to $A$, under the (obviously necessary) assumptions : $$\forall g \in G,\, \, \,  f(g,g)\,=\,2\,D(g), \, \, \, \, \, \, \forall  g, h\, \in G,\, \, f(g,h)\,=\,f(h,g).$$ We claim that
given $D$ and $f$ satisfying this condition, the axiom of multiplicativity of $\Det$ is equivalent to
the following set of conditions : \begin{itemize}\ps
\item[(i)] $D$ is a group homomorphism $G \longrightarrow A^*$ (in particular $D(1)=1$),  \ps
\item[(ii)] for all $g, h, h' \in G$, $f(hg,h'g)=f(h,h')D(g)$,
 \ps
\item[(iii)] for all $g, g', h, h' \in G$,
$f(hg,h'g')+f(hg',h'g)=f(h,h')f(g,g')$.\ps
\end{itemize}
Indeed, assuming that $D$ is a group homomorphism, condition (ii) means that
$\D(xg)=\D(x)\D(g)$ for all $x \in \Z[G]$ and $g \in G$. Assuming this relation,
condition (iii) means that $\D(xy)=\D(x)\D(y)$ for all $x,y \in R$.
Obviously, this multiplicativity property extends automatically to $\Det_B$
for all commutative $A$-algebras $B$.  \ps
We can write these conditions in a slightly different way. Define $T: G \longrightarrow A$ by the formula
$$T(g)=f(g,1).$$ 
Applying (iii) to $g'=h'=1$, we see that $T(1)=2$ and for all $g, h \in G$ we 
$$f(g,h)=T(h)T(g)-T(hg),$$ and in particular $T(gh)=T(hg)$. Morever, 
$f(g,h)=D(h)T(gh^{-1})$ by (ii). 

\begin{lemma}\label{d2aussi} The above map $\Det \mapsto (T,D)$ induces a bijection between the set of $A$-valued determinants of $G$ of dimension $2$ and the set of pairs of functions $(T,D) : G \rightarrow A$ such that $D : G \rightarrow A^*$ is a group homomorphism, $T : G \rightarrow A$ is a function with $T(1)=2$, and which satisfy for all $g, h\in G$: \begin{itemize}
\item[(a)] $T(gh)=T(hg)$,\ps
\item[(b)] $D(g)T(g^{-1}h)-T(g)T(h)+T(gh)=0.$
\end{itemize}\end{lemma} 

\ps
\noindent The lemma follows easily once we observe that assuming (ii),
it is enough to check (iii) for $g'=h'=1$. Note that applying (iii) to
$(h,h',g,g')=(g_1,1,g_2,g_3)$ we obtain $\forall g_1,g_2,g_3 \in G$ \begin{center} {\tiny $T(g_1)T(g_2)T(g_3)-T(g_1)T(g_2g_3)-T(g_2)T(g_1g_3)-T(g_3)T(g_1g_2)+T(g_1g_2g_3)+
T(g_1g_2g_3)=0$},\end{center}
which is the pseudocharacter relation of dimension $2$ for $T$. We will see in Prop.~\ref{casedim2} the following converse
result : {\it Assume that $2$ is invertible in
$A$. Let $T: G \longrightarrow A$ be a map such that $T(1)=2$, $T(gh)=T(hg)$ for all $g,h \in G$,
and that satisfies the $2$-dimensional pseudocharacter identity. If we set 
$D(g)=\frac{T(g)^2-T(g^2)}{2}$, then $(D,T)$ defines a determinant of $G$ of
dimension $2$.} The non-trivial part is to show that $D$ is a group
homomorphism. }\ps 
\end{example}

\subsection{Some polynomial identities}\label{polid} Let $R$ be an
$A$-algebra, $B$ a commutative $A$-algebra, and $\Det \in {\cal
M}_A^d(R,B)$.  For each $r \in R$, we define the {\it characteristic
polynomial} $\chi(r,t) \in B[t]$ of $r$ by the formula
$$\chi(r,t):=\Det(t-r)=:\sum_{i=0}^d (-1)^i \Lambda_i(r) t^{d-i}.$$ This
formula defines $A$-polynomial laws $\Lambda_i : R \longrightarrow B$ which
are homogeneous of degree $i$, for $i\geq 0$.  We have $\Lambda_0=1$,
$\Lambda_d=\Det$, $\Lambda_i=0$ for $i\geq d+1$, and $\Lambda_1$ is an
$A$-linear map, that we shall also denote by $\Tr$ and call the {\it trace}
associated to $\Det$.  \ps When $B=A$, in which case $\Det$ is a
determinant, this defines as well a homogeneous $A$-polynomial law of degree
$d$ $$\chi(r) : R \longrightarrow R, r \mapsto
r^d-\Lambda_1(r)r^{d-1}+\Lambda_2(r)r^{d-2}+\dots+(-1)^d\Lambda_d(r).$$ 
If $n \geq 0$ is an integer, we shall denote by $I_{n,d}$ the set
of $\alpha=(\alpha_1,\alpha_2,\dots,\alpha_n) \in \N^n$ such that
$\sum_{i=1}^n \alpha_i = d$. We will need to consider for each $\alpha \in
I_{n,d}$ the $A$-polynomial law
$\chi_\alpha : R^n \longrightarrow R$ defined by the
following identity in $R[t_1,\dots,t_n]$ : $$\chi(t_1 r_1 + \dots + t_n
r_n)=\sum_{\alpha} \chi_\alpha(r_1,\dots,r_n) t^\alpha,$$
where $t^\alpha=\prod_{i=1}^n
t_i^{\alpha_i}$. 

\begin{example} {\rm \begin{itemize}
\item[(i)] Let us go back to the case $R=A[X]$ (Example~\ref{exziplies} (i)). We already identified $\Gamma_A^d(R)$ with the $A$-algebra $A[\Sigma_1,\dots,\Sigma_d]$, so any homogeneous multiplicative $A$-polynomial law $\Det : R \longrightarrow B$ of degree $d$ is uniquely determined by the image $\Sigma_i(\Det)$ of $\Sigma_i$ in $B$. Unravelling the definitions, we see that $\Sigma_i(\Det)=\Lambda_i(X)$, hence the claim in Example~\ref{exziplies} (i).\ps
\item[(ii)] If $\Det : R \longrightarrow B$ is a homogeneous multiplicative $A$-polynomial law 
of degree $d$, and $r \in R$, we can restrict it to $A[X]$ via the $A$-algebra homomorphism $A[X] \longrightarrow R, \, X \mapsto r$. 
We get this way, and by the previous example, all the possible identities satisfied by determinants of polynomials over a single element of $R$. For example, the {\it Newton relations} hold,  {\it i.e.} for all $r \in R$ we have the following equality in $B[[t]]$ :  
\begin{equation}\label{newtrel} -t \frac{\frac{\partial}{\partial t}\Det(1-tr)}{\Det(1-tr)} = \sum_{n \geq 1} \Tr(r^n) t^n.
\end{equation}
\end{itemize}}
\end{example}

All the functions defined above satisfy a number of polynomial identities, we collect some of them in the following lemma.

\begin{lemma}\label{formules}\label{idpoly} \begin{itemize}
\item[(i)] For all $r, r' \in R$, $\Det(1+rr')=\Det(1+r'r)$.\ps
\item[(ii)] For all $r_1, \dots, r_n \in R$ and $i\geq 0$, $\Lambda_i(r_1+r_2+\dots+r_n)$ satisfies Amitsur's formula\footnote{See formula~(\ref{amitsurhomogene}).}.
\item[(iii)] $\Tr$ satisfies the $d$-dimensional ($B$-valued) pseudocharacter
identity.  \ps
\item[(iv)] If $B=A$, then for all $r, r_1,\dots,r_n \in R$ and all $\alpha
\in I_{n,d}$, $\Det(1+\chi_\alpha(r_1,\dots,r_n)r)=1$.
\end{itemize}
\end{lemma}

\medskip

Let $r, r' \in R$. We want to check that $\Det(1+rr')=\Det(1+r'r)$.  Note
that if $r$ is invertible\footnote{By invertible we shall always mean on
both sides.} in $R$, then this follows from the multiplicativity of $\Det$ and
the commutativity of $B$ :
$$\Det(1+rr')=\Det(r)\Det(r^{-1}+r')=\Det(r^{-1}+r')\Det(r)=\Det(1+r'r).$$
We reduce to this case as follows.  Set $r'=1+u$ and let us work in $R[t]$. 
We claim that $$\Det(1+(1+tu)r)=\Det(1+r(1+tu)) \in B[t],$$ which will
conclude the proof by evaluating $t$ at $1$.  But this is an equality of
polynomials in $t$ with degree less than $d$, so it is enough to show that
they coincide in $B[t]/(t^{d+1})$.  But $1+tu$ is invertible in $R\otimes_A
A[t]/(t^{d+1})$ hence we are reduced to the previous argument.\ps\ps Let us
now prove Amitsur's formula.  We mimic here (and actually for (iii) and (iv)
below also) the beautiful argument of \cite{rs}\footnote{We are grateful to
Emmanuel Breuillard for pointing out this reference to us.}.  \ps Let $n\geq
1$ be any positive integer, $X=\{x_1 < x_2 < \dots < x_n\}$ a totally
ordered alphabet, and $X^+$ the monoid of words in $X$ equipped with the
induced (total) lexicographic ordering $\leq$, with the convention that
$\emptyset < x_i$ for each $i$.  Recall that a word $w \in X^+$ is a {\it Lyndon word}
if $w \leq w'$ for any suffix\footnote{Recall that $w'$ is a suffix of $w$
if $w=mw'$ for some word $m$.} $w'$ of $w$ (see \cite[Ch.  5]{loth}). 
Denote by $\cal L$ the set of Lyndon words.  By Lyndon's theorem, any word
$w$ writes uniquely as a product of Lyndon words $w=w_1w_2\dots w_m$ where
$w_1 \geq w_2 \geq \dots \geq w_m$ (Lyndon factorization of $w$).  This allows to define a
{\it sign} map
$$\varepsilon:
X^+ \longrightarrow \{\pm 1\}$$ as follows. If $w \in X^+$ is a Lyndon word, set
$\varepsilon(w)=(-1)^{\ell(w)-1}$ where $\ell(w)$ is the length of the word
$w$. If $w \in X^+$ is any word, with Lyndon factorization $w=w_1w_2\dots
w_m$, we set $\varepsilon(w)=\prod_{i=1}^m \varepsilon(w_i)$. \ps We fix now some elements $r_1,\dots,r_n$ in $R$, and consider the
$A$-algebra $$A_m=A[t_1,\dots,t_n]/(t_1,\dots,t_n)^m.$$ Lyndon's theorem
writes then as the following equality in $R \otimes_A A_m$ $$
\frac{1}{1-(t_1r_1+\dots+t_nr_n)}=\prod_{w}\frac{1}{1-w}, $$ where the
product on the right hand side is taken over the finite set of Lyndon words of length
$<m$ on the alphabet $\{ t_1r_1 < \cdots < t_nr_n\}$, chosen in the decreasing order.  Applying
$\Det$ and inverting, we get the following equality in $B \otimes_A A_m$
\begin{equation}\label{amitsurformula}\Det(1-(\sum_{j=1}^n t_j
r_j))=\prod_{w \in \cal L} \left(\sum_{i=0}^d (-1)^i \Lambda_i(w)\right)
\end{equation} where the product on the right hand side is now taken over all the
Lyndon words on the $t_ir_i$, which is a well defined element in
$B[[t_1,\dots,t_n]]$. Moreover, the term on the left is the image via
$B[t_1,\dots,t_n] \longrightarrow B\otimes_A A_m$ of the polynomial
$\Det(1-(\sum_{j=1}^n t_j r_j)) \in B[t_1,\dots,t_n]$ which does not depend
on $m$.  As a consequence, the formula (\ref{amitsurformula}), also called
{\it Amitsur's formula}, holds in $B[[t_1,\dots,t_n]]$, for any integer $n$. 
If $i\geq 0$ is any integer, the homogeneous part of degree $i$ of this
equality is (for any $n$)

\begin{equation}\label{amitsurhomogene}
\Lambda_i(t_1r_1+\dots+t_nr_n)=\sum_{\ell(w)=i}
\epsilon(w)\Lambda(w),\end{equation}
where the sum is extended over the $n^i$ words $w$ on the $t_jr_j$ with length $i$,
and if $w=w_1^{l_1}\dots w_q^{l_q}$ is the Lyndon factorization of $w$ with $w_1>w_2>\cdots>w_q$,
where
$$\Lambda(w)=\Lambda_{l_q}(w_q)\cdots\Lambda_{l_2}(w_2)\Lambda_{l_1}(w_1).$$ 
Indeed, observe that for each such word, we have $\ell(w)=i=\sum_{k=1}^q l_k \ell(w_k)$,
thus $\varepsilon(w)=(-1)^{(\sum_{k=1}^q l_k)-i}$. Equality (\ref{amitsurhomogene}) holds a priori in $B[[t_1,\dots,t_n]]$ but both sides belong to $B[t_1,\dots,t_n]$, hence it 
obviously holds in $B[t_1,\dots,t_n]$. Sending each $t_i$ to $1$, we finally get Amistur's formula for $\Lambda_i(r_1+\dots+r_n)$ (in
$B$). 
\ps\ps
Let us check now part (iii) of the Lemma. Let us look at Amitsur's formula~(\ref{amitsurhomogene}) with $i=n=d+1$, and 
consider its homogeneous component
with degree $1$ in each $t_j$. We see at once that it is exactly the
$d+1$-dimensional pseudocharacter identity for $\Lambda_1=\Tr$.
\ps\ps

\begin{remark}\label{mult}{\rm Assume more generally that $B$ is any associative $A$-algebra (non necessarily commutative) and $\Det \in {\cal M}_d(R,B)$. Then the definition 
of the $\Lambda_i$ also makes sense in this extended context and the same proof as above shows that Amitsur's formula~\ref{amitsurhomogene} still holds (the {\it increasing ordering} chosen in the definition of $\Lambda(w)$ is important in this case). However, assertion (iv) only makes sense when $A=B$.}
\end{remark}

To prove assertion (iv), it amounts to show that $\Lambda_i(\chi_\alpha(r_1,\dots,r_n)r)=0$ for all $r, r_1, \dots, r_n \in R$ and $i\geq 1$. We will prove it now only for $i=1$. As $\Lambda_1$ is $A$-linear, and replacing $R$ by $R[t_1,\dots,t_n]$, it is enough to show that $\Lambda_1(\chi(r)r')=0$ for all $r, r' \in R$. Let us look at Amitsur's formula~(\ref{amitsurhomogene}) with $i=d+1$, $n=2$ and $(r_1,r_2)=(r,r')$, consider its homogeneous component
with degree $d$ in $t_1$ and $1$ in $t_2$. Each word in the sum has the form $r^ar'r^b$ with $a+b=d$, whose Lyndon factorization is $(r^ar')(r)^b$, and whose sign is $(-1)^a$. As $\Lambda_{d+1}=0$ we get an equality  
$$0 = \sum_{a+b=d} (-1)^a \Lambda_1(r^ar')\Lambda_b(r)=\Lambda_1(\left(\sum_{a=0}^d (-1)^a\Lambda_{d-a}(r)r^a\right)r')$$ 
what we wanted to show. \ps\ps
We still have to complete the proof of identity (iv), but before let us give a simple consequence of what we already proved.

\begin{cor}\label{corimage} Let $\Det$ be an $A$-valued determinant on $G$ (a monoid) of dimension $d$ and
$B \subset A$ the subring generated by the coefficients $\Lambda_i(g)$ of $\chi(g,t)$ for all $g \in G$. Then $\Det$ factors through a (unique) $B$-valued determinant on $G$ of dimension $d$.
\end{cor}

\begin{pf} We have to show that for all $g_1,\dots,g_n \in G$, $\Det(g_1t_1+\dots+g_n t_n) \in B[t_1,\cdots,t_n]$. By Amitsur's
formula~(\ref{amitsurhomogene}) such a determinant is a signed sum of monomials in $\Lambda_i(w)$ where $w$
is a word on the $g_i$, in particular $w \in G$, and we are done.
\end{pf}

We now come back to the proof of part (iv). Although it might be possible to prove it in the style above, we will rather deduce it from a general theorem of Vaccarino. 
Actually, as we shall see, the multiplicativity assumption on $\Det$ is
strong enough to imply that all the polynomial identities between the
$\Lambda_i(w)$ (where $w$ is a word in elements of $R$) which hold for the
determinant of matrix algebras also hold for $\Det$. These identities have to hold in principle 
in the universal ring $\Gamma_\Z^n(\Z\{X\!\})^{\rm ab}$, where
$\Z\{X\!\}=\Z\{x, x \in X\}$ is the free ring over a set $X$ (e.g. $X=\N$), but it might be a bit tedious in practice to compute in this ring. All we will need to know is actually contained in the aforementioned beautiful result
of Vaccarino (relying on results of Donkin \cite{Do} and Zubkov) that we
explain now.  \ps\ps

{\it Vaccarino's result (\cite[Thm 6.1]{vac0}, \cite[Thm 28]{vac1}). } Let $X$ be any set, $\Z\{X\!\}$ as above, and $F_X(d)=\Z[x_{i,j}]$ the ring of polynomials on the variables $x_{i,j}$
for all $x \in X$ and $1\leq i, j \leq d$. We have the natural {\it generic matrices} representation 
$$\rho^{\rm univ}: \Z\{X\!\} \longrightarrow M_d(F_X(d))$$ defined by $x \mapsto
(x_{i,j})_{i,j}$, hence we get by the usual Amitsur's formula (or by Cor.~\ref{corimage}) a natural homogeneous multiplicative polynomial
law of degree $d$ given by $$\det \circ \rho^{\rm univ}: \Z\{X\!\} \longrightarrow
E_X(d),$$ where $E_X(d) \subset F_X(d)$ is the subring generated by the
coefficients of the characteristic polynomials of the $\rho^{\rm univ}(w)$, $w \in \Z\{X\!\}$. 

\begin{thm}\label{thmvac}(Vaccarino) $\det \circ \rho^{\rm univ}$ induces an isomorphism $\Gamma_\Z^d(\Z\{X\!\})^{\rm ab} \isomo E_X(d)$. In particular, $\Gamma_\Z^d(\Z\{X\!\})^{\rm ab}$ is a free $\Z$-module.
\end{thm}

See \cite[Thm. 6.1]{vac0} and \cite[Thm 28]{vac1}. Set $X=R$ and consider the canonical map $\pi : \Z\{X\!\} \longrightarrow R$. 
Vaccarino's theorem shows that for any determinant $\Det: R \longrightarrow A$ of dimension $d$, there is a unique ring homomorphism $\varphi_X : E_X(d) \longrightarrow
A$ such that for all $w \in \Z\{X\!\}$, 
\begin{equation}\label{vactheorem} \varphi_X(\det(\rho^{\rm univ}(w))) =
\Det(\pi(w)).\end{equation}
More generally, it asserts that $\varphi_X \circ (\det \circ \rho^{\rm univ}) =
\Det \circ \pi$ is an equality of $\Z$-polynomial laws. Via $\varphi_X$, we may view the $A$-module $R$ as an $E_X(d)$-module and $D$ becomes an  $E_X(d)$-polynomial law. We get then a commutative square of $E_X(d)$-polynomial laws : 
$$\xymatrix{ R \ar@{->}[rr]^\Det & & A \\ 
  \Z\{X\!\}\otimes_\Z E_X(d)  \ar@{->}[u]^{\pi \otimes \varphi_X} \ar@{->}[rr]^{\hspace{1cm}\det\circ \rho^{\rm univ}} & & E_X(d) \ar@{->}[u]_{\varphi_X}} $$

Now, all the assertions of Lemma~\ref{formules} follow at once from this
diagram and the classical formulae in matrix rings (here $M_d(F_X(d))$). 
For example for part (iv), the Cayley-Hamilton theorem shows that for
$r_1,\dots,r_n \in X=R$, $\rho^{\rm univ}(\chi_{\alpha}(r_1,\dots,r_n))=0$,
so $$\Det(1+r\chi_{\alpha}(r_1,\dots,r_n))=\det(\rho^{\rm
univ}(1+r\chi_{\alpha}(r_1,\dots,r_n)))=1.$$

\begin{remark}\label{questionun} {\rm \small Actually, part (iv) would follow from
an apparently weaker version of Thm.\ref{thmvac} : For any (finite) set
$X$, $\Gamma_\Z^d(\Z\{X\})^{\rm ab}=\TS_\Z^d(\Z\{X\})^{\rm ab}$ is
torsion free as abelian group\footnote{Indeed, it is enough to show (iv) when
$A=\Gamma_\Z^d(\Z\{X\})^{\rm ab}$, $R=A\{X\}$, and $\Det : R \longrightarrow
A$ is the universal determinant. Fix $\alpha$, $r$, $r_1,\dots,r_n$ as in
the statement and set $x=\chi_\alpha(r_1,\dots,r_n)r$. We showed in the proof above 
that $\Lambda_1(xy)=0$ for all $y \in R$, and in particular that
$\Lambda_1(x^m)=0$ for all $m\geq 1$. By the Newton
relations~(\ref{newtrel}), this implies that $i\Lambda_i(x)=0$ for $i\geq 1$, hence
$\Lambda_i(x)=0$ as $A$ is torsion free.}. Unfortunately, as pointed out to us by Vaccarino, this is 
actually equivalent to Thm.\ref{thmvac} in view of Procesi's results. }
\end{remark}

\subsection{Faithful and Cayley-Hamilton determinants}\label{faithCH} 

Let us first introduce the notion of {\it Kernel} of a polynomial
law. Let $M$ and $N$ be two $A$-modules and $P \in {\cal P}_A(M,N)$. Define $\Ker(P) \subset M$, 
as the subset whose elements are the $x \in M$ such that 
$$\forall B \in {\rm Ob}(\cal
C_A), \, \, \, \forall \, b \in B, \, \, \, \forall m \in M\otimes_A B,\, \,
\, P( x\otimes b +m)=P(m).$$  Equivalently, $x \in
\Ker(P)$ if and only if for any integer $n$ and any $m_1,\dots,m_n \in M$, the element
$P(tx+t_1m_1+\dots+t_nm_n) \in N[t,t_1,\dots,t_n]$ is independent of $t$
({\it i.e.} lies in $N[t_1,\dots,t_n]$).
By definition, $\Ker(P)$ is an $A$-submodule of $M$. We say that $P$ is
{\it faithful} if $\Ker(P)=0$.
\par

\begin{lemma}\label{kergen}\label{kerTeasy}\begin{itemize}
\item[(i)] $\Ker(P)$ is the biggest
$A$-submodule $K \subset M$ such that $P$ admits a factorization $P = \widetilde{P} \circ \pi$
with $\pi$ is the canonical $A$-linear surjection $M \longrightarrow M/K$ and $\widetilde{P} \in \cal
{\cal P}_A(M/K,N)$. \ps
\item[(ii)] $\widetilde{P} : R/\Ker(P) \longrightarrow S$ is faithful.\ps
\item[(iii)] If $B$ is a commutative $A$-algebra,
then $${\rm Im}(\Ker(P)\otimes_A B \rightarrow M\otimes_A B) \subset
\Ker(P \otimes_A B).$$
\end{itemize}
\end{lemma}

\begin{pf} Assertion (iii) follows from the transitivity of tensor product. 
Moreover, it is clear that if $P = \widetilde{P} \circ \pi$ for some
$A$-submodule $K \subset M$ as in the statement, then $K \subset \Ker(P)$. We need to check that for any $A$-submodule $K \subset
\Ker(P)$, $P$ factors through a polynomial map
$\widetilde{P}: M/K \longrightarrow N$. \par
Let $B$ be a commutative $A$-algebra
and consider $K_B:={\rm Im}(K \otimes_A B \longrightarrow
M \otimes_A B)$. Then $(M/K)\otimes_A B \isomo (M\otimes_A B)/K_B$ and 
$K_B \subset \Ker(P\otimes_A B)$ by part (iii). In particular, the map $P_B : M \otimes_A B \longrightarrow N \otimes_A B$ satisfies $P_B(k+m)=P_B(m)$ for
any $M \in M\otimes_A B$ and $k \in K_B$, hence we obtain a well-defined map $\widetilde{P}_B: (M/K) \otimes_A B \longrightarrow N \otimes_A B$ via the formula
\begin{equation}\label{factordeter} \widetilde{P}_B((\pi \otimes_A B)(m))=P_B(m), \, \, \forall m \in M\otimes_A B.\end{equation}
We check at once that the collection of maps $\widetilde{P}_B$ with $B$ varying defines an element $\widetilde{P} \in {\cal P}_A(M/K,N)$. \par
If $K \subset \Ker(P)$ and $P=\widetilde{P} \circ \pi$, it follows from formula (\ref{factordeter}) that $\Ker(\widetilde{P})=\Ker(P)/K$, 
hence (ii).
\end{pf}

Of course, (\ref{factordeter}) shows that if $P$ is homogeneous of some degree $n$ and $P = \widetilde{P}
\circ \pi$ as in the lemma, then so is $\widetilde{P}$.

\begin{lemma}\label{kergen2} Let $R$ and $S$ be two $A$-algebras and $P \in {\cal M}_A^d(R,S)$.
\begin{itemize}\item[(i)] $\Ker(P)=\{ r \in R, \, \, \forall B\,\,, \forall r' \in R\otimes_A B, \, \, 
P(1+rr')=1\} =\{ r \in R, \, \, \forall B\,\,, \forall r' \in R\otimes_A
B, \, \,
P(1+r'r)=1\}$.\ps 
\item[(ii)]  $\Ker(P)$ is a two-sided ideal of $R$, it is proper if $d>0$ and $R \neq 0$. It is the biggest
two-sided ideal $K \subset R$ such that $P$ admits a factorization $P =
\widetilde{P} \circ \pi$  
with $\pi$ is the canonical surjection $R \longrightarrow R/K$
and $\widetilde{P} \in {\cal M}^d_A(R/K,S)$. 
\end{itemize}
\end{lemma}

\begin{pf} Denote by $J_1(P)$ and $J_2(P)$ the two sets on the right in the two equalities in part (i). Let $r \in
\Ker(P)$, $B$ a commutative $A$-algebra, and $m=1+h \in R\otimes_A B$. We want
to show that the elements $P(1+r(1+th))$ and $P(1+(1+th)r)$ of $S\otimes_A B[t]$ are the unit element. As
they are polynomial of degree $d$ in $t$, it is enough to check that this holds
in $S\otimes_A B[t]/(t^{d+1})$. But $1+th$ is invertible in $R\otimes_A B[t]/(t^{d+1})$ thus the
multiplicativity assumption implies that
$$P(1+r(1+th))=P((1+th)^{-1}+r)P(1+th)=P((1+th)^{-1})P(1+th)=P(1)=1,$$
and for the same reason $P(1+(1+th)r)=1$, so $\Ker(P) \subset J_1(P)$, $J_2(P)$. 
The same argument shows conversely that $J_i(P) \subset \Ker(P)$, hence
$\Ker(P)=J_1(P)=J_2(P)$. \ps
By (i), $\Ker(P)$ is a two-sided ideal of $R$. As $P(1)=1$ we have $P(1-t)=(1-t)^d$,
thus $1 \notin \Ker(P)$ if $d>0$. Part (ii) follows from formula (\ref{factordeter}) as in the proof of Lemma~\ref{kergen} (i).
\end{pf}

Observe that Lemma~\ref{kergen2} (i) shows that $$\Ker(P)=\{ r \in R, \, \, \forall B\,\,, \forall r' \in R\otimes_A B, \, \, \forall i \geq 1,\, \, 
\Lambda_i(rr')=0\}.$$ 
Moreover, $r \in \Ker(P)$ if for any $r_1,\dots,r_n \in R$, we have
$$P(1+r(t_1r_1+t_2r_2+\dots,+t_nr_n))=1.$$ 
\noindent When $S=A$ is an infinite domain, then $\Ker(P)=\{ r \in R, \, \, \forall r' \in R, \, \, 
P(1+rr')=1\}$. \ps\ps

\ps\medskip
Assume now that $S=A$, {\it i.e.} that $\Det : R \longrightarrow A$ is a
determinant of dimension $d$. We denote by $\CH(\Det) \subset R$ the two-sided ideal of $R$ generated
by the coefficients of
$$\chi(t_1r_1+\dots+t_nr_n) \in
R[t_1,\dots,t_n],$$ with $r_1,\dots,r_n\in R$, $n \geq 1$ (i.e. by the elements $\chi_\alpha(r_1,\dots,r_n)$ defined in \S~\ref{polid}). We say that $\Det$ is {\it Cayley-Hamilton} if $\CH(\Det)=0$. Equivalently, $\Det$ is Cayley-Hamilton if the polynomial law $\chi : R \longrightarrow R$ is identically zero. In this case, we will say also that $(R,\Det)$ is a {\it Cayley-Hamilton $A$-algebra of degree $d$}. Note that by definition, if $\Det : R \longrightarrow A$ is Cayley-Hamilton and if $B$ is any
commutative $A$-algebra, then $\Det\otimes_A B : R\otimes_A B
\longrightarrow B$ is also Cayley-Hamilton. \ps\ps The Cayley-Hamilton property behaves rather well under several operations, which
is in general not the case of the faithful property. 

\begin{example}\label{exampleCH} {\rm \begin{itemize}
\item[(i)] If $R$ is an Azumaya algebra of rank $d^2$ over $A$ and
$\Det$ is the reduced norm, then $\Det$ is Cayley-Hamilton and faithful. \ps
\item[(ii)] If $D$ is Cayley-Hamilton and $S \subset R$ is any $A$-subalgebra, then the
restriction of $D$ to $S$ is obviously Cayley-Hamilton. However, the analogous
assertion with Cayley-Hamilton replaced by faithful does not hold. For
example, if $T_d(A) \subset M_d(A)$ is the $A$-subalgebra of upper triangular
matrices, then $\det: T_d(A) \longrightarrow A$ is Cayley-Hamilton, but not
faithful. An easy computation shows that $\Ker(\det)$ is the kernel of the
natural diagonal projection $T_d(A) \longrightarrow A^d$ in this case.
\end{itemize}}
\end{example}

\begin{lemma}\label{CHfaith} $\Ker(\Det)$
contains $\CH(\Det)$. In particular, if $\Det$ is faithful then $R$ is
Cayley-Hamilton. 
\end{lemma}

\begin{pf} As $\Ker(\Det)$ is a two-sided ideal by Lemma~\ref{kergen2} (ii),
the first assertion follows from the description of $\Ker(\Det)$ given in
Lemma~\ref{kergen2} (i) and from Lemma~\ref{formules} (iv). The second assertion
follows from the first one. 
\end{pf} 

The next paragraph is a digression about the notion of Cayley-Hamilton representations, the reader urgently interested in the 
proofs of the results stated in the introduction may directly skip to section~\ref{finitenesssesh}.

\subsection{The $\cCH$ category of Cayley-Hamilon representations}\label{CHGcat}
Let us consider the counterpart of these notions on the space $X(G,d)={\rm Spec}(\Z(G,d))$ defined in \S~\ref{defdet}. Consider the tautological (universal) 
determinant of dimension $d$ $$\Det^{u}: \Z(G,d)[G] \longrightarrow \Z(G,d).$$
The {\it universal Cayley-Hamilton algebra} $$R(G,d):=\Z(G,d)[G]/\CH(\Det^{u})$$
is equipped with a natural group homomorphism $\rho^{u} : G \longrightarrow R(G,d)^\ast$. This morphism has the following nice universal property.\ps
Define a {\it Cayley-Hamilton $A$-representation (or $\CH$-representation for short) of $G$ of dimension $d$} as a triple $(A,(R,D),\rho)$ where $A$ is a commutative ring, $(R,D)$ is a Cayley-Hamilton $A$-algebra for the 
determinant $D : R \longrightarrow A$ of dimension $d$, and $\rho : G \longrightarrow R^{*}$ is a group homomorphism. Of course, 
usual representations give rise to $\CH$-representations, but there are many more in general. \ps
Consider the category $\cCH$ 
whose objects are the $\CH$-representations of $G$ of dimension $d$, and with arrows 
$$(A_{1},(R_{1},D_{1}),\rho_{1}) \longrightarrow (A_{2},(R_{2},D_{2}),\rho_{2})$$ 
the pairs $(f,g)$ where $f: A_{1} \longrightarrow A_{2}$ and $g: R_{1} \longrightarrow R_{2}$ are ring homomorphisms such that if $\iota_{i}: A_{i} \longrightarrow R_{i}$ is 
the $A_{i}$-algebra structure on $R_{i}$, then $g \circ \iota_{1} = \iota_{2} \circ f$, $f \circ D_{1}  = D_{2} \circ g$, and $\rho_{2}= g \circ \rho_{1}$.

\begin{prop} $(\Z(G,d),(R(G,d),\Det^{u}),\rho^{u})$ is the initial object of $\cCH$.
\end{prop}

\begin{pf} Let $(A,(R,D),\rho)$ be a $\CH$-representation of $G$ of dimension $d$. The group homomorphism 
$\rho: G \longrightarrow R^*$ is induced by a unique $A$-algebra homomorphism $\widetilde{\rho}: A[G] \longrightarrow R$ and $\Det \circ \widetilde{\rho}$ is then an 
$A$-valued determinant on $G$ of dimension $d$. We get this way a unique ring homomorphism $f : \Z(G,d) \longrightarrow A$, hence a ring homomorphism 
$\Z(G,d)[G] \longrightarrow A[G] \longrightarrow R$. As $(R,D)$ is Cayley-Hamilton, it factors through a ring homomorphism 
$g : R(G,d) \longrightarrow R$, and we check at once that $(f,g)$ has all the required properties.
\end{pf}

The Cayley-Hamilton $\Z(G,d)$-algebra $R(G,d)$ is the global section of a quasi-coherent sheaf of Cayley-Hamilton algebras $\cal R(G,d)$ on $X(G,d)$. Its formation commutes with arbitrary base changes (contrary to the faithful quotient in general) : for any morphism $\Spec(A) \longrightarrow X(G,d)$, corresponding to a determinant $D : A[G] \longrightarrow A$, then the natural surjective map 

\begin{equation}\label{CHrepuniv} A[G]\longrightarrow \cal R(G,d) \otimes_{\Z(G,d)} A\end{equation}
induces an isomorphism $A[G]/\CH(\Det) \isomo \cal R(G,d) \otimes_{\Z(G,d)} A$. \ps \ps

\begin{remark}\label{CHvsrep}(CH-representations versus representations) {\rm In general, given a point $\Spec(A) \longrightarrow X(G,d)$, {\it i.e.} a determinant $D : A[G] \longrightarrow A$, there is no representation $\rho : A[G] \longrightarrow M_d(A)$ such that $D = \det \circ \rho$ (see e.g. \cite[Thm. 1.6.3]{bch}). However, we have a natural candidate for a substitute which is the CH-representation (\ref{CHrepuniv}), {\it i.e.} $$G \longrightarrow (A[G]/\CH(D))^\ast.$$ 
Thus it is an important task to study the sheaf ${\cal R}(G,d)$ of CH-algebras. It turns out to be extremely nice over some specific loci of $X(G,d)$. For instance, we will show in Corollary~\ref{corstructurethm} that it is a sheaf of Azumaya algebras of rank $d^{2}$ over the {\it absolute irreducibility locus} of $X(G,d)$ ; in particular, \'etale-locally on this (open) subspace $D^{\rm univ}$ is the determinant of a true representation (unique, surjective). \par
The situation is more complicated over the reducible locus. In Theorem~\ref{structurethm} we will study more generally the algebra $\cal R(G,d) \otimes \OO_x^{\rm hens}$ when $x \in X(G,d)$ is reducible but {\it multiplicity free} : it is a {\it generalized matrix algebra} in the sense of \cite[\S 1.3]{bch} (and all such algebras occur somehow this way ; when $d! \in A^*$, this result follows from \cite[Thm. 1.4.4]{bch}).}\end{remark}

\begin{remark}\label{embedpr}({\it The embedding problem}) 
{\rm The embeding problem is to decide whether the CH-algebra $(R(G,d),D^{\rm univ})$ admits a CH-embedding in $(M_d(B),\det)$ for some commutative ring $B$. A result of Procesi \cite{Proc2} asserts that it holds after tensoring by $\Q$, but the result over $\Z$ still seems to be open (see \cite{vac2}). The problem is local on $X(G,d)$, and there are some partial known results. For instance, we will show in Theorem~\ref{structurethm} that it holds at $x \in X(G,d)$ (i.e for $R(G,d) \otimes \OO_x$) whenever $x$ is multiplicity free (compare with \cite[Prop. 1.3.13]{bch}).
}\end{remark}

\subsection{Determinants and pseudocharacters}\label{detpseudo} 
We end this paragraph by a comparison between determinants and
pseudocharacters. Let us start with the following result, whose conclusion will be actually sharpened below.

\begin{prop}\label{Qalg} The map $\Det
\mapsto \Tr$ defined in \S~\ref{polid} induces an injection between the set of
$d$-dimensional $A$-valued determinants on $R$ and the set of $d$-dimensional
$A$-valued
pseudocharacters on $R$. When $A$ is a $\Q$-algebra, it is a bijection.

\end{prop}

\begin{pf} Let $\Det : R \longrightarrow A$ be a determinant of dimension $d$. By Lemma~\ref{formules} (iii), $\Tr$ is a
$d$-dimensional pseudocharacter on $R$ (note that $\Tr(1)=d$). Moreover, 
the Newton relations (\ref{newtrel}) show that for each commutative $A$-algebra $B$, each $r \in R
\otimes_A B$ and $i\leq d$, $\Lambda_i(r)$ lies in the $\Z[1/d!]$-algebra
generated by $\Tr(r^j)$ for $j\leq i$, hence $\Tr$ determines
$\Det$. \ps
Let $T$ be a $d$-dimensional $A$-valued pseudocharacter, it remains to show that
it has the form $\Tr$ for some $\Det$. By the Newton relations again, there 
is a unique element $$P \in \Z[1/d!][S_1,...,S_d]$$ such that for any commutative
ring $B$ and $r \in M_d(B)$, we have $P(\dots,\tr(r^i),\dots)=\det(r)$. Of course, if we ask $S_i$
to have degree $i$, then $P$ is homogeneous of degree $d$. Moreover, $P(d,d(d-1)/2,...,d,1)=1$.
We consider then the $A$-polynomial law $\Det : R \longrightarrow A$ defined by
$\Det=P(\dots,T(r^i),\dots)$. It is homogenenous of degree $d$ and satisfies $\Det(1)=1$. By
construction, it is enough to check that $\Det(rr')=\Det(r)\Det(r')$ for all
commutative $A$-algebra $B$ and $r, r' \in R\otimes_A B$. By
construction, $\Det_B(r)=P(\dots,(T\otimes_A B)(r^i),\dots)$ for all $r \in R\otimes_A B$, so we may assume that $A=B$. By a result of Procesi~\cite{Proc2}, there is a
commutative $A$-algebra $C$ with $A \longrightarrow C$ injective and an
$A$-algebra homomorphism $$\rho : R \longrightarrow M_d(C)$$
such that $\tr\circ \rho =T$. But then $\det(\rho(x))=\Det(x)$ is
multiplicative, and we are done. 
\end{pf}

\begin{remark}\label{rmkQd}{\rm The proposition might hold under the weaker assumption $d! \in
A^*$ but we don't know how to prove it in general, namely : we don't know how to show that the obvious $A$-polynomial law of degree $d$ attached to a pseudocharacter 
$T : R \longrightarrow A$ is multiplicative (compare with~\cite[Remark 1.2.9]{bch}). However, using structure theorems for pseudocharacters 
over fields and over local rings instead of \cite{Proc2}, we know that this holds in either of the following situations : \begin{itemize}
\item[(i)] $A$ is reduced, \ps
\item[(ii)] For all $x \in {\rm
Specmax}(A)$, and $k$ an algebraic closure of the residue field at $x$, the
induced pseudocharacter $T\otimes_A k$ is multiplicity free (use~\cite[Prop. 1.3.13]{bch} and \cite[Thm.1.4.4]{bch}).\ps
\end{itemize}
In general, it would be enough (actually equivalent) to know that if $G$ is the free monoid over two letters $\{a,b\}$,
and $T : \Z[G] \longrightarrow A$ is the universal pseudocharacter on $G$ of
dimension $d$ (with $d! \in A^*$), then $A$ (which is easy to describe by generators and
relations) is torsion free over $\Z$. The next result is an evidence for the general case.}
\end{remark}

\begin{prop}\label{casedim2} Assume : \begin{itemize} 
\item[(i)] either that $2$ is invertible in $A$ and $d=2$, 
\item[(ii)] or that $(2d)!$ is invertible in $A$,\end{itemize}
then the map $\Det \mapsto
\Tr$ defined in \S~\ref{polid} induces a bijection between the set of $d$-dimensional $A$-valued
determinants on $R$ and the set of $d$-dimensional $A$-valued
pseudocharacters on $R$. 
\end{prop}

\begin{pf} We first show (i). For $x, y \in R$ set $f(x,y)=T(x)T(y)-T(xy)$ and $D(x)=f(x,x)/2$.
Then $f: R \times R \longrightarrow A$ is 
an $A$-bilinear map and $D : R \longrightarrow A$ is a quadratic $A$-map with associated
bilinear map $f$. In particular, $D$ defines a quadratic $A$-polynomial law
$R \longrightarrow A$ which satisfies $D(1)=1$ (see Example~\ref{examplehomogene} (ii)). We have to check that
$D(xy)=D(x)D(y)$ for all $x, y \in R$. We check at once as in
Example~\ref{exdim2} that
it suffices to show that for all $x,x',y,y' \in R$, we have
\begin{equation}\label{relaprouver}f(xy,x'y')+f(xy',x'y)=f(x,x')f(y,y').\end{equation}

For $m\geq 1$, $\sigma \in \got{S}_m$ and $x=(x_1,\dots,x_m) \in X^m$, set
$T^{\sigma}(x)=T(x_{i_1}...x_{i_r})$ if $x$ is the cycle
$(i_1,\dots,i_r)$, and $T^{\sigma}(x)=\prod_i T^{c_i}(x)$ if
$\sigma=\prod_i
c_i$ is the cycle decomposition of $\sigma$. For example for $m=3$, the
$2$-dimensional pseudocharacter relation reads
$s_3(T):=\sum_{\sigma \in \got{S}_3} \varepsilon(\sigma)T^{\sigma}=0$ on
$R^3$, where $\epsilon$ is   
the signature on $\got{S}_m$. We have to show that this relation implies
(\ref{relaprouver}) if $2$ is invertible    
in $A$. \par \smallskip
Let us fix now $m=4$ and
consider the order $8$ subgroup $H \subset \got{S}_4$ generated by
$H_0=\langle (1,2), (3,4)\rangle$ and $(1324)$. Let $s: H \rightarrow
\{\pm 1\}$ denote the unique character which coincides with the signature
$\varepsilon$ on $H_0$ and such that $s((1324))=1$. Condition
(\ref{relaprouver}) reads 
\begin{equation}\label{equationdet} \forall x \in R^4, \, \, \sum_{h \in H}
s(h)
T^h(x)=0.\end{equation}
Let $B=\Z[1/2][\got{S}_4]$ be the group ring of $\got{S}_4$ over $\Z[1/2]$
and consider the two elements of $B$  
$$p:=\frac{1}{8}\sum_{h \in H} s(h) h, \, \, \, q=\sum_{g \in
\got{S}_3} \varepsilon(g)g$$
where $\got{S}_3$ is viewed as the subgroup of $\got{S}_4$ fixing $\{4\}$.
Note that $p$ is an idempotent of $B$. To prove that the pseudocharacter relation 
implies (\ref{equationdet}), it is enough to show that
$p
\in BqB$ (see e.g. remarks (i) to (v) following Thm. 4.5 of \cite{Proc}).
For that it is
actually enough to show that for any field $k$ in which $2$ is invertible
then $p \in B_k q B_k$,
where $B_k:=B\otimes_\Z k=k[\got{S_4}]$. We fix such a field.
\par\smallskip

 Let $k^4$ be the natural permutation
representation of $\got{S}_4$.
As $2 \in k^*$ we have $k^4=1\oplus \St$, and we check at once that $\St$
absolutely irreducible.
Let $V=\St \otimes \varepsilon$. We have\footnote{Note that the vector
$(1,1,-1,-1)$ (resp. $(1,1,1,-3)$) generates the representation $s\varepsilon$ under
the action of $H$    
(resp. is invariant under the action of $\got{S}_3$). } 
\begin{equation}\label{decomposition} {\rm Ind}_H^{\got{S}_4} s = V,\, \, \,
{\rm Ind}_{\got{S}_3}^{\got{S}_4} \varepsilon = \varepsilon \oplus
V.\end{equation}
As $|\got{S}_4|/\dim(V)=8 \in k^*$, $V$ is a projective $B_k$-module and we
may find a central idempotent $e \in B_k$ acting on $V$ as the identity, and
as $0$ on the $k$-representations of $\got{S}_4$ not containing $V$.
Moreover, the $k$-algebra $e B_k$ (with unit $e$) is isomorphic to 
$\End_k(V) \simeq M_3(k)$ as $V$ is absolutely irreducible. As ${\rm
Ind}_H^{\got{S}_4} s = V$, the idempotent $p$ acts non-trivially in a
$k$-representation $U$ of $\got{S}_4$ if and only if $U$ contains $V$.
Applying this to $U=B_k(1-e)$ we obtain $p(1-e)=0$, so $p \in eB_k$. 
But one easily sees that $q(V)\neq 0$, for instance $q \cdot (1,0,0,-1)=(2,2,2,-6)$. It follows that $eq \neq 0$ so 
$B_k q B_k \supset  B_k eq B_k = e B_k$ by simplicity of $e B_k$, thus $p
\in B_k q B_k$. \par \smallskip
  The second statement is actually a formal consequence of Procesi's results \cite{Proc}. Let us consider 
  the full polarization of the polynomial map $\det(g)\det(h)-\det(gh)$ on $M_{2d}^{2}$, it is given by some element $p \in \Z[\got{S}_{2d}]$ 
  (see below for an explicit formula of a partial polarization), and as above we have to show that $p \in BqB$ where 
  $q = \sum_{\sigma \in \got{S}_{{d+1}}} \varepsilon(\sigma)\sigma$ and $B=\Z[1/(2d!)][\got{S}_{2d}]$. By the second fundamental theorem 
  of invariants of set of matrices \cite{Proc}, we know that this holds over $\Q$, so $mp \in BqB$ for some $m \in \Z$. As $B$ is isomorphic to a direct product 
  of matrix rings over $\Z[1/(2d)!]$, and as $\frac{q}{(d+1)!}$ is an idempotent of $B$, it turns out that $B/BqB$ is torsion free, and we are done.

\end{pf}
\smallskip

We end this paragraph by giving an explicit {\it $(d,d)$-partial} polarization of the homogeneous (of degree $2d$) 
polynomial map 
\begin{equation}\label{detprodmap} (g,h) \mapsto \det(g)\det(h)-\det(gh), \, \,  M_{d}(A)^{2} \rightarrow A \end{equation}
when $d! \in A$, which extends the relation (\ref{equationdet}) obtained in dimension $2$.  By this we mean an $A$-multilinear map 
$\varphi : M_{d}(A)^{2d} \rightarrow A$ which is symmetric only in the first $d$ (resp. last $d$) variables, and such that 
$\varphi(g,g,\dots,g,h,h,\dots,h)=(d!)^{2}(\det(g)\det(h)-\det(gh))$ for any $(g,h) \in M_{d}(A)^{2}$.
\par \smallskip

Let $H_{0} \subset \got{S}_{2d}$ be the subgroup preserving $\{1,\dots,d\}$ (thus $H_{0} \simeq \got{S}_{d}^{2}$) ;  the 
element $$\tau = \prod_{i=1}^{d}(i\, \, d+i)$$ has order $2$ and normalizes $H_{0}$, thus $H=\langle H_{0}, \tau\rangle$ is a subgroup of order $2(d!)^{2} \in A^{*}$. The signature on $H_{0}$ being $\tau$-invariant, 
there exists a unique character 
$s : H \rightarrow \{\pm 1\}$ such that $s(\tau)=-1$ and such that $s$ coincides with the signature on $H_{0}$. We define an $A$-multilinear map 
$\varphi : M_{d}(A)^{2d} \rightarrow A$ by $$\varphi_{T}=\sum_{\sigma \in H} s(\sigma)T^{\sigma},$$ where $T$ is the trace.

\par \smallskip 

\begin{prop}\label{bipoldet} $\varphi_{T}$ is a $(d,d)$-partial polarization of  (\ref{detprodmap}). \end{prop}

\begin{pf} Note that $\varphi_{T}$ is $H_{0}$-invariant by construction. Let us first consider the multilinear invariant map 
$\psi : M_{d}(A) \rightarrow A$ associated to the element $u=\sum_{\sigma \in \got{S}_{d}} 
\epsilon(\sigma) \sigma \in \Z[\got{S}_{d}]$, that is $\psi=\sum_{\sigma \in \got{S}_{d}} \epsilon(\sigma) T^{\sigma}$. As is well-known, 
$\psi$ is the full polarization of $\det$, being the trace of $(g_{1},\dots,g_{n})$ on $\frac{u}{d!}(V^{\otimes_{A} d})=\Lambda^{d}(V)$, where 
$V=A^{d}$. We deduce from this an expression for a partial polarization of $(g,h) \mapsto \det(gh)$ by $(d,d)$-polarizing each term of the form 
$T^{\sigma}(g_{1}h,g_{2}h,...,g_{d}h)$ as
\begin{equation}\label{polpartielle}\sum_{\sigma' \in \got{S}_{d}} T^{\sigma}(g_{1}h_{\sigma'(1)},g_{2}h_{\sigma'(2)},\dots,g_{d}h_{\sigma'(d)}).
\end{equation}
It only remains to identify the associated elements of $\Z[\got{S}_{2d}]$. \par 
Writes $H_{0}=H_{1} . H_{2}$ where $H_{1}$ is the subgroup fixing 
$d+1$ to $2d$ 
and $H_{2}$ fixes $1$ to $d$, and identify $H_{2}$ with $\got{S}_{d}$ under the bijection $\{1,\dots,d\} \isomo \{d+1,\dots,2d\}$, $i \mapsto i+d$. 
A simple cycle computation shows that 

$$T^{\sigma}(g_{1}h_{\sigma'(1)},g_{2}h_{\sigma'(2)},\dots,g_{d}h_{\sigma'(d)})=T^{\sigma''}(g_{1},g_{2},\dots,g_{d},h_{1},h_{2},\dots,h_{d})$$
where $\sigma''=\sigma\sigma' \tau {\sigma'}^{-1} \in \got{S}_{2d}$. The key fact is that for 
$(i_{1}\, i_{2}\, \cdots i_{r}) \in H_{1}$ any cycle, and for $j_{1}$, $j_{2}$, \dots, $j_{r}$ any 
distinct elements in $\{d+1,\dots,2d\}$, then  

$$ (i_{1}\, i_{2}\, \cdots i_{r}) \,(i_{1}\, j_{1})\, (i_{2},\, j_{2})\, \cdots (i_{r},\, j_{r})\, = \, (i_{1}\, j_{1}\,i_{2}\,j_{2}\,\cdots i_{r}\, j_{r}). $$

As a consequence, we get a $(d,d)$-polarization of (\ref{detprodmap}) as the multilinear invariant associated to the element 
$$p=\sum_{\sigma \in H_{0}} \varepsilon(\sigma) \sigma - \sum_{(\sigma,\sigma') \in H_{1}\times H_{2}} 
\epsilon(\sigma)\sigma \sigma'\tau{\sigma'}^{-1} \in \Z[\got{S}_{2d}].$$
A simple change of variables $(\sigma_{1},\sigma_{2})=(\sigma\tau\sigma'\tau,{\sigma'}^{{-1}})$ identifies this map with $\varphi_{T}$.
\end{pf}

\section{Structure and finiteness theorems}\label{finitenesssesh}

In this section, we will give some necessary conditions ensuring that a
determinant $\Det : R \longrightarrow A$ of dimension $d$ 
is the determinant of a true representation $R \longrightarrow M_d(A)$. As explained in Remark~\ref{CHvsrep}, we will get these results by first proving some structure theorems for certain Cayley-Hamilton algebras. \par 
In an independent last paragraph, we will also state and discuss 
a result of Vaccarino and Donkin asserting in particular that $\Z(G,d)$ is
finite type over $\Z$ when the group (or monoid) $G$ is finitely generated.

\subsection{Some preliminary lemmas} \label{prellemmas}

Let $S_1$ and $S_2$ be two $A$-algebras, $B$ a commutative $A$-algebra, $d$ an integer, and 
let $p_j : S_j\longrightarrow B$ be a multiplicative $A$-polynomial
law which is homogeneous of degree $d_i$, with $d_1+d_2=d$. Then we check at
once that the $A$-polynomial map $$p_1p_2: (x_1,x_2) \mapsto p_1(x_1)p_2(x_2), \, \, \, 
S_1 \times S_2 \longrightarrow B,$$ which is
homogeneous of degree $d_1+d_2$, is again multiplicative. We will call
$p_1p_2$ the product of $p_1$ and $p_2$. This operation induces a natural
$A$-algebra homomorphism 
\begin{equation}\label{prodformule}\Gamma_A^d(S)^{\rm ab} \longrightarrow \prod_{i=0}^d \Gamma_A^i(S_1)^{\rm ab}
\otimes_A \Gamma_A^{d-i}(S_2)^{\rm ab}.\end{equation}

Recall that an $A$-algebra is called finite diagonal if it is isomorphic to $A^n$ (with coordinate-wise addition and multiplication)  for some integer $n\geq 1$.

\begin{lemma}\label{prodlemma}\begin{itemize}
\item[(i)] The map (\ref{prodformule}) is an $A$-algebra isomorphism. \ps
\item[(ii)] If $S$ is a finite diagonal $A$-algebra, then so is
$\Gamma_A^d(S)=\Gamma_A^d(S)^{\rm ab}$.\ps
\item[(iii)] Assume ${\rm Spec}(B)$ is connected and $B\neq 0$. Then any multiplicative
homogeneous $A$-polynomial law $S_1
\times S_2 \longrightarrow B$ of degree $d$ is the product $p_1p_2$ of two unique multiplicative
homogeneous $A$-polynomial laws $p_i : S_i \longrightarrow B$ with degree
$d_i$, and we have $d_1+d_2=d$.
\end{itemize}
\end{lemma}

\begin{pf}  Note that if $B$ is any $A$-algebra (non necessarily commutative), and if $p_{j} : S_{j} \longrightarrow B$ 
are multiplicative $A$-polynomial laws of degree $d_{j}$ such that the images of $p_{1}$ and $p_{2}$ commute in the obvious sense, then 
$(p_{1} p_{2})(x_{1},x_{2}):=p_1(x_{1})p_{2}(x_{2})$ still defines an $A$-multiplicative polynomial law of degree $d_{1}+d_{2}$. This defines 
a natural $A$-algebra homomorphism 
\begin{equation}\label{prodformulenab}\Gamma_A^d(S) \longrightarrow \prod_{i=0}^d \Gamma_A^i(S_1)
\otimes_A \Gamma_A^{d-i}(S_2).\end{equation}
of which (\ref{prodformule}) results by abelianization, thus it is enough to check that (\ref{prodformulenab}) is an isomorphism. By definition, 
the projection of (\ref{prodformulenab}) to the $i$-th factor corresponds to the homogeneous $A$-polynomial law of degree $d$ 
$$S_{1} \times S_{2} \longrightarrow 
\Gamma_{A}^{i}(S_{1})\otimes_{A} \Gamma_{A}^{{d-i}}(S_{2}),\, \, \, (s_{1},s_{2}) \mapsto {s_{1}}^{[i]}\otimes s_{2}^{[d-i]},$$
(which is incidentally obviously multiplicative) hence is exactly the map defined more generally by Roby in \cite[\S 9]{roby} for any pair 
of $A$-modules $(S_{1},S_{2})$, and which is an $A$-linear isomorphism by \cite[Thm. III.4]{roby}, which proves (i). More precisely, we showed that as $A$-algebras there is an isomorphism
\begin{equation}\label{prodform}\Gamma_A^d(S_{1} \times S_{2}) \isomo \prod_{i=0}^d \Gamma_A^i(S_1)
\otimes_A \Gamma_A^{d-i}(S_2).\end{equation}
\ps
Note that $\Gamma_A^i(A) = A\cdot 1^{[i]} \simeq A$
for each $i\geq 0$. In particular, if $S_1=A$ then (\ref{prodform}) shows that
$\Gamma_A^d(S) \isomo \prod_{i=0}^d \Gamma_A^i(S_2)$ as $A$-algebras. Part (ii) follows then
by induction. \ps
We now show assertion (iii). It follows from (i) and the following general fact. Consider a finite number of rings with unit $C_1$, \dots, $C_m$, and set $C=\prod_{i=1}^m C_i$. Let $B$ be a nonzero commutative ring with unit, with connected spectrum, and let $f : C \rightarrow B$ be a ring homomorphism. Then $f$ factors through the projection $C \rightarrow  C_j$ for a unique $j$. Indeed, let $c_i$ be the central idempotent of $C$ whose $j$-th component is $0$ if $i \neq j$, and the unit of $C_i$ if $j=i$. Set $e_i=f(c_i)$. Then $\{e_i, 1\leq i \leq m\} \subset B$ is a set of $m$ idempotents of $B$ such that $\sum_i e_i=1$ and $e_ie_j=0$ if $i \neq j$. As $\Spec(B)$ is connected, it follows that $e_i=1$ or $0$ for each $i$, exclusively as $B\neq 0$, and that there is a necessarily unique $j$ such that $e_j=1$, and the claim follows. 
\end{pf}

\begin{example} {\rm Assume that $G = G_1 \times G_2$, then $\Z(G,d) \isomo
\prod_{i+j=d} \Z(G_1,i)\otimes_\Z \Z(G_2,j)$.}
\end{example}

Let $R$ be an $A$-algebra. Recall that an element $e \in R$ is said to be idempotent if $e^2=e$, in which case $1-e$
is also idempotent. The subset $eRe \subset R$ is then
an $A$-algebra whose unit element is $e$ and $eRe \oplus (1-e)R(1-e)$ is an
$A$-subalgebra of $R$ isomorphic to $eRe \times (1-e)R(1-e)$. We say that a family of idempotents $\{e_i\}$ is
orthogonal if $e_ie_j=0$ if $i\neq j$. Let $\Det : R \longrightarrow A$ be a
determinant of dimension $d$ and assume $A \neq 0$.

\begin{lemma} \label{idempotents} Assume that $\spec(A)$ is connected and let
$e \in R$ be an idempotent.
\begin{itemize}
\item[(1)]  The polynomial map $\Det_e : eRe \longrightarrow A$, $x \mapsto
\Det(x+1-e)$, is a determinant whose dimension $r(e)$ is $\leq d$. \ps
\item[(2)] We have $r(1-e)+r(e)=d$. Moreover, the restriction of $\Det$ to
the $A$-subalgebra $eRe\oplus (1-e)R(1-e)$ is the product determinant
$\Det_e \Det_{1-e}$.\ps
\item[(3)] If $\Det$ is Cayley-Hamilton (resp. faithful), then so is
$\Det_e$. \ps
\item[(4)] Assume that $\Det$ is Cayley-Hamilton. Then $e=1$ (resp. $e=0$) if and only if
$\Det(e)=1$ (resp. $r(e)=0$). Let $e_1, \dots, e_s$ be a family of (nonzero) 
orthogonal idempotents of $R$. Then $s\leq d$, and we have the inequality $\sum_{i=1}^s r(e_i)\leq
d$, which is an equality if and only if $e_1+e_2+\dots+e_s=1$.
\end{itemize}
\end{lemma}

\begin{pf} Set $S_1=eRe$, $S_2=(1-e)R(1-e)$, and consider the $A$-subalgebra $S=S_1\oplus S_2 \subset
R$. Then $e$ is a central idempotent in $S$, hence the map $x \mapsto (ex,(1-e)x)$ is an $A$-algebra isomorphism $S
\isomo S_1 \times S_2$. Lemma~\ref{prodlemma} (iii) applied to the
restriction of $\Det$ to $S$ shows parts (1) and (2). \ps
Assume that $\Det$ is faithful. Let $x \in \Ker(\Det_e)$, $B$ a
commutative $A$-algebra and $y \in R\otimes_A B$. Note that 
\begin{equation} eRe \otimes_A
B= e(R\otimes_A B)e\end{equation} is a direct summand of $R\otimes_A B$. We have (using Lemma~\ref{idpoly}) 
$$\Det(1+xy)=\Det(1+exey)=\Det(1+eyex)=\Det(1-e+e+eyex)=\Det_e(e+eyex)=1,$$
so $x \in \Ker(\Det)$, and $x=0$. \ps
Assume that $\Det$ is Cayley-Hamilton. If $x \in R$, then $\chi^\Det(x,x)=0$. For $x \in
eRe\oplus(1-e)R(1-e)$, we know from part (2) that 
\begin{equation}\label{polye}\chi^\Det(x,t)=\chi^{\Det_e}(ex,t)\chi^{\Det_{1-e}}((1-e)x,t).\end{equation}
For $x \in eRe$, apply the Cayley-Hamilton identity to $x$ and $x+1-e$. We
get that $$P(x)x^{d_2}=P(x)(x-1)^{d_2}=0$$ in $R$, where $P=\chi^{\Det_e}(ex,t)
\in A[t]$ is the characteristic polynomial of $x$ in $eRe$ with respect to
$\Det_e$. But the ideal of $A[t]$ generated by $t^{d_2}$ and $(t-1)^{d_2}$
is $A[t]$, so $P(x)=0$. Applying this argument to $R \otimes_A B$  for all commutative
$A$-algebras $B$, we get that $\Det_e$ is Cayley-Hamilton. \ps
Let us show assertion (4). If $e^2=e$, then $\chi(e,e)-\Det(e) \in Ae \subset
R$. If moreover $\Det$ is Cayley-Hamilton and $\Det(e)=1$, then $ae=1$ in $R$ for some $a \in A$, thus
$1=ae=ae^2=e$. If $r(e)=0$, then $\Det((1-e)+x)$ is a determinant of degree
$0$ on $eRe$, so it is constant and equal to $1$. But then $\Det(1-e)=1$ and
$e=0$ by the previous case. For the last property, set $e_{s+1}:=1-(e_1+\dots+e_s)$. Note that $r(e_i)\leq d$ for each
$i\leq {s+1}$ and $\sum_{i=1}^{s+1} r(e_i)=d$ applying part (2) $s$ times.
We conclude as $r(e_j)=0$ implies $e_j=0$. 
\end{pf}

\begin{exo}\label{zipliesbis} {\small (Another proof or Ziplies's result
\cite{ziplies}) {\rm Let $R$ be an Azumaya algebra of rank $n^2$ over $A$,
show as follows that the reduced norm $N$ induces an isomorphism $$\varphi :
\Gamma_A^d(R)^{\rm ab}\isomo \Gamma^{d/n}_A(A)=A$$ if $n$ divides $d$, and
that $\Gamma_A^d(R)^{\rm ab}=0$ otherwise, using only Lemmas~\ref{prodlemma} and~\ref{idpoly} (i) and
(ii).  Using a faithfully flat commutative $A$-algebra $C$ such that
$R\otimes_A C \isomo M_n(C)$, we may assume $R=M_n(A)$, in which case
$N=\det$.  Let $E_{i,j}$ be the usual $A$-basis of $R$ and $\Det : R
\longrightarrow B$ be any homogeneous multiplicative $\Z$-polynomial law of
degree $d$ ($B$ a commutative $A$-algebra).  Using Lemma~\ref{prodlemma} and
the fact that the $E_{i,i}$ are conjugate under $\GL_n(A)$, show that $n$
divides $d$ and that $\Det_e : E_{1,1}A \longrightarrow B$ is a homogeneous
multiplicative $A$-polynomial law of degree $d/n$.  This provides an
$A$-algebra morphism $$ \iota : A=\Gamma^{d/n}_A(A) \longrightarrow
\Gamma_A^d(M_n(A))^{\rm ab}$$ such that $\varphi \circ \iota={\rm id}$.  To
conclude, it is enough to check that that $\Det_e$ determines $\Det$
uniquely.  Note that $\Det(1-tE_{i,i})=\Det(1-tE_{1,1})=(1-t)^{d/n}$ and for
$i\neq j$,
$$\Det(1-tE_{i,j})=\Det(1-tE_{i,i}E_{i,j})=\Det(1-tE_{i,j}E_{i,i})=1.$$ and
conclude using Amitsur's formula.}} \end{exo}

\medskip
\begin{lemma}\label{CHdim1} Let $\Det : R \longrightarrow A$ be a Cayley-Hamilton
determinant of dimension $1$. Then $R=A$ and $\Det$ is the identity.
\end{lemma}

\begin{pf} By assumption $x=\Tr(x)=\Det(x)$ for all $x \in R$, so the $A$-linear map
$\Tr=\Det : R \longrightarrow A$ is an $A$-algebra isomorphism.
\end{pf}

We now study the Jacobson radical (denoted by $\Rad$) of an algebra with a determinant. We shall need the Nagata-Higman theorem \cite{higman}, that we recall now. 
Let $d$ be an integer and let $k$ be a field such that either ${\rm char}(k)=0$ or ${\rm char}(k)>d$. Let  $R$ be an algebra without unit
over $k$, and assume that $x^d=0$ for all $x \in R$. Then there is an integer $N(d)\leq 2^d-1$ (independent of $R$) such that for all $x_1,\dots,x_{N(d)}$ in $R$, the product 
$x_1 \dots x_{N(d)}$ vanishes.

\begin{lemma}\label{rad1} Assume that $\Det : R \longrightarrow A$ is
Cayley-Hamilton of dimension $d$. \ps\begin{itemize}
\item[(i)] $\Rad(R)$ is the largest  
two-sided ideal $J \subset R$ such that $\Det(1+J) \subset A^*$,\ps
\item[(ii)] $\Ker(\Det) \subset \Rad(R)$,\ps

Assume from now on that $A$ is a field. \ps
\item[(iii)] For all $x \in \Ker(\Det)$ we have $x^d=0$. In particular, if $d!$ is invertible in $A$ then $\Ker(\Det)^{N(d)}=0$. \ps
\item[(iv)] $x \in R$ is nilpotent if and only if $\Det(t-x)=t^d$. Morever, $\rad(R)$ consists of nilpotent elements.\ps
\item[(v)] If $J \subset R$ is a two-sided ideal such that $J^n=0$ for some $n\geq 1$, then $J \subset \Ker(\Det)$ (here it is not necessary to assume that $\Det$ is Cayley-Hamilton).
\end{itemize}
\end{lemma} 
\begin{pf} By the Cayley-Hamilton identity, if $x \in R$, then 
$x$ is invertible in $R$ if and only if $\Det(x)$ is invertible in $A$,
hence (i). Assertion (ii) follows as $\Det(1+\Ker(\Det))=1$.  \par
	Assume that $A=k$ is a field. If $x \in \Ker(\Det)$, then $\chi(x,t)=t^d$ thus $x^d=0$ as $\Det$ is Cayley-Hamilton. When $d!$ is invertible in $k$, 
	the Nagata-Higman theorem applies and proves (iii). If $x \in R$ is nilpotent, then $1+tx$ is
invertible in $R$, hence $\Det(1+tx)$ is invertible in $k[t]$, so
$\Det(1+tx)=1$. The converse follows from the Cayley-Hamilton identity, which even shows that $x^d=0$. Assume that $x \in \Rad(R)$. For all $y \in k[x]$, $1+yx$ is invertible in $R$, so $\Det(1+yx) \in k^*$ and the Cayley-Hamilton identity implies that $1+yx$ is actually invertible in $k[x]$. In particular, $x \in {\rm Rad}(k[x])$. This implies that $x$ is nilpotent as $k[x]$ is a finite dimensional $k$-algebra, hence (iv) follows.\ps
Let $J \subset R$ be as in (v) and $x \in J$. If $y \in R[t_1,\dots,t_n]$, we see that $xy$ is nilpotent, so $\Det(1+txy) \in k[t_1,\dots,t_n,t]$ is invertible, hence constant equal to $1$, and $x \in \Ker(\Det)$.
\end{pf}

The following lemma strenghten part (iv) of the previous one.

\begin{lemma}\label{rad2} Assume that $k$ is a field and let
$\Det : R \longrightarrow k$ be a determinant 
of dimension $d$.\ps
 \begin{itemize}
\item[(i)] If $K/k$ is a separable algebraic extension, then the natural
injection $R \longrightarrow R \otimes_k K$ induces isomorphisms
$\Rad(R)\otimes_k K \isomo \Rad(R\otimes_k K)$ and $\Ker(\Det)\otimes_k K
\isomo \Ker(\Det \otimes_k K)$.\ps
\item[(ii)] Assume that $\Det$ is Cayley-Hamilton. Then $\ker(\Det)=\Rad(R)$.\ps
\end{itemize}
\end{lemma}

\begin{pf} We first check (i). The assertion concerning the Jacobson radical
is well-known. Moreover, the injection of the statement
induces an injection (Lemma~\ref{kerTeasy})$$\Ker(\Det)\otimes_k K
\longrightarrow \Ker(\Det \otimes_k K),$$
and it only remains to check its surjectivity. Enlarging $K$ if necessary,
we may assume that $K/k$ is normal. Let $G:={\rm Gal}(K/k)$ acts semilinearily on $R \otimes_k
K$. By Galois descent, each $G$-stable $K$-subvector space $V$ of $R\otimes_k
K$, has the form $V^G \otimes_k K$ where $V^G \subset R$ is $k$-vector space
of fixed points. We claim that $\Ker(\Det \otimes_k K)$ is
$G$-stable. Indeed, if we let $G$ act on $K[t_1,\dots,t_s]$ by
$\sigma(\sum_\alpha a_\alpha t^\alpha)=\sum_\alpha
\sigma(a_\alpha)t^\alpha$, and then $K[t_1,\dots,t_s]$-semilinearily on $R\otimes_k K[t_1,\dots,t_s]$, then for any $r \in R\otimes_k K[t_1,\dots,t_s]$ we
have
$$\D(\sigma(r))=\sigma(D(r)),$$
from which the claim follows at once. For the same reason, we see that $\Ker(\Det \otimes_k
K)^G \subset \Ker(\Det)$, which concludes the proof.\ps
We now prove assertion (ii). We already know from Lemma~\ref{rad1} (ii) that $\Ker(\Det) \subset \Rad(R)$. By extending the scalars to a separable
algebraic closure of $k$ and part (i), we may assume that $k$ is infinite.
In this case (see \S~\ref{faithCH}), 
\begin{equation}\label{kerinfini}\Ker(\Det)=\{x \in R, \, \, \forall y \in R,
\Det(1+xy)=1\}.
\end{equation} By Lemma~\ref{rad1}, $\rad(R)$ is a two-sided ideal of $R$ consisting of
nilpotent elements $x$, for which $\Det(1+x)=1$, hence (\ref{kerinfini})
implies that $\rad(R) \subset \ker(\Det)$. 
\end{pf}
\ps
\begin{example}\label{exinsep} {\rm In part (i) above it is necessary to assume that $K/k$ is separable. Indeed, let 
$k$ be a field of characteristic $p>0$, $K/k$ a purely inseparable extension of $k$ 
such that for some $p^{\rm th}$-power $q\geq 1$, $x^q \in k$ for all $x \in K$.
Then the $k$-polynomial map $F^q: K \longrightarrow k$, defined by $$x \mapsto x^q$$ for any $x \in K\otimes_k B$ with $B$ a commutative $k$-algebra, is a determinant of dimension $q$, necessarily faithful as $K$ is a field. However, $K \otimes_k K$ is not reduced when $q>1$, in which case $F^q \otimes_k K$ is not faithful by Lemma~\ref{rad1} (v).}
\end{example}
\ps
In what follows, $A$ is a local ring with maximal ideal $m$ and residue field $k:=A/m$. We will
denote by $\Rb$ the
$k$-algebra $R\otimes_A k= R/mR$, and by
$\Detb$ the induced determinant $\Det \otimes k : \Rb \longrightarrow k$.
\begin{lemma}\label{radical} Assume that $\Det$ is Cayley-Hamilton. \begin{itemize}\ps
\item[(i)]The kernel of the canonical
surjection $R \longrightarrow  \Rb/\ker \Detb$ is $\rad(R)$. \ps
\item[(ii)] If $m^{s}=0$ for $s\geq 1$ an integer, and if $d!$ is invertible in $A$, then $\rad(R)^{N(d)s}=0$.\ps
\end{itemize}
\end{lemma}

\begin{pf} Let $J$ be the two-sided ideal of the statement (i), we check first that $J \subset {\rm rad}(R)$. 
It is enough to check that $1+J \subset R^*$, {\it i.e. } $\Det(1+J) \in
A^*$, but this obvious as $\Det(1+J) \in 1+m$ by definition. In particular, $m R \subset {\rm rad}(R)$, hence to check the converse we may 
(and do) assume that $A=k$ and even that $\Det$ is faithful. But then ${\rm rad}(R)=0$ by Lemma~\ref{rad2} (ii). \par
 Assume that $m^{s}=0$. Replacing $R$ by $R/mR$ if necessary, we may assume that $A=k$ is a field and we have to show that $\rad(R)^{N(d)}=0$. Here $N(d)$ is the integer coming from the Nagata-Higman theorem, in particular $N(d)\leq 2^d-1$. We conclude by the equality $\Rad(R)=\Ker(\Det)$ and Lemma~\ref{rad1}.
\end{pf}

\medskip

\subsection{Determinants over a field} In all this paragraph, $k$ is a field and $R$ a $k$-algebra. We fix $\overline{k}$ an algebraic closure of $k$, and by $k^{\rm sep} \subset \overline{k}$ a separable algebraic closure of $k$.

\begin{thm}\label{thmcorps} Assume that $k$ is algebraically closed. For any
$d$-dimensional determinant $\Det : R \longrightarrow k$, there exists a
semisimple representation $\rho : R \longrightarrow M_d(k)$ such that $\Det = \det
\circ \rho$. \ps
Such a $\rho$ is unique up to isomorphism, and $\Ker \rho =
\Ker(\Det)$.
\end{thm}

\begin{cor} Let $G$ be a group (or a unital monoid), then for any
algebraically closed field $k$, $X(G,d)(k)$ is in natural bijection with the
isomorphism classes of $d$-dimensional semisimple $k$-linear representations
of $G$.
\end{cor}

Let us prove the first part of the theorem. By replacing $R$ by $R/\ker(\Det)$ if necessary, we may assume that $\Det$
is faithful. By Lemmas~\ref{idempotents}
and~\ref{radical}, $R$ satisfies the assumptions of the 
following general fact from classical
noncommutative ring theory\footnote{This result is presumably well-known ; it is close to some old results of Kaplanski (see \cite[Chap. 6.3]{herstein}, as well as \cite{Proc0} for a related use). We have learnt most of it from
Rouquier~\cite[Lemme 4.1]{rou}.}.

\begin{lemma}\label{superlemme} Let $k$ be a field, $R$ a $k$-algebra with
trivial Jacobson radical, and $n\geq 1$ an integer. Assume that each element of
$R$ (resp. of $R\otimes_k k^{\rm sep}$) is algebraic over $k$ (resp. $k^{\rm sep}$) of degree less than $n$,  and that the length of families of orthogonal idempotents of $R\otimes_k k^{\rm sep}$ is
also bounded by $n$. Then 
$$R \isomo \prod_{i=1}^s M_{n_i}(E_i)$$ 
where $E_i$ is a division $k$-algebra which is finite dimensional over its center $k_i$. Moreover, each $k_i$ is a finite separable extension of $k$, unless maybe $k$ has characteristic $p>0$, in which case $k[k_i^q]$ is separable over $k$ where $q$ is the biggest power of $p$ less than $n$. \ps
In particular, $R$ is semisimple. It is finite dimensional over $k$ in each of the following three cases : $k$ is a perfect field, or $k$ has characteristic $p>0$ and $[k:k^p] < \infty$, or $p>n$.
\end{lemma}

\begin{pf} Let $A$ be a commutative $k$-algebra such that each element of $A$ is algebraic over $k$ of degree less than $n$, and that the length of families of orthogonal idempotents of $A$ is
also bounded by $n$. If $k$ has characteristic $p>0$ we define $q$ as in the statement, and we set $q=1$ else. Then we check at once that there is a $k$-algebra isomorphism $$A \isomo \prod_{i=1}^r A_i$$
where $r\leq n$ and where $A_i$ is a field whose maximal separable $k$-subextension $A_i^{\rm et}$ is finite dimensional over $k$ (with dimension $\leq n$), and satisfies $A_i^q \subset A_i^{\rm et}$. These facts apply in particular to the center $Z$ of $R$. We get moreover that $\dim_k(Z) <\infty$ in the three cases discussed in the last assertion on the statement.
\par
We prove now that $R$ is semisimple. Let $M$ be a simple $R$-module, $E$ the division $k$-algebra
$\End_R(M)$. We claim first that $M$ is finite dimensional over $E$. Indeed, by Jacobson's density theorem, we know that either $M$ is finite dimensional over $E$ and $R
\longrightarrow \End_E(M) \simeq M_s(E^{\rm opp})$ is surjective, or for each $r\geq 1$
there is a $k$-subalgebra $S_r \subset R$ and a surjective $k$-algebra
homomorphism $S_r \longrightarrow M_r(E^{\rm opp})$, but this second possibility is impossible as elements of $R$ (hence of $S_r$) are algebraic over $k$ of bounded degree by assumption.\par
We claim now that there are only finitely many (pairwise non isomorphic) simple
$R$-modules $M_1, \dots, M_s$, which will conclude the proof. Indeed, assuming this claim and as ${\rm Rad}(R)=0$, there is a natural injective homomorphism 
\begin{equation}\label{weder}R \longrightarrow \prod_{i=1}^s
\M_{n_i}(E_i^{\rm opp}),\end{equation}
$E_i=\End_R(M_i)$, which is surjective as the $M_i$ are pairwise
non isomorphic and simple, hence (\ref{weder}) is an isomorphism. It remains to check the
claim. If $M_1,\dots,M_s$ are any pairwise non isomorphic simple
$R$-modules, and $E_i=\End_R(M_i)$, the morphism (\ref{weder}) is still
surjective. As ${\rm Rad}(R)=0$ and $R$ is algebraic over $k$, the family of orthogonal
idempotents of the right hand side lifts in $R$ 
(\cite[\S 4, ex.5(b)]{bki}), hence $s\leq n$ by assumption, and we are done. \ps
It only remains to show that $R$ is a finite type $Z$-module. As $Z\otimes_k k^{\rm sep}$ is faithfully flat over $Z$, it is enough to check that $R\otimes_k k^{\rm sep}$ is a finite type $Z\otimes_k k^{\rm sep}$-module. Note that the $k^{\rm sep}$-algebra $R\otimes_k k^{\rm sep}$ satisfies also the assumptions of the lemma, hence is semisimple by what we proved till now. Moreover, its center is easily seen to be $Z \otimes_k k^{\rm sep}$. By the Wedderburn-Artin theorem, $$R\otimes_k k^{\rm sep} \isomo \prod_{j=1}^t M_{d_j}(k_j)$$
where $k_j$ is a division $k^{\rm sep}$-algebra, which is moreover algebraic over $k^{\rm sep}$ here. The Jacobson-Noether theorem implies that such a division algebra is commutative, hence each $k_j$ is a field extension of $k^{\rm sep}$, which concludes the proof. 
\end{pf}

Going back to the proof of Theorem~\ref{thmcorps}, we get that $R$ is
isomorphic to a finite product of matrix $k$-algebras, $$R \isomo \prod_{i=1}^s
M_{n_i}(k).$$
In particular, fixing such a $k$-algebra isomorphism, $\Det$ appears as a
determinant of such an algebra. By Lemma~\ref{prodlemma} (${\rm Spec}(k)$ is certainly connected), there are unique determinants
$\Det_i: M_{n_i}(k) \longrightarrow k$, say of dimension $d_i$, such that
$\Det$ is the product of the $\Det_i$ and $d=\sum_i d_i$. \ps

\begin{lemma}\label{zipliesbisautre} If $\Det : M_n(k) \longrightarrow k$ is a determinant of
dimension $d$, then $d=mn$ is divisible by $n$ and and $\Det$ is the
$m^{\rm th}$-power of the
usual determinant (here $k$ is actually any commutative ring, and $(M_n(k),\det)$ can be replaced by any Azumaya algebra equipped with its reduced norm). 
\end{lemma}

\begin{pf} Indeed, by Ziplies theorem~\cite[Thm. 3.17]{ziplies} (or Ex.~\ref{zipliesbis}), any such
determinant is a composition of the usual determinant with a multiplicative
$k$-polynomial law $k \longrightarrow k$. It is clear that any such law is of
the form $x \mapsto x^m$ for some integer $m\geq 0$.
\end{pf}

As a conclusion, we may write $d_i=m_in_i$, and if $M_i$ is the simple module of $R$ corresponding to
$M_{n_i}(k)$, then $\Det$ is the determinant of the semisimple representation
$\oplus_{i=1}^s M_i^{m_i}$. As a semisimple representation is well known to be uniquely determined by
its characteristic polynomials (Brauer-Nesbitt's theorem), this representation is unique up to
isomorphism. As $\rho$ is obviously injective, the second assertion on 
$\Ker \rho$ follows. This concludes the proof of
Theorem~\ref{thmcorps}.\ps\ps

We now investigate the case of a general field $k$, starting with the following useful observations. \ps
Let $K$ be a field extension of $k$ and denote by $k' \subset K$ the maximal separable $k$-subextension of $K$. Assume that $k'$ is finite over $k$. If $p:={\rm char}(k)>0$ assume also that there exists a finite power $q$ of $p$ such that $K^q \subset k'$. We define the {\it exponent} $(f,q) \in \N^2$ of $K/k$ by $f=[k':k]$, $q=1$ if $K=k'$, and $q$ is the smallest power of $p={\rm char}(k)>0$ as above if $K\neq k'$. \par
Let $S$ be a central simple $K$-algebra with rank $n^2$ over $K$ and reduced norm $N : S \longrightarrow K$, let $N_{k'/k} : k' \longrightarrow k$ be the usual norm (i.e. the determinant of the regular $k$-representation) and $F^q : K \longrightarrow k$ the $q^{\rm th}$-Frobenius law (see Ex.~\ref{exinsep}). Then we have a natural determinant 
$${\det}_S : S \longrightarrow k$$
of dimension $nqf$ defined by $\det_S =N_{k'/k} \circ F^q \circ N$. 

\begin{thm}\label{thmcorps2} Let $\Det : R \longrightarrow k$ be a determinant
of dimension $d$. Then as a $k$-algebra $$R/{\ker(\Det)} \isomo \prod_{i=1}^s S_i$$ 
where $S_i$ is a simple $k$-algebra which is of finite dimension $n_i^2$ over its center $k_i$, and where $k_i/k$ has a finite exponent $(f_i,q_i)$.\ps Moreover, under such an isomorphism, $\Det$ coincides with the product determinant
$$\Det=\prod_{i=1}^s
{{\det}_{S_i}}^{m_i}, \, \, \, d=\sum_i m_in_iq_if_i,$$ where $m_i$ are some uniquely determined integers.\ps

In particular, $R/\ker(\Det)$ is semisimple. It is finite dimensional over $k$ if and only if each $k_i$ is. This always occurs in each of the following three cases : $k$ is perfect, or $k$ has characteristic $p>0$ and $[k:k^p]< \infty$, or $d<p$.
\end{thm}

By Lemmas~\ref{rad2} (i) and~\ref{superlemme}, it only remains to show the following lemma.

\begin{lemma} Let $K/k$ be a field extension with finite exponent $(f,q)$ and $S$ a central simple $K$-algebra which is finite dimensional over $K$. Then any determinant $S \longrightarrow k$ has the form
$\det_S^{m}$ for some unique integer $m\geq 0$.
\end{lemma}
\begin{pf} Let $\Det : S \longrightarrow k$ be a determinant of dimension $d$ and $n^2:=\dim_K(S)$. Note that if $\Det=\det_S^m$, then we have by homogeneity  $d=fmnq$ thus $m$ is unique if it exists. Moreover, note that by Prop.~\ref{representabilite}, if two determinants $D_1, D_2: R \longrightarrow A$ of dimension $d$ are such that $D_1 \otimes_A B = D_2 \otimes_A B$ for some commutative $A$-algebra $B$ with $A \rightarrow B$ injective, then $D_1=D_2$. We will apply this below when $B$ is a field extension of a field $A$. \par 
Assume first that $k$ is separably closed (hence so is $K$); by the Noether-Jacobson theorem $S \isomo M_n(K)$ for some $n\geq 1$. Set $A:=K\otimes_k K$ and consider the kernel $I$ of natural split surjection $A \longrightarrow K$; $I$ is generated as $A$-module by the $x\otimes 1 - 1 \otimes x$, which are nilpotent of index $\leq q$, thus any finite type $A$-submodule of $I$ is nilpotent. Lemma~\ref{rad1} (iv) implies then that any determinant $M_n(A) \longrightarrow K$ factors through $\pi: M_n(A) \rightarrow M_n(A/I)=M_n(K)$. Applying this to $\Det \otimes_k K$, we get a determinant $M_n(K) \longrightarrow K$, which is an integral power of the usual determinant by Lemma~\ref{zipliesbis}, say $\Det \otimes_k K = \det^s \circ \pi$ and $d=ns$. A necessary condition is that $\det^s(M_n(K)) \subset k$, which implies that $q$ divides $s$. In particular, there is a unique possibility for $\Det \otimes_k K$, hence applying this again to $\Det'=\det_S^{s/q}$, the remark above shows that $\Det=\det_S^{s/q}$.  \par
We now reduce to the previous case. We have $$K\otimes_k k^{\rm sep} \isomo \prod_{i=1}^f K_i$$
where $K_i=K.k^{\rm sep}$ is a separable algebraic closure of $K$ such that $K_i^q \subset k^{\rm sep}$ (and $q$ is still minimal for that property), and where $\Gal(k^{\rm sep}/k)$ permutes transitively the $K_i$. Moreover, 
$$S \otimes_k k^{\rm sep} = S\otimes_K (K \otimes_k k^{\rm sep}) \isomo \prod_{i=1}^f S_i,$$
and $S_i=S\otimes_K K_i$ is central simple of rank $n^2$ over $K_i$. By Lemma~\ref{prodlemma} (iii), each $\Det \otimes_k k^{\rm sep}$ is a product of determinants $S_i\isomo M_n(K_i) \longrightarrow k^{\rm sep}$, which have the form $\det_{S_i}^{m_i}$ by the previous step and $d=n(\sum_{i=1}^f m_i)$.
As $\Det \otimes_k k^{\rm sep}$ is $\Gal(k^{\rm sep}/k)$-equivariant, this implies that $m:=m_i$ is independent of $i$, thus $m=d/nf$. In particular, there is a unique possibility for $\Det \otimes_k k^{\rm sep}$ and thus $\Det=\det_S^m$.
\end{pf}
\ps\ps
\begin{defprop}\label{defabsirr} {\rm Let $\Det : R
\longrightarrow k$ be a determinant of dimension $d$. We say that $\Det$ is
{\it absolutely irreducible} if one of the following equivalent properties is satisfied:
\ps
\begin{itemize}
\item[(i)] The unique semisimple representation $R \longrightarrow
M_d(\overline{k})$ with determinant $\Det$ (which exists by Theorem~\ref{thmcorps}) is irreducible,\ps
\item[(ii)] $(R\otimes_{k} \overline{k})/\ker(\Det\otimes_{k} \overline{k}) \simeq M_{d}(\overline{k})$,\ps
\item[(iii)] $R/\ker(\Det)$ is a central simple $k$-algebra of rank $d^2$,\ps
\item[(iv)] $R/\CH(\Det)$ is a central simple $k$-algebra of rank $d^2$,\ps
\item[(v)] for some (resp. all) subset $X \subset R$ generating $R$ as a $k$-vector space,  
there exists $x_{1},x_{2},\dots,x_{d^{2}} \in X$ such that the abstract $d^{2}\times d^{2}$ matrix $(\Lambda_{1}(x_{i}x_{j}))_{i,j}$ belongs to $\GL_{d^{2}}(k)$. 
\end{itemize}}
\ps \noindent If they are satisfied, then $\CH(\Det)=\Ker(\Det)=\{x \in R, \forall y \in R, \Lambda_1(xy)=0\}$.
\end{defprop}\ps\ps

\begin{pf} It is clear that (ii) implies (i). If $\rho : R\otimes_{k}\overline{k} \longrightarrow M_{d}(\overline{k})$ 
is as in (i), then a standard result of Wedderburn asserts that $\rho$ surjective, and we check at once that 
$\Ker(\rho) \subset \Ker(\Det\otimes_{k}\overline{k})$, hence (ii) follows by Theorem~\ref{thmcorps2}, and (v) (for any $X$) is the 
nondegeneracy of the trace on $M_{d}(\overline{k})$. If (v) holds for some $X$, we see that 
$$\dim_{\overline{k}}\left((R\otimes_{k} \overline{k})/\ker(\Det\otimes_{k} \overline{k})\right) \geq d^{2},$$ hence (v) implies (ii) by Theorem~\ref{thmcorps2}. So far, we showed that (i), (ii) and (v) are equivalent. \par
By Lemma~\ref{rad1}, the kernel of the natural surjective map $R/\CH(\Det) \longrightarrow R/\Ker(\Det)$ is a nilideal (and lemma~\ref{rad2} shows that it is the Jacobson radical), hence a standard argument shows that (iii) $\Leftrightarrow$ (iv) and that if they hold this map is an isomorphism. 
As the formation of $R/\CH(D)$ commutes with arbitrary base changes (hence with $k \rightarrow \bar{k}$), and as a $k$-algebra $E$ is central simple of rank $d^2$ if and only if $E\otimes_k \overline{k}$ has this property over $\overline{k}$, then (iv) $\Leftrightarrow$ (ii).\end{pf}


\ps\ps

Let us give some more definitions.

\begin{definition}{\rm  We say that $\Det : R \longrightarrow k$ is {\it multiplicity free} if $\Det \otimes_{k} \overline{k}$ is the
determinant of a direct sum of pairwise non-isomorphic absolutely irreducible $\overline{k}$-linear
representations of $R$. In the notations of
Theorem~\ref{thmcorps2}, it means that $m_i=q_{i}=1$ for each $i$. \par
We say that $\Det$ is {\it split} if it is the determinant of a representation $R \longrightarrow M_d(k)$. Equivalenty, $\Det$ is split if and only if $R/\ker(\Det)$ is a finite product of matrix algebras over $k$. }
\end{definition}
\ps\ps
We leave as an exercise to the reader to check the equivalences in the definition above. Moreover, we see easily that $\Det : R \longrightarrow k$ is split and absolutely 
irreducible (resp. multiplicity free) if, and only if,  $\Det$ is the determinant of a surjective $k$-representation $R \longrightarrow M_{d}(k)$ 
(resp. of a direct sum of pair-wise non isomorphic absolutely irreducible representations of $R$ defined over $k$).

\begin{example}\label{absirrloc} ({\it The absolute irreducible locus}) {\rm Let $G$ be  group (or a unital monoid) and $d \geq 1$ an integer.
If $x \in X(G,d)$, we say that $x$ is {\it absolutely irreducible} if the induced determinant
$k(x)[G] \longrightarrow k(x)$ has this property, where $k(x)$ is the residue field at $x$. Let $$X(G,d)^{{\rm irr}} \subset X(G,d)$$ be the 
subset of absolutely irreducible points. It is a Zariski open subset. Indeed, for each 
sequence of elements $\underline{g}=(g_{1},\dots,g_{{d^{2}}}) \in G^{d^{2}}$, consider the 
abstract $d^{2} \times d^{2}$ matrix $$m_{\underline{g}}=(\Tr(g_{i}g_{j})) \in M_{d^{2}}(\Z(G,d)),$$
where $\Tr=\Lambda_{1}$ is the trace of the universal determinant of $G$ of dimension $d$, and define $I \subset \Z(G,d)$ as the ideal generated by the $\det(m_{\underline{g}})$ when 
$\underline{g}$ varies in all the sequences as above. Then $X(G,d)^{{\rm irr}}=X(G,d)-V(I)$ by Def.-Prop.~\ref{defabsirr}, hence the claim.}
\end{example}

\ps\ps
\subsection{Determinants over henselian local rings} We now study the
determinants $\Det : R \longrightarrow A$ where $A$ is a local ring. We shall use the notations
of Lemma~\ref{radical}. Let $\Det : R \longrightarrow A$ be a determinant of
dimension $d$, we call $\Detb : \Rb \longrightarrow \overline{k}$ the residual determinant. 

\begin{thm}\label{structurethm} Assume that $\Det$ is Cayley-Hamilton and that $A$ is henselian.
\ps
\begin{itemize}
\item[(i)] If $\Detb$ is split and absolutely irreducible, then there is an
$A$-algebra isomorphism $$\rho: R \longrightarrow M_d(A)$$
such that $\Det = \det \circ \rho$. \ps
\item[(ii)] More generally, if $\Detb$ is split and multiplicity free, then $(R,\Tr)$ is a
generalized matrix algebra in the sense of \cite[\S 1.3]{bch}.\ps
\end{itemize}
\end{thm}

\begin{pf} The proof is almost verbatim the same as in \cite[Lemma 1.4.3]{bch}, replacing
the appeals to \cite[Lemma 1.2.5]{bch} and \cite[Lemma 1.2.7]{bch} by
the ones of \S~\ref{prellemmas}, so we will be a bit sketchy. \ps
By assumption, we have an isomorphism
$$\varphi : R/\ker(\Detb) \isomo \prod_{i=1}^s M_{n_i}(k)$$
such that $\Detb = \det \circ \varphi$ and $\sum_{i=1}^s n_i=d$. Call $\Detb_i$ the determinant of the
representation $R \longrightarrow M_{n_i}(k)$ on the $i^{\rm th}$ factor, so
$\Detb=\prod_{i=1}^s \Detb_i$. \ps
Assume first that we are in case (i), i.e. $s=1$. As $R$ is integral
over $A$, $A$ is henselian, and $\rad(R)=\ker \Detb$ by Lemma~\ref{radical},
we may find some elements $E_{i,j} \in R$ such that
$E_{i,j}E_{k,l}=\delta_{j,k}E_{i,l}$ lifting the usual basis of $M_d(k)$
(see \cite[chap. III, \S 4, exercice 5]{bki}).
Set $e_i=E_{i,i}$. Lemma~\ref{idempotents} shows that
$\Det_{e_i} : e_iR e_i \longrightarrow A$ is a Cayley-Hamilton determinant.
Its dimension is the integer $r(e_i)$ such that
$t^{r(e_i)}=\Det_{e_i}(te_i)=\Det(1-e_i+ t e_i)
\in A[t]$. Projecting this equality in $k[t]$ we get that
$$t^{r(e_i)}=\Detb(1-\overline{e_i}+t\overline{e_i})=t \in k[t]$$
so $r(e_i)=1$. By Lemma~\ref{idempotents} again, $e_1+\dots+e_s=1$ and $\Det_{e_i}
: e_iRe_i \longrightarrow A$ is Cayley-Hamilton of dimension $1$, so
$e_iRe_i=Ae_i$ is free of rank $1$ over $A$ by Lemma~\ref{CHdim1}. But if
$x \in e_iRe_j$, then $x=E_{i,j}(e_jE_{j,i}x) \in A E_{i,j}$ and we check at
once that $R=\oplus_{i,j} A E_{i,j} \simeq M_d(A)$, in which case $\Det$
necessarily coincides with the usual determinant by Ziplies'theorem. \ps
Assume now that we are in case (ii). Let us lift the family of
central orthogonal idempotents of length $s$ of $R/\Ker(\Detb)$ to a family of orthogonal idempotents
$e_1+\dots+e_s=1$ in $R$. Arguing as above we see again that 
$\Det_{e_i} : e_iR e_i \longrightarrow A$ is a Cayley-Hamilton determinant
of dimension $n_i$, which is residually split and absolutely irreducible. By (i) we
get that $e_iRe_i \simeq
M_{d_i}(A)$, and it is immediate to check that $R$ is a generalized matrix
algebra whose trace coincides with the trace of $\Det$.
\end{pf}

We get in particular the following nice corollary (see \S \ref{CHGcat}).

\begin{cor}\label{corstructurethm} Let $G$ be a group (or a unital monoid).\ps \begin{itemize}
\item[(i)] Over $X(G,d)^{\rm irr}$, the Cayley-Hamilton $\OO_X$-algebra $\cal R(G,d)$ is
an Azumaya $\OO_X$-algebra of rank $d^2$ equipped with its reduced norm. \ps
\item[(ii)] For each split $x \in X(G,r)^{\rm irr}$, the pro-artinian completion of $\OO_x$ is canonically isomorphic to the usual deformation ring of the associated absolutely irreducible representation $G \longrightarrow \GL_d(k(x))$ (see e.g. \cite{mazur}). \ps
\end{itemize}
\end{cor}

\begin{pf} Let $x \in X(G,d)^{\rm irr}$ and $A$ the strict henselianization of $\OO_x$. Recall that the formation of the Cayley-Hamilton quotient commutes with arbitrary base change. In particular, $$R(G,d)\otimes_{\Z(G,d)}A \isomo A[G]/\CH(\Det^u\otimes A).$$ Theorem~\ref{structurethm} (i) shows that the $A$-algebra on the right side is isomorphic to $M_{d^2}(A)$, thus $R(G,d) \otimes_{\Z(G,d)} \OO_x$ is an Azumaya algebra of rank $d^2$ as $\OO_x \longrightarrow A$ is faithfully flat. Part (i) follows then from the following abstract result, a variant of which is implicitely used in \cite[Thm. 5.1]{rou}\footnote{We are grateful to R. Rouquier for providing us a proof of this result.}: Let $C$ be a commutative ring, $d\geq 1$ an integer, and $R$ a $C$-algebra. Assume that for all $x \in {\rm Spec}(C)$, then $R_x$ is Azumaya of rank $d^2$ over $C_x$, then $R$ is an Azumaya $C$-algebra (locally free) of rank $d^2$. \par
   Part (ii) follows at once from Theorem~\ref{structurethm} (i), which moreover identifies canonically the universal representation to the natural map $G \longrightarrow (R(G,d)\otimes_{\Z(G,d)} \OO_x)^\ast$. 
\end{pf}

\ps\ps

\subsection{Determinants over $A[\epsilon]$}\label{deformations}

Let us fix a commutative ring $A$ and a determinant $D_0 : R \longrightarrow A$ of dimension $d$. Consider the $A$-algebra $A[\epsilon]$ with $\epsilon^2=0$ ; if $M$ is an $A$-module we will write more generally $M[\epsilon]$ for $M\otimes_A A[\epsilon]$. \ps \ps
We are interested in the set of determinants $D: R[\epsilon] \longrightarrow A[\epsilon]$ {\it lifting $D_0$}, {\it i.e.} such that $D\otimes_A A[\epsilon] = D_0$. Via the identification ${\cal M}_A^d(R,A[\epsilon]) \isomo {\cal M}_{A[\epsilon]}^d(R[\epsilon],A[\epsilon])$, it coincides with the set $\cal T$ of 
$d$-homogeneous multiplicative $A$-polynomial laws $P: R \longrightarrow A[\epsilon]$ which map to $D_0$ via the $A$-algebra homomorphism $\pi: A[\epsilon] \overset{\epsilon \mapsto 0}{\longrightarrow} A$. In other words, $\cal T=(\pi^*)^{-1}(D_0)$ where
$$\pi^*: \Hom_{A{\rm-alg}}(\Gamma_A^d(R)^{\rm ab},A[\epsilon])  \longrightarrow \Hom_{A{\rm-alg}}(\Gamma_A^d(R)^{\rm ab},A), \, \, f \mapsto \pi \circ f.$$
\noindent This expression makes $\cal T$ appear as a relative tangent space, thus $\cal T$ carries a natural structure of $A$-module in the usual way. \ps
Recall that we have a natural $A$-module isomorphism $${\cal P}_A^d(R,A[\epsilon]) \isomo {\cal P}_A^d(R,A)^2, \, \, P \mapsto (P_0,P_1), \,P=P_0+\epsilon P_1,$$
and any $P \in \cal T$ writes by definition as $$P=D_0+\epsilon \Delta$$
for some $\Delta \in {\cal P}_A^d(R,A)$. 

\begin{prop}\label{identV} The map $P \mapsto \Delta$ above induces an $A$-module isomorphism between $\cal T$ and the $A$-submodule of elements $\delta$ of ${\cal P}_A^d(R,A)$ such that for any commutative $A$-algebra $B$ and any $x, y \in R\otimes_A B$, 
$$\delta(xy)=D_0(x)\delta(y)+D_0(y)\delta(x).$$
\end{prop}

\begin{pf} Immediate from the definitions.
\end{pf}

As in the case of determinants, the polynomial map $\Delta$ (associated to some $P \in \cal T$) satisfies a number of polynomial identities. For example $\Delta(1)=0$, $\Delta(xy)=\Delta(yx)$, and $\Delta$ satisfies a variant of Amitsur's formula.\ps

In what follows, an important role will be played by the two-sided ideal $$I:=\Ker(D_0) \subset R.$$
The main reason for this are the following lemmas.  

\begin{lemma}\label{factorI} Assume that $A$ is a field of characteristic $0$ or $>d$. For any $P \in \cal T$, $I^{2N(d)} \subset \Ker(P)$. In other words, $\cal T \subset {\cal P}_A^d(\,R/I^{2N(d)}\,,A)$. 
\end{lemma}

Here, $N(d)$ is the integer coming from the Nagata-Higman theorem, in particular $N(d)\leq 2^d-1$.

\begin{pf} Let $P \in \cal T$ and $D : R[\epsilon] \longrightarrow A[\epsilon]$ the associated determinant. We check at once that via the natural injection $R \longrightarrow R[\epsilon]$, 
we have $$\Ker(P)  = R \cap \Ker(D),$$
so it suffices to show that $I^{2N(d)} \subset \Ker(D)$.\par
Remark that $I \supset \CH(
\Det)$ and consider the Cayley-Hamilton quotient $S=R[\epsilon]/\CH(\Det)$. For $r \in I[\epsilon] \subset R[\epsilon]$, we have by assumption $\Lambda_i(r) \in \epsilon A$ for all $i\geq 1$, thus 
$s^d \in \epsilon A[s]$ for all $s \in J=I[\epsilon]/\CH(\Det)$. The Nagata-Higman theorem implies that $(J/\epsilon J)^{N(d)}=0$ and then that $J^{2N(d)}=0$. 
As a consequence, $I^{2N(d)} \subset \CH(D) \subset \Ker(D)$, and we are done. 
\end{pf}

The next lemma is well-known.

\begin{lemma}\label{isurideux} There is a natural $A$-module isomorphism $$\Hom_R (I/I^2,R/I) \isomo {\rm Ext}_R^1(R/I,R/I).$$
\end{lemma}
\ps
(The $\Hom$ and ${\rm Ext}$ above are understood in the category of (left) $R$-modules.)\ps

\begin{pf} Apply $\Hom_R(-,R/I)$ to the exact sequence of $R$-modules $$ 0 \longrightarrow I \longrightarrow R \longrightarrow R/I \longrightarrow 0,$$
and use that ${\rm Ext}_R^1(R,-)=0$. 
\end{pf}

Let us study now a more specific example where those concepts apply. Assume that $A=k$ is a field, and that $S:=R/I$ is a finite dimensional semisimple $k$-algebra. Recall that by Theorem~\ref{thmcorps2}, this is always the case if $k$ is perfect, or if ${\rm char}(k)=p>0$ and $[k:k^p]<\infty$, or if $d<p$. Let $M_1, \dots, M_r$ denote the simple $S$-modules and $M:=\oplus_{i=1}^r M_r$.

\begin{prop}\label{exfinitedim} Assume either ${\rm char}(k)=0$ or ${\rm char}(k)>d$. If ${\rm Ext}^1_R(M,M)$ is finite dimensional over $k$, then so is $\cal T$.
\end{prop}

\begin{pf} For any semisimple ring $S$, the left $S$-module $S$ is a finite direct sum of simple modules, hence $\Ext^1_R(S,S)$ is a finite dimensional $k$-vector space by assumption. As a consequence, the $S$-module $\Hom_S(I/I^2,S)$ is also finite dimensional over $k$ by Lemma~\ref{isurideux}, which implies that $I/I^2$ has a finite length as $S$-module as $S$ is semisimple, hence $\dim_k I/I^2 < \infty$. But then $R/I^{2N(d)}$ is also finite dimensional over $k$ and we are done by Lemma~\ref{factorI}.
\end{pf}

Assume moreover that $R=k[G]$ and say that $k$ is perfect. Let $\rho: G \longrightarrow \GL_d(\overline{k})$ denote the unique semisimple representation of $G$ such that $\det(1-t\rho(g))=D_0(1-tg)$ forall $g \in G$ (Thm. \ref{thmcorps}), the assumption in the proposition is equivalent to $$\dim_{\overline{k}}H^1(G,{\rm ad}(\rho)) < \infty,$$
which generalizes a well-known result in the case $\rho$ is irreducible (see the remark below).\ps

\begin{remark}\label{remjoel}{\rm It would be interesting to know whether the known improvements of the Nagata-Higman theorem (as Shirshov's height theorem) lead to a generalization of this proposition to fields of characteristic $\leq d$. The arguments above actually give an explicit upper bound for $\dim_k \cal T$, which is however very bad in general\footnote{It could actually be improved by studing more carefully the successive restrictions of elements of $\cal T$ to the subspaces $I^k/I^{2N(d)}$ but still the general bound would be rather bad.}. For example, when $\rho$ is defined over $k$ (say) and irreducible, Theorem~\ref{structurethm} and a standard argument give a natural identification $\cal T \isomo H^1(G,{\rm ad}(\rho))$, which is much finer than what we got by the previous analysis. When $\rho$ is defined over $k$ and multiplicity free (and in the context of pseudocharacters), this space $\cal T$ has recently been studied by Bella\"iche~\cite{bel}. }
\end{remark}

\subsection{Continuous determinants}\label{continuite} For later use we
shall need to study a variant of the notions we have studied till now taking
care of some topology.  \ps Assume that $G$ is a topological group and that
$A$ is a topological ring.  Let $\Det : A[G] \longrightarrow A$ be a
determinant of dimension $n$, or which is the same, a homogeneous multiplicative $C$-polynomial law
$C[G] \longrightarrow A$ of degree $n$ for any subring $C\subset A$. We say that $D$ is {\it
continuous} if for each $\alpha \in I_n$, the map $\Det^{[\alpha]}:
G^n \longrightarrow A$ defined in~\S\ref{notations} is continuous. By Amitsur's formula, 
$\Det$ is continuous if, and only if, $\Lambda_i : G \longrightarrow A$ is continuous
for all $i\leq n$ (same argument as in the proof of Lemma~\ref{corimage}). 
\ps

\begin{example}(Restriction to a dense subgroup) {\rm  Assume that $H
\subset G$ is a dense subgroup, then a continuous determinant on $G$ is uniquely determined by its
restriction to $H$. Indeed, if two such determinants $\Det_1$ and $\Det_2$ coincide on
$\Z[H]$, and if $n$ denotes their common dimension, then for each $\alpha \in I_n$ the continuous maps $\Det_1^{[\alpha]},
\Det_2^{[\alpha]} : G^n \longrightarrow A$ coincide on $H^n$, hence on the
whole of $G^n$, so $\Det_1=\Det_2$.}
\end{example}

\begin{example}\label{patching} ({\it Glueing determinants})\, \,  {\rm In some 
applications to number theory, we are in the following
situation. Let $G$ be a compact topological group, $A$ and $\{A_i,i \in I\}$ topological
rings with $A$ compact, $\iota : A \longrightarrow \prod_i A_i$ a continuous
injective map, $\Det_i : A_i[G]
\longrightarrow A_i$ a continuous determinant on $G$ of dimension $d$. We assume that
for each $g$ in a dense subset $X \subset G$, $(\chi_i(g,t)) \in
A[t]$ (of course $\chi_i(g,t)$ denotes here the characteristic polynomial of
$g$ with respect to $\Det_i$). We claim that there is a continuous determinant $\Det: A[G] \longrightarrow A$
such that $\Det_i = \Det \otimes_A A_i$ for each $i$. Indeed, set $C=\prod_i A_i$
and consider the map $\psi : G \longrightarrow C[t]$, $g \mapsto
(\chi^i(g,t))$. By assumption, $\psi(X) \subset A[t]$. As $A$ is compact,
$\iota(A)$ is a closed subspace of $C$ homeomorphic to $A$, hence $\psi(G)
\subset A[t]$ for $X$ is dense in $G$ and the $\Det_i$ are continuous. The
claim follows then from Corollary~\ref{corimage} and
the dicussion above.}
\end{example}

From now on, we equip $A$, as well as all the commutative $A$-algebras $B$, with the discrete topology, and we assume that $G$ is a profinite group. In this context, a $B$-valued determinant $D$ on $G$ is continuous if, and only if, the characteristic polynomal map 
$$G \longrightarrow B[t], \, \, g \mapsto D(1+tg)$$
\noindent factors through $G \mapsto G/H$ for some normal open subgroup $H$ of $G$. \ps
This leads us to define for each normal open subgroup $H \subset G$ the two-sided ideal of $A[G]$ $$J(H):=\ker \left(A[G] \overset{\rm can}{\longrightarrow} A[G/H]\right)$$
and to equip $A[G]$ with the topology defined by this filtered set of ideals. 

\begin{lemma}\label{equivcont} A $B$-valued determinant $D$ on $G$, viewed as an element  $P \in {\cal M}_A^d(A[G],B)$, is continuous if, and only if, $\Ker(P) \subset A[G]$ is open (that is, contains some $J(H)$). \par
If it is the case, then the natural representation $$G \longrightarrow (B[G]/\Ker(D))^*$$
factors through a finite quotient $G/H$ of $G$ for some open subgroup $H$.\end{lemma}

\begin{pf} If $\Ker(P) \supset J(H)$, then $P$ factors through an element of ${\cal M}_A^d(A[G/H],B)$ hence $D$ is obviously continuous. Assume conversely that $D$ is continuous. As $B$ is discrete and $G$ profinite, there is an open normal subgroup $H \subset G$ such that all the $\Lambda_i : G \longrightarrow B$ factor through $G/H$. As a consequence, Amitsur's formula shows that for $g \in G$ and $h \in H$,
$$D(t(g-gh)+\sum_i t_i g_i)=D(\sum_i t_i g_i)$$
so $g-gh \in \Ker(P)$, and $J(H)=\sum_{g \in G, h\in H}A g(h-1) \subset \Ker(P)$. The last assertion is obvious.
\end{pf}
 \ps

\begin{example}\label{excont} {\rm Assume that $A=k$ is a field and consider the (unique) semisimple representation $$\rho : G \longrightarrow \GL_d(\overline{k})$$ 
such that $\det(1+t\rho(g))=D(1+tg)$ for all $g \in G$ (see Thm.~\ref{thmcorps}). Equip $\GL_d(\overline{k})$ with the discrete topology. Then $\rho$ is continuous if, and only if, $\Det$ is continuous.}
\end{example} 

 We end by discussing continuous deformations of a continuous determinant. We adopt the notations of \S~\ref{deformations} with $R=A[G]$ as above and with $D_0$ a continuous determinant $A[G] \longrightarrow A$ of dimension $d$. Consider the $A$-submodule $$\cal T^c \subset \cal T$$
of continuous liftings of $D_0$. This $A$-module writes $$\cal T^c = \bigcup_H \cal T^H$$
where $H$ varies in the set of all normal open subgroups of $G$ such that $\ker (D_0) \supset J(H)$, and where $\cal T^H$ is defined as the subset of liftings $P$ such that $\Ker(P) \supset J(H)$.\ps
Assume now that $A=k$ and that $k[G]/\ker(D_0)$ is finite dimensional over $k$ (see~\ref{thmcorps2}), and let $\rho : G \longrightarrow \GL_d(\overline{k})$ be the continuous representation associated to $D_0$ as in Example~\ref{excont} above. The following result is a variant of Prop.~\ref{exfinitedim}.

\begin{prop}\label{tcfini} Assume ${\rm char}(k)=0$ or $>d$. If the continuous cohomology group $H^1_c(G,{\rm ad}(\rho))$ is finite dimensional over $\overline{k}$, then $\cal T^c$ is finite dimensional over $k$.
\end{prop}

\begin{pf} It is enough to show that $\dim_k \cal T^H$ is bounded independently of the normal open subgroup $H$ such that $J(H) \subset I:=\ker(D_0)$. Fix such an $H$. By Lemma~\ref{factorI} 
$$\cal T^H \subset {\cal P}_k^d(k[G]/(I^{2N(d)}+J(H)),k),$$
so it is enough to show that $\dim_k(k[G]/(I^{2N(d)}+J(H)))$ is bounded independently of $H$. \par As $J(H) \subset I$, we have for each $n\geq 1$ a natural $k$-linear surjection 
$$ \left(I/(I^2+J(H))\right)^{\otimes_k n}  \longrightarrow (I^n +J(H))/(I^{n+1}+J(H)),$$
hence it is enough to show that $\dim_k(I/(I^2+J(H)))$ is bounded independently of $H$ (recall that $k[G]/I$ is finite dimensional). As $k[G]/I$ is a semisimple $k$-algebra, it is then enough to show that 
$$\dim_k\left(\Hom_{k[G]}(I/(I^2+J(H)),k[G]/I)\right)$$ is bounded independently of $H$. But by Lemma~\ref{isurideux}, 
$$\Hom_{k[G]}(I/(I^2+J(H)),k[G]/I) \isomo \Ext^1_{k[G/H]}(k[G]/I,k[G]/I)$$
and this latter space is naturally a subvector space of the space of continuous $G$-extensions of $k[G]/I$ by itself, which does not depend on $H$, and which is finite dimensional by assumption.
\end{pf}

\begin{remark}\label{castgzero} {\rm The proof above shows moreover that if $H^1_c(G,{\rm ad}(\rho))=0$, then $\cal T^c=0$.} \end{remark}

\subsection{A finiteness result}\label{finitenessresult}
We end this paragraph by the following important finiteness result, 
which follows from works of Donkin ~\cite{Do}, Seshadri~\cite{seshadri} and 
Vaccarino~\cite{vac0}.

\begin{prop}\label{finiteness} Assume that $R$ is finitely generated as $A$-algebra and let
$d \geq 1$ be an integer, then $\Gamma_A^d(R)^{\rm ab}$ is a finite type
$A$-algebra.
\end{prop}

\begin{pf} Let $X$ be a finite set and $\Z\{X\}$ the free ring on
$X$, and set $m=|X|$. By
assumption there is a surjective $A$-algebra homomorphism $$A\{X\}:=A \otimes_\Z
\Z\{X\} \longrightarrow R,$$ hence a surjective $A$-algebra homomorphism
$\Gamma_A^d(A\{X\})^{\rm ab} \longrightarrow \Gamma_A^d(R)^{\rm ab}$, so
we may assume that $R=A\{X\}$. As $\Gamma_A^d(A\{X\})^{\rm ab}$ is
canonically isomorphic to $A \otimes \Gamma_\Z^d(\Z\{X\})^{\rm ab}$, we
may also assume that $A=\Z$. \ps
Let $B=F_X(d)$ as in \S~\ref{polid} be the coordinate ring of
$\M_d^m$ over $\Z$ and $B^H \subset B$ the ring
of invariant elements under the componentwise conjugacy of $H:=\GL_d(\Z)$.
Recall that we have a natural ring homomorphism $$\rho^{\rm univ} :  \Z\{X\} \longrightarrow
M_d(B)$$ sending
$x \in X$ to the matrix $(x_{i,j})_{i,j}$, and that $E_X(d) \subset B$ is the subring generated
by the coefficients of the characteristic polynomials of the elements of
$\rho^{\rm univ}(\Z\{X\})$. Clearly we have $E_X(d) \subset B^H$ and a theorem of 
S. Donkin (\cite[Thm. 1 and \S 3]{Do}) shows that $E_X(d)=B^H$. As 
$\GL_d/\Z$ is reductive (and in particular reduced), and as $H \longrightarrow
{\rm PGL_n}(\C)$ has a Zariski-dense image, a general result of
Seshadri~\cite[Thm. 2]{seshadri} implies that $E_X(d)=B^H$ is a finite type
$\Z$-algebra. By the result of Vaccarino (Thm. \ref{thmvac}) recalled in
\S~\ref{polid}, $\Gamma_\Z^d(\Z\{X\}) \simeq E_X(d)$, and we are done.
\end{pf}

\begin{cor}\label{corfiniteness} Assume that $G$ is finitely generated and fix $d\geq 1$. \ps
\begin{itemize}\item[(i)] $\Z(G,d)$ is a finite type $\Z$-algebra,\ps
\item[(ii)] There exists a finite set $X \subset G$ such that for each commutative ring $A$
and any two $d$-dimensional determinants $\Det_1, \Det_2 : A[G] \longrightarrow
A$, then $\Det_1=\Det_2$ if and only if\footnote{Of course, $\chi^D(x,T)$ denotes
here the characteristic polynomial of $x$ with respect to $\Det$.}
$\chi^{D_1}(x,t)=\chi^{D_2}(x,t)$ for all $x \in X$.\ps
\end{itemize}
\end{cor}

\begin{pf} Part (i) is a special case of the proposition and part (ii)
follows from part (i) and Amitsur's formula (see Lemma~\ref{corimage}) : $E_X(d)$ is generated
as $\Z$-algebra by the coefficients of the $\chi(g,t)$ for $g \in G$.
\end{pf}

\section{The universal rigid-analytic families of pseudocharacters of a
profinite group}\label{chap3}

In this section, $G$ is a profinite group\footnote{The methods of this section could easily be extended to study more general topological groups, as locally profinite ones. }, $p$ is a fixed prime
number and $d\geq 1$ is an integer.\ps
\subsection{The deformation space of a given residual
determinant}\label{defdetthm} Let 
$k$ be a finite field of characteristic $p$ equipped with its discrete topology and $$\Dbar : k[G] \longrightarrow k$$ 
a continuous determinant of dimension $d$. Recall that by Theorem~\ref{thmcorps} (and Example~\ref{excont}), it is
equivalent to give such a determinant and (the isomorphism class of) a continuous semisimple
representation $$\rhob : G \longrightarrow \GL_d(\kb)$$
such that $\det(1+t\rhob(g)) \in k[t]$ for all $g \in G$, the relation being then $\Dbar(x)=\det(\rhob(x))$ for all $x \in k[G]$. \ps

Let $W(k)$ be the ring of Witt vectors of $k$. 
Let $\cal C$ be the category whose objects are the local $W(k)$-algebras which
are finite (as a set) and with residue field isomorphic to $k$, and whose morphisms are $W(k)$-algebra
homomorphisms. If $A \in {\rm Ob}(\cal C)$, we will denote by $m_A$ its
maximal ideal. The given map $W(k) \rightarrow A$ induces then a canonical 
isomorphism $k \isomo A/m_A$, and we shall always identify $A/m_A$ with $k$
using this isomorphism. We shall always equip such an $A$ with the discrete
topology. Moreover any arrow $A \longrightarrow A'$ is
local, {\it i.e.} sends $m_A$ into $m_{A'}$, and is continuous. We define a covariant functor
$$F: {\cal C} \longrightarrow {\rm Ens}$$
as follows. For an object $A$, define $F(A)$ as the set of continuous homogeneous multiplicative $W(k)$-polynomial laws $P: W(k)[G] \longrightarrow A$ of degree $d$ (or equivalently, of continuous $A$-valued determinant $D : A[G] \longrightarrow A$ of dimension $d$), such that $P\otimes_A k
= \Dbar$ (see \S\ref{continuite}). If $\iota : A \longrightarrow A'$ is an arrow in $\cal C$, and $P \in
F(A)$, then we check immediately that $\iota \circ P \in F(A')$, which makes $F$ a functor. \ps

Let us extend the functor $F$ a little bit. Consider more generally the category $\cal C'$ whose objets are the profinite\footnote{By this we shall always mean a directed projective limit of finite rings with surjective transition maps.} local $W(k)$-algebras $A$ with residue field $k$, and whose morphisms are the local continuous $W(k)$-algebra homomorphisms. Denote by $F'(A)$ the set of continuous homogeneous multiplicative $W(k)$-polynomial laws $P : W(k)[G] \longrightarrow A$ of degree $d$ such that $P\otimes_A k=\Dbar$ (here $A \rightarrow k$ is the natural $W(k)$-algebra morphism). As before, $F' : \cal C' \rightarrow {\rm Ens}$ is a covariant functor ; it coincides by definition with $F$ over the full subcategory $\cal C$ of $\cal C'$. It turns out that $F'$ coincides with the natural pro-extension of $F$.

\begin{lemma}\label{prononrep} If $A \isomo \projlim_i A_i$ is a projective limit in $\cal C'$, then the natural map $F'(A) \longrightarrow \projlim_i F'(A_i)$ is a bijection.
\end{lemma}

\begin{pf} If $R$ is any $A$-algebra, the functor ${\cal M}_A^d(R,-)$ from $A$-algebras to sets is representable, hence commutes with any projective limit. As a map $G \longrightarrow \projlim_i A_i$ is continuous if and only if each coordinate map $G \rightarrow A_i$ is continuous, we get the lemma. \end{pf}

\begin{prop}\label{repF} The functor $F'$ is representable. 
\end{prop}

This means that there is a profinite local $W(k)$-algebra $A(\rhob)$ with
residue field $k$, and a determinant $$\Det(\rhob) : W(k)[G]
\longrightarrow A(\rhob),$$ such that for any $A \in {\rm Ob}(\cal C)$ and $D \in
F(A)$, there is a unique continuous $W(k)$-algebra homomorphism $\varphi_D: A(\rhob)
\longrightarrow A$ such that $\Det(\rhob) \otimes_{\varphi_D} A = D$. Such a pair
$(A(\rhob),\Det(\rhob))$ is unique, if exists. 

\begin{pf} Let us show the existence. Consider the $W(k)$-algebra $$B=\Gamma_{W(k)}^d(W(k)[G])^{\rm ab}=W(k)\otimes_{\Z} \Gamma_{\Z}^d(\Z[G])^{\rm ab},$$ the universal multiplicative polynomial law $P^u : W(k)[G] \longrightarrow B$, and let $\psi : B \longrightarrow k$
be the $W(k)$-algebra homomorphism corresponding to $\Dbar$. Say 
that an ideal $I \subset B$ is {\rm open} if $I \subset \Ker(\psi)$,
$B/I$ is a finite local ring and if the induced multiplicative polynomial law $P_I: W(k)[G] \longrightarrow B/I$ obtained as the composition of $P$ with $B \rightarrow B/I$ is continuous. If $I$ and $J$ are open, then so is $I\cap J$, as $B/(I\cap J) \longrightarrow B/I \times
B/J$ is injective, and a homeomorphism onto its image (!), so those ideals
define a topology on $B$. Set
$$A(\rhob):=\projlim_{I {\rm open}} B/I$$
and consider the law $P(\rhob)=\iota \circ P : W(k)[G] \longrightarrow A(\rhob)$ where $\iota : B \longrightarrow A(\rhob)$ is the canonical map. Then $A(\rhob)$ is an object of $\cal C'$ and $$P(\rhob)=(P_I) \in F'(A(\rhob))=\projlim_I F(B/I)$$ by the previous lemma. \par 
If $A$ is an object of $\cal C$ and $P \in F(A)$, then by Prop.~\ref{representabilite} there is a unique $W(k)$-algebra homomorphism $\phi: B \longrightarrow A$ such that $P = \phi \circ P^u$ and $\phi$ mod $m_A$ is $\psi$, hence $\Ker(\phi) \subset \Ker(\psi)$. But $B/\Ker(\phi) \subset A$ is necessarily finite local, and the continuity of $P$ implies that $\Ker(\phi)$ is open, and we are done by Lemma~\ref{prononrep}. 
\end{pf}

\begin{example}\label{exresirr}{\rm  If we assume that $\rhob$ is absolutely irreducible and (say) defined over $k$, then $F$ is canonically 
isomorphic with the usual deformation functor of $\rhob$ defined by Mazur in \cite{mazur}, by Theorem~\ref{structurethm}. }
\end{example}

\begin{remark}\label{remgenaz}{\rm By construction, $A(\rhob)$ is topologically generated by the $\Lambda_i(g)$ for $g \in G$ and $i\geq 1$.  }
\end{remark}

Recall that for a profinite local $W(k)$-algebra $B$, say with maximal ideal $m$
and residue field $k$, the following properties are equivalent :
\ps\begin{itemize}
\item[(i)] there is a continuous $W(k)$-algebra surjection 
$W(k)[[t_1,\cdots,t_h]] \longrightarrow B$,\ps
\item[(ii)] $B$ is noetherian,\ps
\item[(iii)] $\dim_k(m/m^2) < \infty$,\ps
\end{itemize} \ps

As is well known, the {\it tangent space} $F(k[\epsilon])$ has a natural
structure of $k$-vector space, and actually $F(k[\epsilon])=\Hom_k(m_{A(\rhob)}/m_{A(\rhob)}^2,k)$. 
This leads us to consider the following equivalent hypotheses, that we will denote by $C(\rhob)$ : \begin{itemize}\ps\ps
\item[(a)] $\dim_k F(k[\epsilon]) < \infty$,\ps\ps
\item[(b)] $A(\rhob)$ is topologically of finite type as $W(k)$-algebra.\ps\ps\end{itemize}

As an immediate consequence of  Corollary~\ref{corfiniteness} (ii) and Example~\ref{continuite}, $C(\rhob)$ holds if $G$ is topologically of finite type. 
Following Mazur, consider the following weaker condition:

\begin{center} (F) For any open subgroup $H \subset G$, there
are only finitely many continuous group homomorphisms $H \longrightarrow
\Z/p\Z$.
\end{center}

\begin{example}\label{exemgal}{\rm (F) is satisfied if 
$G$ is the absolute Galois group of a local field of characteristic $0$ or if $G=\Gal(K_S/K)$ 
with $K$ a number field, $S$ a finite set of places of $K$ and $K_S$ a maximal algebraic extension 
of $K$ which is unramified outside $S$ (by class field theory and a result of
Hermite). The condition (F) is not satisfied when $G=(\Z/p\Z)^\N$. We leave as an exercise to the reader to check that for a given $G$, $H^1_c(G,{\rm ad}(\rhob))$ is finite 
dimensional for any continuous semisimple representation $G \longrightarrow
\GL_m(\overline{\F}_p)$ (for any $m$) if and only if  (F) holds.

}\end{example}

\begin{prop}\label{finitenesstgsp} Assume either that (F) is satisfied, or that $p>d$ and $H^1_c(G,{\rm ad}(\rhob))$ is finite dimensional.
Then $C(\rhob)$ holds. 
\end{prop}

\begin{pf} In the second case, it follows from Prop.~\ref{tcfini}. When $G$ is topologically of finite type we already explained that $C(\rhob)$ holds. When we only assume (F), we are 
reduced to this case by the following lemma. Indeed, if $F^\ast : \mathcal{C} \rightarrow {\rm Ens}$ is the determinantal deformations of $\rhob$ viewed as a representation of $G/H$, the lemma shows that 
the natural transformation $F^\ast \longrightarrow F$ is an equivalence.
\end{pf}

\begin{lemma}\label{factorisationimage} Let $A$ be a commutative, profinite, local 
$W(k)$-algebra with residue field $k$ and let $\D : A[G] \longrightarrow A$ be a 
continuous determinant deforming $\det(\rhob)$. Then $\D$ factors through $A[G] \rightarrow A[G/H]$ where 
$H \subset \ker(\rhob)$ is the  smallest closed normal subgroup such that $\ker(\rhob)/H$ is pro-$p$.
\end{lemma}

\begin{pf} We have to check that $\D(T-gh)=\D(T-g)$ for all $g \in G$ and $h \in H$. We may assume that $A$ is a finite ring. By Lemma~\ref{equivcont}, we may assume that $\D$ factors 
through a finite quotient $G'$. By Lemma~\ref{radical} and Theorem~\ref{thmcorps2}, the radical of the finite ring
$$B:=A[G']/\Ker(\D)$$ is the 
kernel of the natural extension of $\rhob : k[G] \rightarrow
M_{d}(\overline{k})$. In particular, the image of the natural
continuous group homomorphism $\ker(\rhob) \rightarrow B^{*}$ falls into the $p$-group $1+\rad(B)$, what we had to prove.
\end{pf}
\ps

Assume that $C(\rhob)$ is satisfied, and consider the affine formal scheme $$\cal X(\rhob):=\Spf(A(\rhob))$$
over $\Spf(W(k))$, as well as the rigid analytic space $$X(\rhob):={\cal X}(\rhob)[1/p]$$ attached to ${\cal X}(\rhob)$ by Berthelot. Our next aim is to describe the functors that those two spaces represent.
More generally, we will identify them as component parts of the universal formal (resp.
rigid analytic) determinant of dimension $d$.

\subsection{Formal and rigid analytic determinants}\label{formrig}
\subsubsection{The formal scheme of continuous determinants}\label{prellemnform} We refer to \cite[Ch. 0 \S 7, Ch. 1 \S 10]{ega} for the basics of topological rings and formal schemes. 
\ps
Let us consider $\Z_p$ as a topological ring, equipped with the $p$-adic topology. We denote by $\cal F$ the category whose objects are the admissible topological rings $A$ equipped with a continuous homomorphism $\Z_p \longrightarrow A$, and whose morphisms are continuous ring homomorphisms. Recall that the admissibility of $A$ means that there is a topological isomorphism $$A \isomo \plim A_\lambda$$
where the limit is taken over a directed ordered set $S$ with minimal element $0$, $A_\lambda$ is a discrete ring, and each $A_\lambda \rightarrow A_0$ is surjective with nilpotent kernel.\ps
An object $A$ is said {\it topologically of finite type over $\Z_p$}, if it is a 
quotient of the topological ring\footnote{Recall that $A\langle t \rangle$
is the $A$-subalgebra of $A[[t]]$ of power series $\sum_n a_n t^n$ with $a_n
\rightarrow 0$ (say $A$ is admissible here).} $$\Z_p[[t_1,\dots,t_s]]\langle x_1,\dots,x_r \rangle$$ (for some 
$s$ and $r$) equipped with its $I$-adic topology defined by $I=(t_1,\dots,t_s,p)$. Actually, 
we would not lose much in restricting to the full subcategory of such objects of $\cal F$ but 
it is unnecessary.

 \par

\begin{lemma}\label{finitedisc} Let $A$ be an object of $\cal F$ and let $D: A[G] \longrightarrow A$ be a continuous 
determinant. Denote by $B \subset A$ the closure of the $\Z_p$-algebra generated by the $\Lambda_i(g)$ for $g \in G$ 
and $i\geq 1$. \ps\begin{itemize}
\item[(i)] $B$ is an admissible profinite subring of $A$. In particular, it is finite product of local
$\Z_p$-algebras.\ps
\item[(ii)] Assume moreover that $\iota : A \longrightarrow A'$ is a continuous $\Z_p$-algebra
homomorphism and let $D' : A'[G] \longrightarrow A'$ be the induced determinant and $B' \subset A'$
the ring associated as above. Then $\iota$ induces a continuous surjection $B \longrightarrow B'$.
\end{itemize}
\end{lemma}

\begin{pf} Assume first that $A$ is discrete, in which case the assumption reads $p^nA=0$ for some integer $n\geq 1$. Let $P: (\Z/p^n\Z)[G] \longrightarrow A$ the continuous multiplicative polynomial law associated to 
$D$. By Lemma~\ref{equivcont}, $P$ factors through $(\Z/p^n\Z)[G/H]$ for some normal open subgroup $H \subset G$. 
But $\Gamma_{\Z/p^n\Z}^d((\Z/p^n\Z)[G/H])$ is a finite ring as $G/H$ is finite, hence so is the ring of the statement 
which is (by Amitsur's relations) the image of the natural ring homomorphism 
$$\Gamma_{\Z/p^n\Z}^d((\Z/p^n\Z)[G/H]) \longrightarrow A$$
attached to $P$ (and $D$).\ps
Consider now the general case. Write $A \isomo \plim A_\lambda$ as above and denote by  $\pi_\lambda: A \rightarrow A_\lambda$ the natural projection. Let $P: \Z_p[G] \longrightarrow A$ denote the continuous multiplicative polynomial law associated to $D$, and 
$P_\lambda = \pi_\lambda \circ P$. By the discrete case, the image 
$B_\lambda$ of $B$ in $A_\lambda$ is a finite ring,
hence $$B \isomo \plim B_\lambda$$
is a profinite admissible $\Z_p$-subalgebra. The last part of the first statement holds obviously for any profinite admissible ring $B$ : the radical of $B$ contains any ideal of definition of $B$ by admissibility, hence $B/\rad(B)$ is finite as $B$ is
profinite.\par
By Amitsur's relations, the ring $B$ is the closure of the image of the
natural map $$\Gamma_{\Z_p}^d(\Z_p[G])^{\rm ab} \longrightarrow A$$
given by $D$, so the last assertion follows. 
\end{pf}
\ps
\begin{definition}\label{defgd} {\rm We denote by $|G(d)| \subset \Spec(\Gamma_{\Z_p}^d(\Z_p[G])^{\rm
ab})$ the subset of closed points $z$ with finite residue field, that we
shall denote by $k(z)$. }\end{definition}\ps

For each $z \in |G(d)|$, there is a canonical determinant $$D_z : k(z) [G] \longrightarrow
k(z).$$
By Theorem~\ref{thmcorps2}, and Ex.~\ref{excont}, $|G(d)|$ is in bijection with the set of 
continuous semisimple representations $G \longrightarrow \GL_d({\overline \F_p})$
taken up to isomorphism and Frobenius actions on coefficients\footnote{This
means that we identify such a representation $\rho$ (whose image actually
falls into some $\GL_d(F)$ with $F$ a finite subfield) exactly with the
representations $Q ({\rm Frob}^m \circ \rho) Q^{-1}$ for any $Q \in
\GL_d(\overline \F_p)$ and $m\geq 1$.}.

\begin{definition}\label{resconst}{\rm Let $A$, $D$ and $B \subset A$ be as in the statement of Lemma~\ref{finitedisc}. 
If $B$ is local, we will say that $D$ is {\it residually constant}. \par
If it is so, the radical of the kernel of
the natural surjective ring
homomorphism $$\Gamma_{\Z_p}^d(\Z_p[G])^{\rm ab} \longrightarrow B_0$$
defines a point $z \in |G(d)|$ which is independent of the ideal of definition $I$ of $B$ chosen such that $B_0=B/I$. The field $k(z)$ is canonically isomorphic to
the residue field of $B$ and the determinant $\Dbar$ obtained by reduction of $D$
via $B \rightarrow k(z)$ coincides by definition with $D_z$ : we will say that
{\it $D$ is residually equal to}
$D_z$.}
\end{definition}
 \ps

For instance, if $\Spec(A)$ 
is connected then any $D \in E(A)$ is residually constant. In general, if $D \in E(A)$ then there is a unique finite set $|D| \subset |G(d)|$, a unique decomposition $A = \prod_{i \in |D|} A_i$ in $\cal F$, and a unique collection of determinants $D^i : A_i[G] \longrightarrow A_i$ with $D^i$ residually constant equal to $D_i$, such that $D=(D^i)_i: A[G] \longrightarrow A$. \ps 

Let us define a covariant functor $$ E : \cal F \longrightarrow {\rm Ens}$$
as follows. For an object $A$ of $\cal F$, define $E(A)$ as the set of
continuous determinants $A[G] \longrightarrow A$ of (fixed) dimension $d$,
or which is the same, of continuous homogeneous multiplicative $\Z_p$-polynomial
laws $\Z_p[G] \longrightarrow A$ of degree $d$. If $\iota : A \rightarrow A'$ is a
morphism in $\cal F$ and $P \in E(A)$ is such a law, then $E(\iota)(P):=\iota
\circ P \in E(A')$, which makes $E$ a covariant functor. For each $z \in
|G(d)|$, define $E_z(A) \subset E(A)$ as the subset of determinants which are
residually constant and equal to $D_z$. As the formation of the ring $B$ of
Lemma~\ref{finitedisc} is functorial, $E_z$ is a subfunctor of $E$. \ps

As a start, let us fix some $z \in |G(d)|$ and let $\rhob_z : G
\longrightarrow \GL_d(\overline{k(z)})$ be "the" continuous semisimple representation
such that $\det(1+t\rhob_z(g))=D_z(1+gt)$ for all $g \in G$
(see Ex.~\ref{excont}).

\begin{prop}\label{repE} Assume that $C(\rhob_z)$
holds. Then $E_z$ is representable by an object $A(z)$ of $\cal F$. This object $A(z)$ is a local ring whose residue field is canonically isomorphic to $k(z)$, moreover it is topologically of finite 
type over $\Z_p$. Actually, the $W(k(z))$-algebra $A(z)$ is
canonically topologically isomorphic to $A(\rhob_z)$ of Prop.~\ref{repF}.
\end{prop}

\begin{pf} By Lemma~\ref{finitedisc}, for any object $A$ and any $P : \Z_p[G] \longrightarrow A$
in $E(A)$, $P$ is the composite of a continuous multiplicative polynomial law
$P' : \Z_p[G] \longrightarrow B$ with $B \rightarrow A$. If $P \in E_z(A)$,
then $B$ is a $W(k(z))$-algebra in a natural way and $P'$ extends to a
continuous multiplicative polynomial law $P'': W(k(z))[G] \longrightarrow B$ which
reduces to $\D_z$, thus $P'' \in F'(B)$ where $F'$ is the functor defined in
section~\ref{defdetthm}. As a consequence, there is a unique continuous 
$W(k(z))$-algebra homomorphism $A(\rhob_z) \longrightarrow B$ corresponding
to $P''$. As $C(\rho_z)$
holds, $A(\rhob_z)$ is an object of $\cal F$ which is moreover local and topologically of finite type over $\Z_p$. Unravelling the definitions we get the result.
\end{pf}

It is then essentially formal to deal with $E$ rather than a given $E_z$. For that we need to extend $E$ and the $E_z$ to the category ${\cal FS}/\Z_p$ of all formal schemes over $\Spf(\Z_p)$. \ps
For an object $\cal X$ of ${\cal FS}/\Z_p$, let $\TE(\cal X)$ be the set of continuous determinants 
$\OO(\cal X)[G] \longrightarrow \OO(\cal X)$ of dimension $d$, which makes $$\TE : {\cal FS}/\Z_p \longrightarrow {\rm Ens}$$ a contravariant functor in the obvious way. The restriction of $\TE$ to the full subcategory of affine formal scheme coincides with $E^{\rm opp}$. In the same way,
define a subfunctor $\TE_z \subset \TE$ where $\TE_z(\cal X) \subset \TE(\cal X)$ is the subset of elements $D$ such
that for any open affine $\cal U \subset \cal X$, the image of $D$ in
$\TE(\cal U)=E(\OO(\cal U))$ belongs to $E_z(\OO(\cal U))$. If $\Spf(A)=\bigcup_i {\cal U}_i$ is an affine covering, note that an element $D \in E(A)$ belongs to $E_z(A)$ if, and only if, its image in each $D_i$ belongs to $E_z(\cal O(\cal U_i))$, by Lemma~\ref{finitedisc}. In particular, the restriction of $\TE_z$ to the full subcategory of affine formal scheme coincides with $E_z^{\rm opp}$.

\begin{cor}\label{cortot} Assume that condition (F) holds for $G$ (see~\ref{defdetthm}). Then $\TE$ (resp. $\TE_z$) is representable by 
the formal scheme $\coprod_{z\in |G(d)|} \Spf(A(z))$ (resp. by $\Spf(A(z))$).
\end{cor}

\begin{pf} By definition, if $\cal X$ is a formal scheme then the topology on $\OO(\cal X)$ is the weakest topology such that the $\OO(\cal X) \longrightarrow \OO(\cal U)$ are continuous for each open affine $\cal U$. From this we check at once as in Lemma~\ref{prononrep} that $\TE$ and $\TE_z$ are sheaves for the Zariski topology on ${\cal FS}/\Z_p$. As $\TE_z$ coincides with $E_z^{\rm opp}$ on $\cal F^{\rm opp}$, we have (by Prop.\ref{repE}) a canonical isomorphism $\TE_z \isomo \Spf(A(z))$. The assertion on $\TE$ follows then from Lemma~\ref{finitedisc} (i).
\end{pf}

\subsubsection{Rigid analytic determinants}\label{analyticfinal}
 
Let $\Aff$ be the category of affinoid $\Q_p$-algebras (\cite[Ch. 6]{BGR}). We define again an obvious covariant 
functor $$E^{\rm an} : \Aff \rightarrow {\rm Ens}$$ as follows. If $A$ is an affinoid algebra, 
$E^{\rm an}(A)$ is the set of continuous determinants $A[G] \longrightarrow A$ of dimension $d$, 
and if $\varphi : A \longrightarrow B$ is a $\Q_p$-algebra homomorphism (necessarily continuous) 
and $P : \Z_p[G] \longrightarrow A$ is in $E^{\rm an}(A)$, then we set 
$E^{\rm an}(P)=\varphi \circ P$. Remark that by Prop.\ref{Qalg}, $E^{\rm an}(A)$ also
coincides with the set of continuous pseudocharacters $G \longrightarrow A$ of dimension 
$d$. \ps

Recall that for any object $\cal A$ of $\cal F$ which is topologically of finite
type over $\Z_p$, the algebra $A:=\cal A[1/p]$ is an affinoid algebra and the map $\cal A \longrightarrow 
A$ is continuous and open. We say that $\cal A$ is a model of $A$. Any affinoid algebra admits 
at least one (and in general many) such model, as $\Q_p\langle t_1,\dots,t_n\rangle$ does. If $\cal A$ is a model of $A$ we have a natural 
map $$\iota_{\cal A}: E(\cal A) \longrightarrow E^{\rm an}(A),$$ which is moreover injective 
if $\cal A$ is torsion free over $\Z_p$. If $\cal A'$ is another model of $A$ and if we have 
a continuous ring homorphism $\cal A \longrightarrow \cal A'$, we get a natural map 
$E(\cal A) \longrightarrow E(\cal A')$, whose composite with $\iota_A'$ is $\iota_A$, so we 
get a natural injective map 
$$\iota : \underset{\rightarrow}{\lim} \,E(\cal A) \longrightarrow E^{\rm an}(A),$$
the colimit being over the (directed set of) models $\cal A$ of $A$.
\ps
If $A$ is affinoid, we denote by $A^0 \subset A$ the subset of
elements $a$ with bounded powers (i.e. such that the sequence $a^n$, $n\geq
1$ is bounded in $A$). It is an open $\Z_p$-subalgebra, such that
$A^0[1/p]=A$. When $A$ is reduced, $A^0$ is a model of $A$ (actually, the
biggest torsion free model), but not in general (think about
$A=\Q_p[\epsilon]/(\epsilon^2)$).

 \begin{lemma}\label{lemmerig} Let $A$ be an affinoid algebra and $D \in E^{\rm an}(A)$. \ps
\begin{itemize}
\item[(i)] For any $g \in G$ and any $i\geq 1$, $\Lambda_i(g)
\in A^0$.\ps  
\item[(ii)] The map $\iota$ is bijective.\ps
\item[(iii)] When $A$ is reduced, then $E(A^0) = E^{\rm an}(A)$.
\end{itemize}
\end{lemma}

\begin{pf} We first check (i). Recall that an element of an affinoid algebra $A$ has bounded 
powers if and only if its image in all the residue fields $A/m$ has norm $\leq 1$ : we may 
assume that $A$ is a finite extension of $\Q_p$. Fix $g \in G$ ; up to
replacing $A$ by a finite extension, we may assume that 
$D(t-g)=\prod_{i=1}^{d}(t-x_{i}) \in A[t]$ splits in $A$, and we have to show that each $x_{i}$ has norm $1$. 
As $D(g) \in A^{*}$, each $x_{i}$ is in $A^{*}$. By Newton's relations (or by Theorem~\ref{thmcorps}), 
$D(t-g^{n})=\prod_{i=1}^{d}(t-x_{i}^{n})$ for each $n \in \Z$. By the continuity assumption, 
$D(t-g^n)=\prod_i(t-x_i^n)$ goes to $(t-1)^d$ (in $A^d$) when $n$ tends to $0$ in $\widehat{\Z}$, and it is 
a simple exercise to conclude. \par
We check (ii), it only remains to see the surjectivity of $\iota$. Let 
$D \in E^{\rm an}(A)$ and $\cal A \subset A$ a model of $A$. Consider the compact subset 
$K=\cup_{i=1}^d \Lambda_i(G) \subset A$. As $\cal A$ is open in $A$, $K$ meets only finitely 
many of the translates of $\cal A$. In particular, there exists a finite number of elements 
$k_1,\dots,k_s \in K$ such that $$K \subset \sum_{i=1}^s (k_i + \cal A).$$ 
By part (i), those $k_i$ have bounded powers, thus  
$$\cal{A'}=\cal{A}\langle k_1,\dots,k_s\rangle \subset A$$ 
is a model of $A$ containing $K$. By Amitsur's relations, we obtain that 
$D \in {\rm Im}(\iota_{\cal A'})$, hence (ii). Part (iii) is a consequence
of (ii) and the fact that $A^0$ is the biggest model of $A$ included in $A$.
\end{pf}

For $z \in |G(d)|$ and an affinoid algebra $A$, let us define $E_z^{\rm an}(A)$
as the colimit of the $E_z(\cal A)$ with $\cal A$ a model of $A$. Equivalently, a $D \in E^{\rm an}(A)$ belongs to $E_z^{\rm an}(A)$ if and only if $D=\iota_{\cal A}(D')$ for {\it some} model $\cal A$ and some $D' \in E_z(\cal A)$. Obviously,
this defines a subfunctor $$E_z^{\rm an} : \Aff \rightarrow {\rm Ens}.$$ of
$E^{\rm an}$. Let us first give a
useful alternative description of this functor. Fix an affinoid $A$ and
consider $x \in {\rm Specmax}(A)$, $L$ its residue field (a finite
extension of $\Q_p$), $\OO_L=L^0$ its ring of integers and $k$ the residue
field of $\OO_L$. We have natural maps 
$$E^{\rm an}(A) \longrightarrow E^{\rm an}(L)=E(\OO_L) \longrightarrow E(k),$$
hence a natural {\it reduction map} 
\begin{equation}\label{redmap} \RED_x : E^{\rm an}(A) \longrightarrow |G(d)|.
\end{equation}

We check at once the following characterization of $E_z^{\rm an}$ :

\begin{lemma}\label{characex} $E_z^{\rm an}(A)=\{ D \in E^{\rm an}(A), \, \, \forall x \in
{\rm Specmax}(A), \, \, \RED_x(D)=z\}$.
\end{lemma}

\begin{pf} The inclusion $\subset$ is immediate as $E_z$ is a functor. Conversely, let $D \in 
E^{\rm an}(A)$ belong to the set on the right. By Lemma~\ref{lemmerig},
it comes from an element $D' \in E(\cal A)$ for some model $\cal A \subset A$. Consider the ring $B \subset \cal A$ associated to $D'$ as in Lemma~\ref{finitedisc}, and write it as a product of local rings $B=\prod_{i=1}^n B_i$. In particular, $A = \prod_{i=1}^n A_i$ itself is a product affinoid algebras. If $x_i$ is a closed point of ${\rm Specmax}(A_i)$, with residue field $L_i$, then the kernel of the natural continuous map $B_i \longrightarrow \OO_{L_i}/m_{\OO_{L_i}}$ corresponds to $z$ by assumption on $D'$. As the natural map $\Gamma_{\Z_p}^d(\Z_p[G])^{\rm ab} \rightarrow B/\rad(B)$ is surjective by construction, $B$ is local and $D'\in E_z(\cal A)$. 
\end{pf}

Let us denote by $\An$ the category of rigid analytic spaces over $\Q_p$ (\cite{BGR}). For any rigid space $X$, we endow the $\Q_p$-algebra $\OO(X)$ with the weakest topology such that all the $\Q_p$-algebra homomorphisms $\OO(X) \longrightarrow \OO(U)$, with $U \subset X$ open affinoid, are continuous. Of course, such an $\OO(U)$ is equipped here with its usual Banach topology ; if $X$ itself is affinoid then this weak topology on $\OO(X)$ coincides with its Banach topology. For a general $X$, we check at once that $\OO(X)$ is a complete topological $\Q_p$-algebra (it is even a Frechet space if $X$ is separable), and that the sheaf $\OO_X$ becomes a sheaf of topological $\Q_p$-algebras.    \ps

Define a contravariant functor of continuous determinants $$\TE^{\rm an} : \An \longrightarrow {\rm Ens}$$ as usual : for any rigid space $X$, let $\TE^{\rm an}(X)$ be the set of continuous determinants $\OO(X)[G] \longrightarrow \OO(X)$ of (fixed) dimension $d$. Of course, over the full subcategory of affinoids, $\TE^{\rm an}$ coincides by definition with the opposite of $E^{\rm an}$. 
\par
For $z \in |G(d)|$, define $\TE^{\rm an}_z : \An \longrightarrow {\rm Ens}$ as the following subfunctor of $\TE^{\rm an}$ : $\TE^{\rm an}_z(X)$ is the set of determinants such that for all closed points $x \in X$ (with residue field $k_x$) the induced determinant in $\TE^{\rm an}(\{x\})=E^{\rm an}(k_x)=E(\OO_{k_x})$ is residually equal to $z$. By Lemma~\ref{characex}, $\TE^{\rm an}_z$ is the opposite functor of $E_z^{\rm an}$ over the full subcategory of affinoids.

\begin{thm}\label{mainthm} Assume that condition (F) holds. The functor $\TE^{\rm an}$ (resp. $\TE^{\rm an}_z$) is representable by a rigid analytic space $X$ (resp $X_z$). It is canonically isomorphic to the generic fiber of the formal scheme $\TE$ (resp. $\TE_z$). \par 
Moreover, $X$ is the disjoint union of the $X_z$, $z \in |G(d)|$, and each $X_z$ is isomorphic to a closed subspace of some $h_z$-dimensional open unit ball $\mathbb{B}_{[0,1[}^{h_z}$, $h_z \in \N$. In particular, $X$ is a quasi-Stein space.
\end{thm}

Recall that Berthelot~\cite[0.2.6]{ber} constructed a functor ${\cal FS}'/\Z_p \longrightarrow \An$, $$\cal X \mapsto \cal X^{\rm rig},$$ extending Raynaud's one, where ${\cal FS}'/\Z_p$ is the full subcategory of ${\cal FS}/\Z_p$ whose objects are locally topologically of finite type (see also \cite[Ch. 7]{dj}). The universal property of $\cal X^{\rm rig}$ is given by (\cite[\S 7.1.7.1]{dj})
\begin{equation}\label{propunivrig} \lim_{\cal Y {\rm \,\,\,\,model \,\,of\, \,} Y} {\rm Hom}_{{\cal FS}'/\Z_p}(\cal Y,\cal X) = {\rm Hom}_{\An}(Y,\cal X^{\rm rig}),\end{equation}
where $Y$ is any affinoid. In the case we are interested in of a $\Spf(A)$ with $$A \isomo W(k)[[t_1,\dots,t_h]]/I,$$
then $\Spf(A)^{\rm rig}$ is isomorphic to the closed subspace of the open unit ball of dimension $h$ $$\Spf(A)^{\rm rig} \subset \mathbb{B}_{[0,1[}^h$$
defined by $I=0$. In particular, by Corollary~\ref{cortot}, it is enough to prove the theorem to show that $\TE^{\rm an}$ (resp. $\TE^{\rm an}_z$) represents the generic fiber of $\TE$ (resp. $\TE_z$). As those functors are sheaves for the rigid-analytic Grothendieck topology on $\An$, it is enough to check the universal property over affinoids, in other words (\ref{propunivrig}). But that follows from Lemma~\ref{lemmerig} (ii), QED. \par 

\begin{remark}\label{rempas}{\rm Of course, if we are only interested in $\TE^{\rm an}_z$ for some $z$, and if $C(\rhob_z)$ holds, then the same argument and Prop.~\ref{repE} shows that $\TE^{\rm an}_z$ is representable by $\Spf(A(z))^{\rm rig}$}.
\end{remark}

\section{Complements}\label{complements}
We keep the notations of \S~\ref{chap3}. Let us assume that condition (F) holds for $G$ and denote by $\CX$ the formal scheme $\TE=\coprod_{z \in |G(d)|} \Spf(A(z))$ and $X=\cal X^{\rig}=\TE^{\rm an}$. We shall also denote by $\cal D$ and $D$ the respective universal determinants of $G$ over $\CX$ and $X$.\par
Alternatively, we might fix some $z \in |G(d)|$ and assume only that $C(\rhob_z)$ holds, in which case all what we say below would also apply to the restricted spaces $\CX=\TE_z$ and $X=\CX^{\rig}=\TE_z^{\rm an}$. 

\subsection{Completion at a point}\label{complpoint} Let us fix some (closed) point $x \in X$, with residue field $k_x$ (a finite extension of $\Q_p$), and associated continuous determinant $D(x): k_x[G] \longrightarrow k_x$. As $X$ represents a functor, we get a natural interpretation for the completed local ring $\widehat{\OO}_x$, viewed as a $k_x$-algebra, as pro-representing the functor $F(x)$ of continuous deformations of $D(x)$ to the category local artinian $k_x$-algebras with residue field $k_x$.\ps This applies in particular when $D(x) = \det \circ \rho(x)$ is absolutely irreducible and split, in which case this functor $F(x)$ is canonically isomorphic to the usual deformation functor of $\rho(x)$ in the sense of Mazur by \ref{structurethm} (i)\footnote{Of course, when an absolutely irreducible $D(x)$ is not split, but splits over $L/k_x$, we get such an interpretation for the $L$-algebra $\widehat{\OO}_x\otimes_{k_x}L$.}. \par 

\subsection{The absolutely irreducible locus} \label{univcohCH} For the same reason as in Example~\ref{absirrloc}, the locus $$X^{\rm irr} \subset X$$
whose points $x$ parameterize the absolutely irreducible $D(x)$ is an admissible (Zariski) open subset. In particular, the subfunctor $\TE^{\rm an, irr} \subset \TE^{\rm an}$, parameterizing determinants $D \in \TE^{\rm an}(Y)$ whose evaluation at each closed point of $Y$ is absolutely irreducible, is representable by the rigid analytic space $X^{\rm irr}$.\ps

The universal Cayley-Hamilton algebra on $X$ is the sheaf $$U \mapsto R(U)=\OO(U)[G]/\CH(T(U)),$$
where $T(U) : G \rightarrow \OO(U)$ is the tautological pseudocharacter on the open affinoid $U$. It defines a sheaf on $X$ as the formation of the biggest Cayley-Hamilton quotient 
commutes with any base change. \ps

Let us now prove Proposition $G$ of the introduction\footnote{Recall that if $L$ is a field and if $R$ is an Azumaya algebra over $L$ of rank $d^{2}$, that is a central 
simple $L$-algebra of dimension $d^{2}$, and if $\rho: G \rightarrow R^{*}$ is a group homomorphism, then $\rho$ is said to be absolutely 
irreducible if $$\rho \otimes_{L} \overline{L} : G \rightarrow (R\otimes_{L} \overline{L})^{*}\simeq \GL_{d}(\overline{L})$$
is irreducible. If $\det : R \rightarrow L$ is the reduced norm of $R$, then $\det(\rho)$ is absolutely irreducible if, and only if, $\rho$ is.}. We have to show that 
$E^{\rm irr}$ is represented by $X^{\rm irr}$ and 
$\rho : G \rightarrow R^{*}_{|X^{\rm irr}}$.
First, the reduced trace of Azumaya algebras induces a natural transformation 
$E^{\rm irr} \rightarrow E^{\rm an}$ which factors by definition through $X^{\rm irr} \subset X$. To show that $E^{\rm irr} \rightarrow X^{\rm irr}$ is an isomorphism it is enough to show that for any affinoid 
$\Q_{p}$-algebra $A$, and any $T \in E^{\rm an}(A)$ such that all the evaluations $T_{x}$, for all closed points $x$, are absolutely irreducible, 
there is a unique isomorphism 
class of continuous representations $\rho : G \rightarrow B^{*}$ where $B$ is an Azumaya $A$-algebra of rank $d^{2}$,  namely : the canonical map 
$\rho^{u} : G \rightarrow (A[G]/\CH(T))^{*}$. But this follows from Theorem~\ref{structurethm} as in Corollary~\ref{corstructurethm} (i) ($\rho^{u}$ is continuous as $T$ is and the reduced trace of an Azumaya algebra is nondegenerate).


\section{An application to Galois deformations}\label{sectexemple}

Let $G$ be the Galois group of a maximal algebraic extension of $\Q$ unramified outside $\{2,\infty\}$ and consider
$$\rhob : G \longrightarrow \GL_2(\F_2)$$
the {\it trivial} representation. Our main aim here is to study the generic fiber $X(\rhob)$ of the universal deformation 
of $\det(\rhob) : \F_{2}[G] \rightarrow \F_{2}$ as in \S\ref{defdetthm}, 
and more precisely its {\it odd} locus, {\it i.e.} the closed and open subspace $$X(\rhob)^{\rm odd} \subset X(\rhob)$$ where a complex 
conjugation $c \in G$ has determinant $-1$. By class field theory, 
the (separated) abelianization $G^{\rm ab}$ of $G$ is isomorphic to $\Z_{2}^{*}$, thus condition $C(\rhob)$ is satisfied and $X(\rhob)$ makes sense.

\begin{thm}\label{caspartp2} $X(\rhob)^{\rm odd}$ is the open unit ball of dimension $3$ over $\Q_{2}$. \end{thm}

\begin{remark}\label{remp2}{\rm  By a well-known result of Tate~\cite{tate}, $\rhob$ is the unique continuous semisimple representation 
$G \longrightarrow \GL_{2}(\overline{\F}_{2})$, so $|G(2)|=\{\rhob\}$ 
and $X(\rhob)$ is actually the universal $2$-dimensional $2$-adic analytic pseudocharacter of $G$.}
\end{remark}

In order to prove Theorem~\ref{caspartp2}, we shall consider the subfunctor $$F^{{\rm odd}} \subset F$$ of the deformation functor $F$ of $\det(\rhob)$ which is
defined by the condition that the characteristic polynomial of $c$ be
$T^{2}-1$ (here and below we shall use the notations of \S\ref{defdetthm}). 
This subfunctor $F^{{\rm odd}}$ is pro-representable by the quotient of $$A(\rhob)^{{\rm
odd}}=A(\rhob)/(f,g-1),$$ 
where $\D(\rhob)(T-c)=T^{2}-fT+g \in A(\rhob)[T]$. Moreover, we check at once (following the proofs of Lemmas~\ref{lemmerig} (ii) and~\ref{finitedisc})
that the generic fiber of $F^{\rm odd}$ is $X(\rhob)^{\rm odd}$. As a consequence, it is enough to show that 
\begin{equation}\label{formsmaodd}A(\rhob)^{\rm odd} \simeq \Z_{2}[[x,y,z]].\end{equation} \ps \smallskip

We start with a tangent space computation that we explain in its natural generality. 
In the following lemma, $G$ is any profinite group, $A$ is a discrete commutative ring such that $2A=0$, 
and $\D_{0} : A[G] \longrightarrow A$ is {\it the trivial determinant of dimension $2$},  so 
$$D_{0}(T-g)=(T-1)^{2} \in A[T],\, \, \, \,  \forall g \in G.$$ 
We denote by $G^{2}$ the closed subgroup of $G$ generated by the squares of the elements of $G$. This is a normal 
subgroup containing\footnote{Indeed, $xyx^{-1}y^{-1}=(xy)^{2}(y^{-1}x^{{-1}}y)^{2}y^{-2}$.} the commutators of $G$, so $G/G^{2}$ is a profinite
$\F_{2}$-vector space. We shall be interested in the $A$-module $\cal T$ of
continuous deformations of $\D_{0}$ to $A[\epsilon]$ (see~\S\ref{deformations}). By Lemma~\ref{d2aussi}, 
any $\D \in \cal T$ may be written uniquely as a pair 
$$(2+\epsilon \tau,1+\epsilon \delta)$$ for some maps $\tau,\delta : G \longrightarrow A$.

\begin{lemma}\label{computetang} The map $\D \in \cal T  \mapsto (\tau,\delta)$ is an $A$-linear isomorphism onto the $A$-module of pairs of continuous maps
$(t,d) : G/G^{2} \longrightarrow A$ where $t(1)=0$ and $d$ is a group homomorphism.
\end{lemma}

\begin{pf} Let $\D=(2+\epsilon \tau,1+\epsilon \delta) \in \cal T$ be an $A[\epsilon]$-valued determinant.  As $2A=0$, condition (b) in 
Lemma~\ref{d2aussi} is reduced to
$$\tau(g^{-1}h)=\tau(gh), \, \, \forall g,h \in G,$$ 
or which is the same $\tau(h)=\tau(g^{2}h)$ for all $g,h \in G$. The lemma 
follows then from Lemma~\ref{d2aussi}. \end{pf}

Let us go back now to the case of the Galois group $G$, for which $G/G^{2} \simeq 
\F_{2} \times \F_{2}$. By Lemma~\ref{computetang}, and taking into account the {\it odd} condition, the tangent space 
$F^{\rm odd}(\F_{2}[\epsilon])$ is isomorphic to the $\F_{2}$-vector space of pairs of maps $(\tau,\delta)$ with 
$\tau(1)=\tau(c)=0$ and $\delta : G/G^{2} \rightarrow \F_{2}$ a group homomorphism with $\delta(c)=0$, 
so $$\dim_{\F_{2}}(F^{\rm odd}(\F_{2}[\epsilon]))=3.$$
\noindent In particular, if $m$ is the maximal ideal of $A(\rhob)^{{\rm odd}}$, then $$\dim_{\F_{2}}(m/m^{2}) \leq 4,$$ and it only 
remains to show that the Krull dimension of $A(\rhob)^{{\rm odd}}$ is at least $4$, or better that the Krull dimension 
of $A(\rhob)^{\rm odd}[\frac{1}{2}]$ is at least $3$. For that it is enough to show that for some (closed) point 
$x \in X(\rhob)^{{\rm odd}}$, the completed local ring $\widehat{\OO}_{X,x}$ has Krull dimension at least $3$. Indeed, the Krull dimension of 
a local noetherian ring does not change after completion, and $\widehat{\OO}_{X,x}$ is (canonically) isomorphic to the completion 
of $A(\rhob)^{\rm odd}[\frac{1}{2}]$ at its maximal ideal defined by $x$ (see~\cite[Lemma 7.1.9]{dj}).\ps

Consider for instance the point $x$ 
parameterizing the Galois representation $$\rho_{\Delta} : G \longrightarrow \GL_{2}(\Q_{2})$$
attached by Deligne to Ramanujan's $\Delta$ modular form. This representation is irreducible, odd, with trivial residual associated determinant (actually, any such representation would allow us to conclude below). 
By \S\ref{complpoint}, $\widehat{\OO}_{X,x}$ is the universal deformation ring (in the sense of Mazur) of $\rho_{\Delta}$. But it is a well-known 
observation of Mazur \cite{mazur} that the Krull dimension of such a deformation ring is at least 
$$\dim_{\Q_{2}} H^{1}(G,{\rm ad}(\rho_{\Delta}))-\dim_{\Q_{2}}H^{2}(G,{\rm ad}(\rho_{\Delta}))=3$$
by the global Euler characteristic formula of Tate (and as $\rho_{\Delta}$ is odd). This concludes the proof of (\ref{formsmaodd}), and of Theorem~\ref{caspartp2}.
\ps

\begin{remark}\label{rmpfthm} {\rm \begin{itemize}\item[(i)] If $g \in G$ is any element such that $g$ and $-1$ generate topologically 
$G^{\rm ab} \isomo \Z_{2}^{*}$, we actually showed that $A(\rhob)^{{\rm odd}} = \Z_{2}[[\Tr(g)-2,\Tr(cg)-2,\Det(g)-1]]$, where $\Tr$ and $\Det$ denote the universal trace and determinant.\ps
\item[(ii)] A maybe more elementary method to show the smoothness of $F^{\rm odd}$ would have been to study abstractly the relations occuring in the process of lifting determinants. 
\end{itemize}}\end{remark}

To end the proof of Theorem H of the introduction, we still have to study the other (less interesting) components $X(\rhob)^{\pm}$ over which the universal trace of $c$ is $\pm 2$. As there are continuous 
characters $\chi : G \longrightarrow \{\pm 1\}$ such that $\chi(c)=-1$, $X(\rhob)^{+}$ and 
$X(\rhob)^{-}$ are isomorphic, so we focus on $X(\rhob)^{+}$. We claim that over $X(\rhob)^{+}$, the 
universal pseudocharacter $\Tr$ factors through $G^{\rm ab}/\langle c \rangle=\Z_{2}^{*}/\{\pm 1\}$. It is enough to show that : \begin{itemize}
\ps
\item[(a)] Over the whole of $X(\rhob)$, $\Tr$ factors through the maximal pro-$2$ quotient $P$ of
$G$,\ps
\item[(b)] Over $X(\rhob)^+$, $\Tr$ factors through $G/H$ where $H$ is the closed normal subgroup of $G$ generated by $c$. \end{itemize}
\ps
Indeed, assuming (a) and (b), $\Tr$ factors through the quotient of $P$ by the image $H'$ of $H$ in $P$. But $(P/H')^{\rm ab}=\Z_{2}^{*}/\{\pm 1\}$ is monogenic, so 
$P/H' = \Z_{2}^{*}/\{\pm 1\}$ by Frattini's argument. 

Part (b) above follows from the fact that $$e:=(1-c)/2$$ is an idempotent of $\Q_{2}[G]$ such that $\Tr(e)=0$, so 
$e \in \Ker(\Tr)$ by \cite[Lemme 1.2.5 (5)]{bch}, thus $\Tr(cg)=\Tr(g)$ for all $g \in G$.

Part (a) is a consequence of  lemma~\ref{factorisationimage}  (recall that  in this lemma, $G$ is any 
profinite group, $k$ a finite field of characteristic $p$, and
$\rhob : G \rightarrow \GL_{d}(\overline k)$ is any continuous semisimple representation 
such that $\det (T-\rhob(g)) \in k[t]$ for all $g \in G$).

As a consequence, we may replace $G$ by its quotient $G'=\Z_{2}^{*}/\{\pm 1\} \simeq \Z_{2}$ to study $X(\rhob)^{+}$, which is now a trivial exercise. 
Consider the (pro-representable) subfunctor 
$$F^{+}=:\Spf(A(\rhob)^{+}) \subset F$$ of deformations of $\det(\rhob)$ {\it 
as determinants on $G'$}. Its generic fiber is $X(\rhob)^{+}$, and we claim that 
$F^{+} \simeq \Spf(\Z_{2}[[x,y]])$. Indeed, as $G'/{G'}^{2} \simeq \F_{2}$, Lemma~\ref{computetang} shows that 
$$\dim_{\F_{2}}F^{+}(\F_{2}[\epsilon])=2.$$
It remains to show that de Krull dimension of $A(\rhob)^{+}$ is at least three. Consider two copies 
$\chi_{i} : \Z_{2} \longrightarrow \Z_{2}[[T_{i}]]^{*}$, $i=1,2$, of the universal $2$-adic
character of $\Z_{2}$ ($1$ being sent to $1+T_{i}$), and set $$D:=\chi_{1}\chi_{2}, \, \, \, \, \, \Z_{2}[G] \longrightarrow \Z_{2}[[T_{1},T_{2}]].$$ This 
$2$-dimensional determinant takes its values in the subring\footnote{Note that this ring coincides with $\Z_{2}[[T_{1},T_{2}]]^{\got{S}_{2}}$.} $\Z_{2}[[x,y]]$ where $x=T_{1}+T_{2}$ and $y=T_{1}T_{2}$. The induced map 
$$ A(\rhob)^{+} \longrightarrow \Z_{2}[[x,y]]$$ 
is clearly surjective, hence an isomorphism, which concludes the proof of Theorem H. \ps\ps

\begin{remark}{\rm We showed that the universal pseudocharacter on $X(\rhob)^{\pm}$ is
everywhere absolutely reducible  : precisely, it becomes a sum of two characters over a covering of $X(\rhob)^\pm$ of degree $2$
by the $2$-dimensional open unit ball. In the same vein, 
it is easy to determine the reducible locus of $X(\rhob)^{\rm odd}$ : in terms of the coordinates $x=\Tr(g)-2$, $y=\Tr(cg)-2$ and $z=\det(g)-1$ 
(see Remark~\ref{rmpfthm}), it is given by the relation 
$$x^{2}-y^{2}=4(1-x+y+z).$$
Moreover, we could show that over $X(\rhob)^{\rm odd}$ there exists a 
continuous representation $G \rightarrow \GL_2(\OO)$ whose trace is the universal pseudocharacter
$\Tr$.
}
\end{remark}

\noindent
\address{{\sc Ga\"etan Chenevier \\ Centre de Math\'ematiques Laurent Schwartz \\ \'Ecole Polytechnique \\ 91128 Palaiseau Cedex \\ FRANCE }}

\end{document}